\documentclass[10pt]{iopart}
\usepackage{graphics,epsfig,graphicx,color}

\expandafter\let\csname equation*\endcsname\relax
\expandafter\let\csname endequation*\endcsname\relax
\usepackage{amsmath}

\usepackage{mathrsfs,amsfonts}
\usepackage{fullpage}
\usepackage{url}
\usepackage{hyperref}
\graphicspath{{./FIGS/}}

\newtheorem{rem}{Remark}[section]
\newtheorem{prop}{Proposition}[section]
\newtheorem{lem}{Lemma}[section]

\numberwithin{equation}{section}
\DeclareMathAlphabet{\itbf}{OML}{cmm}{b}{it}
 \DeclareMathAlphabet\mathbfcal{OMS}{cmsy}{b}{n}

\renewcommand{\hat}{\widehat}
\renewcommand{\tilde}{\widetilde}
\def\RR{\mathbb{R}}
\def\bx{{{\itbf x}}}

\def\by{{{\itbf y}}}
\def\bz{{{\itbf z}}}
\def\bu{{{\itbf u}}}

\def\bde{\boldsymbol{\delta}}

\def\be{{\itbf e}}

\def\bU{{\itbf U}}
\def\bv{{\itbf v}}
\def\ss{{(s)}}
\def\bka{\boldsymbol{\kappa}}

\def\bV{{\itbf V}}
\def\bR{{\itbf R}}

\def\bS{{\itbf S}}
\def\bQ{{\itbf Q}}
\def\bM{{\itbf M}}
\def\bW{{\itbf W}}
\def\bLa{\boldsymbol{\Lambda}}

\def\bD{{\itbf D}}
\def\cP{{\mathcal P}}
\def\cI{{\mathcal I}}
\def\cA{{\mathcal A}}

\def\cG{{\mathcal G}}

\def\RM{{\scalebox{0.5}[0.4]{ROM}}}

\def\TR{{\scalebox{0.5}[0.4]{TR}}}
\def\cPR{{\boldsymbol{\cP}^{\RM}}}
\def\BP{{\scalebox{0.5}[0.4]{BP}}}
\def\PS{{\scalebox{0.5}[0.4]{PS}}}
\def\RTM{{\scalebox{0.5}[0.4]{RTM}}}

\def\bsig{{\boldsymbol{\sigma}}}
\def\lb{\left <}
\def\rb{\right >}
\def\cF{\mathcal{F}}

\def\om{\omega}
\def\la{\lambda}
\def\12{{\frac{1}{2}}}

\begin{document}

\title{Reduced order model approach for imaging with waves}

\author{Liliana Borcea$^1$, Josselin Garnier$^2$,
Alexander V. Mamonov$^3$, and J\"{o}rn Zimmerling$^1$}

\hyphenation{mamonov}
\hyphenation{mzaslavsky}

\address{$^1$Department of Mathematics, University of Michigan,
Ann Arbor, MI 48109-1043}
\address{$^2$Centre
  de Math\'ematiques Appliqu\'ees, Ecole Polytechnique, 91128
  Palaiseau Cedex, France} 
\address{$^3$Department of Mathematics, University of Houston,
  3551 Cullen Blvd Houston, TX 77204-3008 }
	\eads{\mailto{borcea@umich.edu}, \mailto{josselin.garnier@polytechnique.edu}, \mailto{mamonov@math.uh.edu}, \mailto{jzimmerl@umich.edu}}

%%===============================================================================
\begin{abstract}
We introduce a novel, computationally inexpensive approach for imaging with 
 an active array of sensors, which probe an unknown medium with a  pulse and 
 measure the resulting waves. The imaging function uses a data driven 
estimate of the ``internal wave" originating from the vicinity of the imaging point and propagating to the 
sensors through the unknown medium. We explain how this estimate can be obtained
using a reduced order model (ROM) for the wave propagation. We analyze the imaging function,
connect it to the time reversal process  and describe how its resolution depends on the 
aperture of the array, the bandwidth of the probing pulse and the medium through which the waves propagate. We also show how the internal wave can be used for selective focusing of waves at points in the imaging region. This 
can be implemented experimentally and can be used for pixel scanning imaging. We assess the performance of the imaging methods with  numerical simulations and compare them to the conventional reverse-time migration method and 
the  ``backprojection'' method introduced recently as an application of the same ROM. \\

\end{abstract}


\submitto{\IP}

%\maketitle

%===============================================================================
\section{Introduction}
\label{sec:intro}
This paper is concerned with  an application of reduced order modeling to imaging reflective structures in a known, non-scattering  host medium, from data gathered by an
active array of $m$ sensors that emit probing pulses and measure the resulting waves.

Model reduction  is an important topic in computational science, 
which traditionally has been concerned with  
finding a low-dimensional reduced order model (ROM)
that approximates the response (observables) of a given dynamical system for a set of inputs \cite{antoulas2001survey,benner2015survey}.
In wave-based imaging the inputs and observables are controlled and measured by the sensors in the array, but the dynamical system is  not given, as it is governed by the wave equation with unknown coefficients like the wave speed. Thus,  we need a data-driven ROM. 

Data driven reduced order modeling is a growing field which combines ideas from traditional model reduction and learning \cite{brunton2019data}. Much of it is concerned with studying dynamical systems using the Koopman operator theory \cite[Chapter 1]{mauroy2019koopman}. Dynamical system identification based on the Koopman operator has been proposed for instance in \cite{brunton2016discovering} and  \cite[Chapter 13]{mauroy2019koopman}. However, these approaches are difficult to use in imaging because they assume knowing the full state of
the dynamical system, aka the snapshot of the wave, at a finite set of time instances. Since we only know the wave at the sensor locations, which are far from the imaging region, a new way of 
learning the wave propagation from the data is needed.

To our knowledge, the first sensor array data driven ROM for wave propagation was introduced in
\cite{druskin2016direct} for the one-dimensional wave equation. The extension to higher dimensions was obtained in  \cite{borcea2019robust} and was analyzed in \cite{borcea2020reduced}. The latter study showed that wave propagation can 
be viewed as a discrete  time dynamical system governed by a ``propagator operator", where the time step $\tau$ is the data sampling interval. This propagator operator maps the 
wave from the states at instants $(j-1)\tau$ and $j \tau$ to the future state at time $(j+1)\tau$, 
for any $j \in \mathbb{N}$. The ROM in \cite{borcea2019robust,borcea2020reduced} is an algebraic analogue of the  dynamical system. Its evolution is controlled by  an $nm\times nm$ propagator matrix, given by the Galerkin 
projection of the propagator operator on the function space spanned by the first $n$ snapshots of the 
wave, assuming that the array records for the duration $(2n-1) \tau$. What distinguishes the ROM construction from the many other Galerkin projections  in the literature, is that it is obtained only from the measurements of the
snapshots at the sensors in the array i.e., it does not require knowing the approximation space.

The ROM  introduced in  \cite{borcea2019robust,borcea2020reduced} has been used so far to: (1) 
Approximate the Fr\'{e}chet derivative of the reflectivity to sensor array data map \cite{DtB}. This gives the single scattering (Born) forward map used in conventional  imaging 
in radar \cite{curlander1991synthetic,cheney2009fundamentals}, seismic inversion \cite{biondi20063d} and elsewhere. (2) Obtain a fast converging, iterative inverse scattering method for the acoustic impedance  in a medium with known and smooth wave speed \cite{borcea2020reduced}. (3) Develop a non-iterative  ``backprojection" imaging method  that is free of multiple scattering artifacts \cite{druskin2018nonlinear}.

We propose yet another application of the ROM: Estimate the ``internal" acoustic wave that originates from the vicinity of the imaging point and propagates through the unknown medium to the array of sensors. This idea has been tried before for Schr\"{o}dinger's equation in the spectral (frequency) domain
\cite{borceainternal,druskin2021lippmann}, where the solutions are smooth functions that are easier to approximate than the internal waves  in this paper. We use the internal waves  for two novel imaging methods: The first is a computationally inexpensive approach designed to sense rapid changes of the wave speed in the vicinity of the imaging point. Its imaging function
is connected to the point spread function of the time reversal process, and we explain how the aperture of the array, the bandwidth of the probing pulse  and the medium through which the waves propagate affect the resolution. The second method can be implemented experimentally. It controls the excitation from the array in order to focus waves at the imaging points, and then uses a matched field approach to image with the resulting backscattered wave  in a pixel scanning manner.  

The paper is organized as follows: We begin in section \ref{sect:setup} with the mathematical formulation of the imaging problem and review briefly from \cite{borcea2020reduced} the relevant 
facts about the ROM, needed in the next sections. The estimation of the internal wave is described 
in section \ref{sect:intern}. The first imaging method based on this internal wave is introduced and analyzed in section \ref{sect:normg}. We also give there a comparison with the backprojection 
imaging method introduced in \cite{druskin2018nonlinear}. The second, pixel scanning imaging method is described in section 
\ref{sect:pixel}. We use numerical simulations in section \ref{sect:num} to assess the performance of the imaging methods and to compare them with the backprojection approach \cite{druskin2018nonlinear} and with the conventional, reversed time migration method 
\cite{biondi20063d}. We end with a summary in section \ref{sect:sum}.

\section{Formulation of the imaging problem and the ROM}
\label{sect:setup}

We are interested in imaging reflective structures in a non-scattering and known host medium occupying the bounded domain $\Omega$, using data gathered by an active array of sensors located at $\bx_s$, for $s =1, \ldots, m$.  We suppose  that the aperture of the array is planar  in three-dimensions or linear 
in two-dimensions, and call ``range"  the spatial coordinate in the direction orthogonal to it. The coordinates in the plane (line in two-dimensions) parallel to the aperture are called ``cross-range".

The $s^{\rm th}$ sensor probes the medium with a pulse $f(t)$ and thus generates the wave $w^\ss(t,\bx)$, the solution of the wave equation
\begin{equation}
\label{eq:F1}
\partial_t^2 w^\ss(t,\bx) + A(c) w^\ss(t,\bx) = f'(t) \delta_{\bx_s}(\bx), \qquad t \in \RR, ~\bx \in \Omega,
\end{equation}
with quiescent initial condition
\begin{equation}
w^\ss(t,\bx) \equiv 0, \qquad t \ll 0, ~~ \bx \in \Omega.
\label{eq:F2}
\end{equation}
We assume henceforth a real valued  pulse $f(t)$, that is an even function supported in the short interval $(-t_f,t_f)$ and has non-negative Fourier transform\footnote{If the source emits an arbitrary pulse $\mathfrak{f}(t)$, we can convolve the echos received at the array with $\mathfrak{f}(-t)$. Mathematically, this  is equivalent to having the even pulse $f(t) =  \mathfrak{f}(t) \star_t \mathfrak{f}(-t)$, with Fourier transform $\hat f(\om) = 
\big| \hat{\mathfrak{f}}(\om) \big|^2 \ge 0$.
} 
\begin{equation}
\hat f(\om) = \int_{\RR} dt \, f(t) e^{i \om t} = \int_{\RR} dt \, f(t) \cos(\om t) \ge 0,
\label{eq:Ff}
\end{equation} 
which is negligible in the complement of the set $(\om_c-B,\om_c+B) \cup (-\om_c-B,-\om_c+B)$,
where $\om_c$ is the center (carrier) frequency and $B = O(1/t_f)$ is the bandwidth. 

The reflective structures are modeled in \eqref{eq:F1} by rough changes (jumps) of the wave speed $c(\bx)$, with respect to  the known and smooth reference wave speed $c_o(\bx)$ of the host medium. These changes are supported in the imaging domain $\Omega_{\rm im}$, which is a subset of $\Omega$ lying at large distance from the array. The unknown 
$c(\bx)$ appears as a coefficient in the self-adjoint, second order elliptic operator
\begin{equation}
A(c) = -c(\bx) \Delta \big[ c(\bx) \cdot \big],
\label{eq:F3}
\end{equation}
with homogeneous boundary conditions at $\partial \Omega$. The self-adjointness  of $A(c)$
is convenient for the operator calculus in \cite{borcea2020reduced} and the next sections, but we note that  \eqref{eq:F1} can be written in the standard wave equation form for the acoustic pressure 
$
p^\ss(t,\bx) = c(\bx) w^\ss(t,\bx).
$
Since $c(\bx)$ is known and equal to $c_o(\bx)$ at the sensor locations, the measurements of the pressure, 
sampled in time at interval $\tau$, define  the array data 
\begin{equation}
\mbox{data} = \left\{w^\ss(t,\bx_r), ~s,r=1,\ldots, m, ~~ t = j \tau, ~~ j =0,\ldots, 2n-1\right\}.
\label{eq:F4}
\end{equation}

\begin{figure}[t!]
\centering
\includegraphics[width=0.38\textwidth]{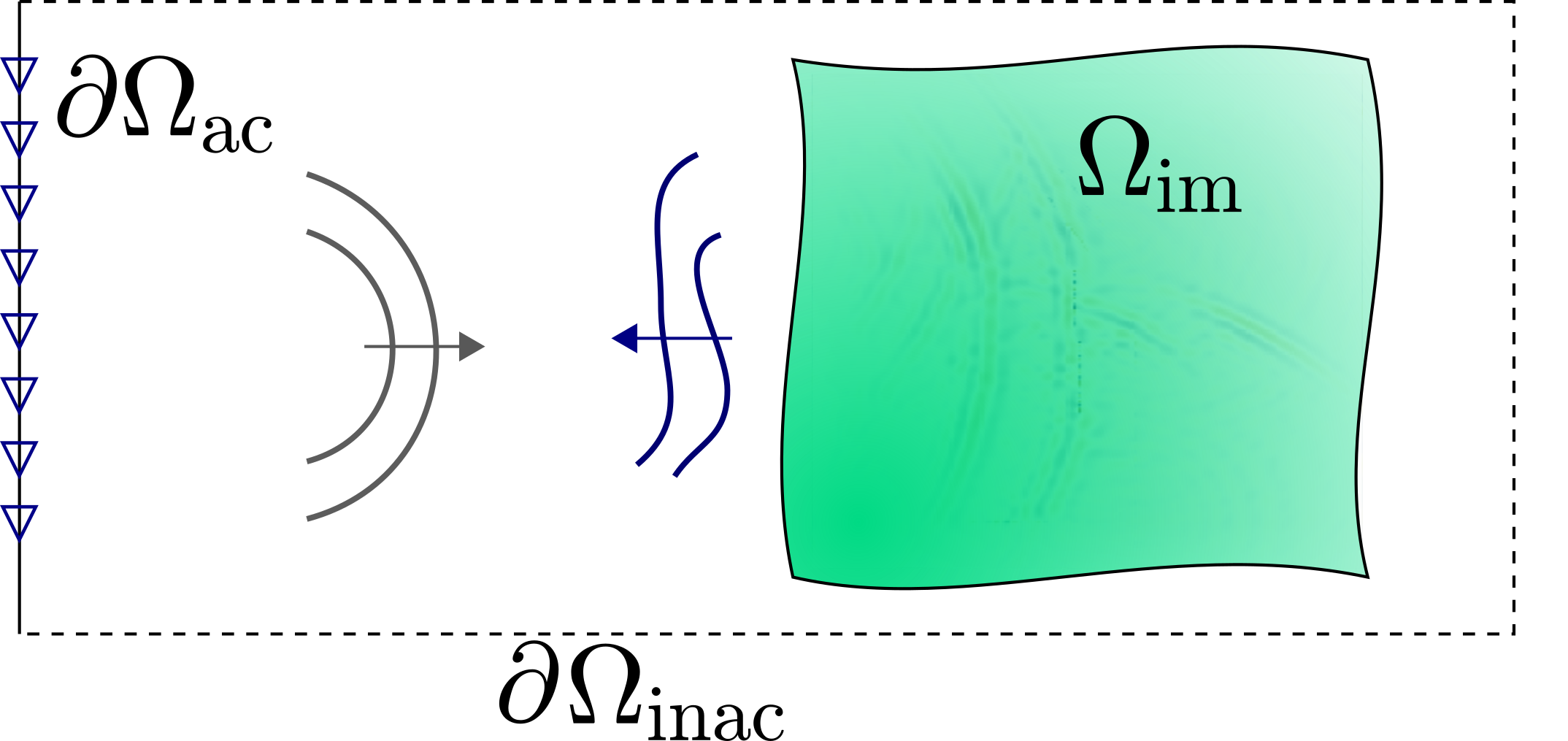}
\vspace{-0.1in}\caption{Illustration of the setup: An array of sensors (indicated with triangles) lying near the accessible boundary
$\partial \Omega_{\rm ac}$ probes a medium with incident waves and measures the backscattered waves.  The inaccessible
boundary  $\partial \Omega_{\rm inac}$ is  drawn with the dashed line. The sought after reflectors  are supported in the remote subdomain $\Omega_{\rm im}$. }
\label{fig:setup}
\end{figure}

The domain $\Omega$ may be physical or the truncation of an infinite domain, justified by 
 %causality 
{hyperbolicity}
 and the finite duration $(2n-1) \tau$ of the measurements. In either case, we divide the boundary in two parts: The ``accessible" boundary $\partial \Omega_{\rm ac}$, named so because it lies 
 in the immediate vicinity of the array,  and the ``inaccessible" boundary $\partial \Omega_{\rm inac} = 
 \partial \Omega \setminus \partial \Omega_{\rm ac}$. The accessible boundary is useful for the ROM construction because  the waves propagate only on one side of the array\footnote{If such a boundary does not exist, the medium should be known and homogeneous on the other side of the array, so that 
 the waves there can be removed with some additional processing.}, as illustrated in Fig. \ref{fig:setup}. We model it as sound hard, using the homogeneous Neumann boundary condition 
 \begin{equation}
 \partial_n w^\ss(t,\bx) = 0, \qquad t \in \RR, ~~ \bx \in \partial \Omega_{\rm ac},
\label{eq:F5}
\end{equation}
where $\partial_n$ denotes the normal derivative. 
The inaccessible boundary is modeled as sound soft,
\begin{equation}
w^\ss(t,\bx) = 0, \qquad t \in \RR, ~~ \bx \in \partial \Omega_{\rm inac},
\label{eq:F6}
\end{equation}
and if it is due to the truncation of an infinite domain, it is sufficiently far away from the sensors to affect the waves over the duration of the measurements.

The imaging problem is to estimate the support of the large and localized variations  $c(\bx)-c_o(\bx)$ of the wave speed, from the data \eqref{eq:F4} collected by the array. 

\subsection{Review of the ROM for wave propagation}
\label{sect:setup.1}
Here we review briefly from \cite{borcea2020reduced} the relevant facts 
about the ROM, needed to state the new results. 

\subsubsection{The dynamical system for wave propagation.}
We work with the even  in time wave
\begin{equation}
w_e^\ss(t,\bx) := w^\ss(t,\bx) + w^\ss(-t,\bx),
\label{eq:F7}
\end{equation}
which satisfies $w_e^\ss(t,\bx) = w^\ss(t,\bx)$  for $ t > t_f$, due to causality and the initial condition \eqref{eq:F2}. During the short duration $t_f$ of the pulse, the wave senses only the vicinity 
of the sensor location $\bx_s$, where the wave speed equals the known $c_o(\bx)$.  
Thus, the second term in \eqref{eq:F7} can be calculated and we can work with the data
matrices
\begin{equation}
\bD_j = \big( D_j^{(s,r)} \big)_{s,r=1,\ldots, m}, \quad D_j^{(s,r)} := w_e^\ss(j \tau,\bx_r), ~~ j=0, \ldots, 2n-1.
\label{eq:F9}
\end{equation}

Because  $f(t)$ is even, it is easy to obtain from \eqref{eq:F1} and \eqref{eq:F5}--\eqref{eq:F6} that $w_e^\ss(t,\bx) $ satisfies
\begin{align}
\partial_t^2 w_e^\ss(t,\bx) + A(c) w_e^\ss(t,\bx) &= 0, \qquad t>0, ~\bx \in \Omega, \label{eq:F10} \\
 \partial_n w_e^\ss(t,\bx) &= 0, \qquad t >0, ~~ \bx \in \partial \Omega_{\rm ac},  \label{eq:F11} \\
w_e^\ss(t,\bx) &= 0, \qquad t >0, ~~ \bx \in \partial \Omega_{\rm inac},  \label{eq:F12}
\end{align}
with initial conditions derived in \cite[Appendix A]{borcea2020reduced}
\begin{equation}
w_e^\ss(0,\bx) = \hat f\big(\sqrt{A(c)}\big)\delta_{\bx_s}(\bx), \qquad 
\partial_tw_e^\ss(0,\bx) = 0, \qquad \bx \in \Omega.
\label{eq:F13}
\end{equation}
We define throughout functions of the operator $A(c)$ in the standard way, using its spectral decomposition deduced from  \cite[Theorem 4.12]{mclean}. Specifically, if we  denote
by $\{\theta_l >0, ~l \geq1\}$ the eigenvalues, ordered like  $0 <\theta_1 \le \theta_2 \le \ldots$ and 
satisfying $\lim_{l \to \infty} \theta_l = \infty$,  and by $\{y_l(\bx), ~l \geq 1\}$ the eigenfunctions, 
which form an orthonormal basis of $L^2(\Omega)$
with the appropriate boundary conditions, 
then we have
\begin{equation}
\hat f\big(\sqrt{A(c)}\big)\delta_{\bx_s}(\bx) := \sum_{l=1}^\infty 
\hat f\big(\sqrt{\theta_l}\big) y_l(\bx) y_l(\bx_s).
\label{eq:F14}
\end{equation}
The pulse is band-limited, so the sum is for   $l \le l_{\max}$, where $\sqrt{\theta_{l_{\max}+1}} > \om_c+B.$ 

Note that the solution of \eqref{eq:F10}--\eqref{eq:F13} is
\begin{align}
w_e^\ss(t,\bx) &= \cos \big(t\sqrt{A(c)}\big)\hat f\big(\sqrt{A(c)}\big)\delta_{\bx_s}(\bx) 
= \sum_{l=1}^\infty 
\cos\big(t\sqrt{\theta_l}\big)\hat f\big(\sqrt{\theta_l}\big) y_l(\bx) y_l(\bx_s)\nonumber \\
&=\hat f^\12\big(\sqrt{A(c)}\big)\cos \big(t\sqrt{A(c)}\big)\hat f^\12\big(\sqrt{A(c)}\big)\delta_{\bx_s}(\bx), 
\label{eq:F15}
\end{align}
and that the data matrices \eqref{eq:F9} can be written in symmetric inner product form as follows
\begin{align}
D_j^{(r,s)} &= \int_\Omega \delta_{\bx_r}(\bx)  \cos \big(j \tau\sqrt{A(c)}\big)\hat f\big(\sqrt{A(c)}\big)\delta_{\bx_s}(\bx)  = \lb  \delta_{\bx_r}^f, \cos \big(j \tau\sqrt{A(c)}\big)\delta^f_{\bx_s}\rb.
\label{eq:F16}
\end{align}
Here we used that functions of $A(c)$ commute, and denoted by
\[
\lb \phi,\psi \rb := \int_{\Omega} d \bx \, \phi(\bx) \psi(\bx), \qquad \forall \phi, \psi \in L^2(\Omega),
\]
the $L^2(\Omega)$ inner product. We also  introduced the ``sensor functions"
\begin{equation}
\delta^f_{\bx_s}(\bx) := \hat f^\frac{1}{2}\big(\sqrt{A(c)}\big)\delta_{\bx_s}(\bx), \qquad s = 1,\ldots,m,
 \label{eq:F18}
\end{equation}
and used the assumption \eqref{eq:Ff} to define the square root of $\hat f$.

The notation in \eqref{eq:F18} reminds us that  $\delta^f_{\bx_s}(\bx)$ is a pulse dependent, blurry version of the  Dirac $\delta_{\bx_s}(\bx)$. Indeed,  comparing \eqref{eq:F18} with \eqref{eq:F13}, we note 
that $\delta^f_{\bx_s}(\bx)/2$ is the initial state of the solution of \eqref{eq:F1}--\eqref{eq:F2} and \eqref{eq:F5}--\eqref{eq:F6}, when the sensor emits the pulse
\begin{equation}
\check f^{\frac{1}{2}}(t) := \int_{\RR} \frac{d \om}{2 \pi} \hat f^{\frac{1}{2}}(\om) e^{-i \om t} =
\int_{\RR} \frac{d \om}{2 \pi} \hat f^{\frac{1}{2}}(\om) \cos(\om t),
\label{eq:F19}
\end{equation}
which is real valued, even and satisfies $f(t) = \check f^{\frac{1}{2}}(t) \star_t \check f^{\frac{1}{2}}(t).$
Because of causality and the finite wave speed, $\delta^f_{\bx_s}(\bx)$ is supported in 
a ball centered at $\bx_s$, with radius of order $c_o(\bx_s) t_f \ll \mbox{dist}\big(\bx_s,
\Omega_{\rm im}\big)$ and it  can 
be computed using the operator $A(c_o)$ in the host medium
\begin{equation}
\delta^f_{\bx_s}(\bx) = \hat f^\frac{1}{2}\big(\sqrt{A(c)}\big)\delta_{\bx_s}(\bx) = \hat f^\frac{1}{2}\big(\sqrt{A(c_o)}\big)\delta_{\bx_s}(\bx) .
 \label{eq:F20}
\end{equation}

Let us group all the sensor functions in the $m-$dimensional row vector field
\begin{equation}
\bde^f(\bx) = \left(\delta^f_{\bx_1}(\bx),\ldots,\delta^f_{\bx_m}(\bx) \right),
\label{eq:F21}
\end{equation}
and define the ``snapshots" as the $m-$dimensional row vector fields
\begin{equation}
\bu_j(\bx) = \left( u_j^{(1)}(\bx), \ldots, u_j^{(m)}(\bx) \right):= \cos \big( j \tau \sqrt{A(c)}\big) \bde^f(\bx), \qquad j \ge 0.
\label{eq:F22}
\end{equation}
These are the states of the discrete time dynamical system governed by the 
``propagator operator"
\begin{equation}
\cP := \cos \big(\tau \sqrt{A(c)}\big).
\label{eq:F23}
\end{equation}
Indeed, using a trigonometric identity of the cosine, we get  that the states evolve like
\begin{align}
\bu_{j+1}(\bx) &= 2\cP \bu_{j}(\bx) -\bu_{j-1}(\bx), \qquad j \ge 0, ~~ \bx \in \Omega, \label{eq:F24} 
\end{align}
starting from 
\begin{align}
 \bu_{0}(\bx) =  \bde^f(\bx), \qquad  
%\bu_{1}(\bx) =
\bu_{-1}(\bx) = \cP\bde^f(\bx), 
\qquad  \bx \in \Omega. \label{eq:F25} 
 \end{align}

\subsubsection{Data driven ROM construction.} 
\label{sect:ROMrev}
The ROM is the algebraic analogue of the  dynamical system 
\eqref{eq:F24}, 
\begin{align}
\bu_{j+1}^\RM &= 2\cPR \bu_{j}^\RM -\bu^\RM_{j-1}, \qquad j \ge 0,  \label{eq:F26} 
\end{align}
with propagator matrix $\cPR \in \RR^{nm \times nm}$ and states $\bu_j^{\RM} \in \RR^{nm \times m}$.
It corresponds to the Galerkin projection of \eqref{eq:F24} on the $nm-$dimensional function space spanned by
the first $n$ snapshots \eqref{eq:F22}. Using linear algebra notation, we write this space as
\begin{equation}
\mathscr{S}:=\mbox{range} \, \bU(\bx), \qquad \bU(\bx) = \big(\bu_0(\bx), \ldots, \bu_{n-1}(\bx)\big),
\label{eq:F27}
\end{equation}
and note that $\bU(\bx)$ is a $nm-$dimensional row vector field and that it is unknown. Nevertheless, the construction in \cite[Section 2]{borcea2020reduced}
shows that it is possible to get the ROM \eqref{eq:F26} from what we know: the initial snapshot
$\bde^f(\bx)$ in \eqref{eq:F21} and the $m\times m$ data matrices
\begin{equation}
\bD_j = \int_{\Omega} d \bx \, \bde^{f}(\bx)^T \bu_j(\bx) =: \langle\langle \bde^f, \bu_j \rangle \rangle, \qquad
j=0,\ldots,2n-1,
\label{eq:F28}
\end{equation}
where $T$ denotes the transpose and we introduced the notation $\langle \langle \cdot, \cdot \rangle \rangle$ for the integral of the outer product of $m-$dimensional row vector functions.

Key to the ROM construction are the data driven, symmetric, positive definite ``mass" matrix $\bM:= \int_\Omega d\bx  \bU(\bx) ^T \bU (\bx)  \in \RR^{nm \times nm}$ with $m \times m$ blocks
\begin{align}
\bM_{j,l} &:= \langle \langle \bu_j,\bu_l \rangle \rangle \nonumber \\ &= \langle \langle \cos \big(j \tau \sqrt{A(c)} \big) 
\bde^f,\cos \big(l \tau \sqrt{A(c)} \big)\bde^f \rangle \rangle \nonumber \\
&= \langle \langle  
\bde^f,\cos \big(j \tau \sqrt{A(c)} \big)\cos \big(l \tau \sqrt{A(c)} \big)\bde^f \rangle \rangle \nonumber \\
&= \frac{1}{2}\left[  \langle \langle  
\bde^f,\cos \big((j +l)\tau \sqrt{A(c)} \big)\bde^f \rangle \rangle + \langle\langle  
\bde^f,\cos \big(|j -l|\tau \sqrt{A(c)} \big)\bde^f \rangle \rangle \right] \nonumber \\
&= \frac{1}{2} \left(\bD_{j+l} +\bD_{|j-l|} \right), \qquad j,l =0, \ldots, n-1,
\label{eq:F29}
\end{align}
and the ``stiffness" matrix $\bS:= \int_\Omega d\bx \bU(\bx) ^T \cP \bU(\bx)  \in \RR^{nm \times nm}$ with $m \times m$ blocks
\begin{align}
\bS_{j,l} &:= \langle \langle \bu_j,\cP \bu_l \rangle \rangle \nonumber \\ 
&= \frac{1}{2} \langle \langle \bu_j, \bu_{l+1}+\bu_{l-1} \rangle \rangle \nonumber \\ 
&= \frac{1}{4} \left(\bD_{j+l+1} +\bD_{|j-l-1|} + \bD_{|j+l-1|} +\bD_{|j-l+1|} \right), \qquad j,l =0, \ldots, n-1.
\label{eq:F30}
\end{align}
Here we used definitions \eqref{eq:F22} and \eqref{eq:F28}, the self-adjointness of $A(c)$ and a trigonometric identity for the cosine.

The ROM propagator is defined in \cite[Section 2.2.1]{borcea2020reduced} as follows: Let $\bR \in \RR^{nm \times nm}$ be the block upper triangular matrix obtained from the block 
Cholesky factorization of the mass matrix
\begin{equation}
\bM = \bR^T \bR.
\label{eq:F31}
\end{equation}
This matrix $\bR$ can be used to write the Gram-Schmidt orthogonalization of $\bU(\bx)$ 
\begin{equation}
\bU(\bx) = \bV(\bx) \bR,
\label{eq:GS}
\end{equation}
which defines the orthonormal, causal basis of the approximation space \eqref{eq:F27}, 
\begin{equation}
\bV(\bx)=\big(\bv_0(\bx), \ldots, \bv_{n-1}(\bx) \big), 
\label{eq:F31p}
\end{equation}
with  $m-$dimensional row vector field components $\bv_j(\bx)$, $j =0,\ldots, n-1$,   called the ``orthonormal snapshots". Then, we have
\begin{equation}
\cPR :=
{ \bR^{-T} \bS \bR^{-1} } = \int_\Omega d\bx \bV(\bx)^T \cP \bV (\bx) = \big( \langle \langle \bv_j,\cP \bv_l \rangle \rangle \big)_{j,l =0,\ldots,n-1},  \label{eq:F32}
\end{equation}
where the index $-T$ denotes the inverse and transpose. The first equality in this equation is used to compute $\cPR$ from the data, and the second equality shows that it is a projection 
of the operator $\cP$. 

The first $n$ ROM states satisfy
\begin{equation}
\bR = \begin{pmatrix} \bR_{0,0} & \bR_{0,1} & \bR_{0,2} & \ldots & \bR_{0,n-1}\\
{\bf 0} & \bR_{1,1} & \bR_{1,2} & \ldots & \bR_{1,n-1} \\
{\bf 0} & {\bf 0} & \bR_{2,2} & \ldots & \bR_{2,n-1} \\
\vdots & \vdots & \vdots & \ldots &\vdots \\
{\bf 0}& {\bf 0} &{\bf 0} & \ldots & \bR_{n-1,n-1} 
\end{pmatrix} = \big(\bu^\RM_0, \ldots \bu^\RM_{n-1}\big) = \int_\Omega d\bx \bV(\bx)^T \bU(\bx) =\big( \langle \langle \bv_j,\bu_l \rangle \rangle \big)_{j,l =0,\ldots,n-1}, \label{eq:F33}
\end{equation}
and we note how the algebraic structure of $\bR$ captures the causal wave propagation: The $nm \times m$ column blocks of $\bR$ are indexed using the time instants 
$j \tau$, for $j = 0, \ldots, n-1$, while the $m \times nm$ row blocks of $\bR$ are indexed according to the range locations reached by the wavefront at these instants.  The first column block of $\bR$, which equals $\bu^{\RM}_0$, has all but the first block equal to zero, because the true snapshot $\bu_0(\bx)$ is supported near the array. The 
second column block of $\bR$, which equals $\bu^{\RM}_1$, has an additional nonzero block because the true snapshot $\bu_1(\bx)$ reaches some range in the medium, and so on.

\section{The  internal wave}
\label{sect:intern}
We now use the data driven ROM reviewed above to
estimate an internal wave that originates from the vicinity of an arbitrary point $\by \in \Omega_{\rm im}$ and propagates to the array through the true, unknown medium.

The best estimate that we could hope for would be
\begin{equation}
g^{\rm ideal}(t,\bx;\by):= \cos \big(t \sqrt{A(c)}\big) \delta^{f}_{\by}(\bx),
\label{eq:In1}
\end{equation}
with the initial state 
\begin{align}
g^{\rm ideal}(0,\bx;\by) = \delta^{f}_{\by}(\bx) &= \hat{f}^\12\big(\sqrt{A(c)} \big) \delta_\by(\bx) \approx \hat{f}^\12\big(\sqrt{A(c)} \big) 
 \bV \bV^T \delta_{\by}(\bx) ,
 \label{eq:In2}
\end{align}
given by the approximation 
\begin{equation}
\bV \bV^T\delta_{\by}(\bx) := \sum_{j=0}^{n-1} \bv_j(\bx) \bv_j^T(\by) = \bV(\bx) \bV^T(\by),
\label{eq:In2p}
\end{equation}
of $\delta_\by(\bx)$ in the space \eqref{eq:F27}, blurred a little by the pulse dependent operator $\hat{f}^\12\big(\sqrt{A(c)}\big)$. 
The same reasoning used for  the sensor functions \eqref{eq:F18} applies to  \eqref{eq:In2}
and shows that it is supported in the ball centered at $\by$, with $O(c(\by) t_f)$  radius. 
We are interested in evaluating the wave \eqref{eq:In1} at  $t > 0$ and the  sensor locations $\bx_r$, for $r = 1, \ldots, m$. However, we  cannot get exactly $g^{\rm ideal}(t,\bx_r;\by)$, because we do not know the approximation space \eqref{eq:F27}. The next proposition shows that  we can compute instead 
\begin{equation}
g(t,\bx;\by):= \cos \big(t \sqrt{A(c)}\big)  \delta_{\by}^{f,\RM}(\bx), \qquad
\delta_{\by}^{f,\RM}(\bx):= \hat f^{\12} \big(\sqrt{A(c)}\big) \delta_{\by}^{\RM}(\bx),
\label{eq:In3}
\end{equation}
at $\bx = \bx_r$, for $r =1,\ldots,m,$ 
where the approximation \eqref{eq:In2p}  of $\delta_\by(\bx)$  is replaced by the ``ROM point spread function"
\begin{equation}
\delta_{\by}^\RM(\bx):= \bV \bV_o^T \delta_{\by}(\bx) = \sum_{j=0}^{n-1} \bv_j(\bx) \bv_{o,j}^T(\by) = \bV(\bx) \bV_o^T(\by),
 \label{eq:In4}
\end{equation}
calculated with the orthonormal snapshots in the reference medium with known wave speed $c_o(\bx)$
\begin{equation}
\bV_o(\by)=\big(\bv_{o,0}(\by), \ldots, \bv_{o,n-1}(\by) \big).
\label{eq:F31pref}
\end{equation}

We will see in the next section that the tighter the focus of  $\delta_{\by}^\RM(\bx)$ around $\by$, the better the imaging using the internal wave $g(t,\bx;\by)$.  So when can we expect such a result? The answer lies in how well we can approximate the snapshots  in $\bU(\bx)$  in the reference space calculated for the known $c_o(\bx)$,
\begin{equation}
\mathscr{S}_o:=\mbox{range} \, \bU_o(\bx), \qquad \bU_o(\bx) = \big(\bu_{o,0}(\bx), \ldots, \bu_{o,n-1}(\bx)\big).
\label{eq:RefSp}
\end{equation}
If it is true that the approximation error is small, then the Gram-Schmidt orthogonalization, which is a stable procedure, gives 
\begin{equation}
\bv_j(\bx) \approx \bv_{o,j}(\bx), \qquad j = 0, \ldots, n-1,
\label{eq:approxV}
\end{equation}
and the ROM point spread function \eqref{eq:In4} is an approximation
of \eqref{eq:In2}.

We discuss in  \ref{ap:A} two setups where we can analyze explicitly the approximation 
\eqref{eq:approxV}: in a layered medium and in a waveguide. The error in these two cases is controlled
by the time step $\tau$, the separation between the sensors and the array aperture size. In more general settings we only have numerical evidence that if $\tau$ and the sensor separation are small enough and the aperture is large enough, then the ROM 
point spread function $\delta_{\by}^\RM(\bx)$ is peaked at $\by$. 

\begin{prop}
\label{prop.1}
Let $\bR$ be the block upper triangular Cholesky factor of the data driven mass matrix $\bM$, with block entries  given by \eqref{eq:F29}. The estimated internal wave \eqref{eq:In3} evaluated at the sensor locations and at the time instants $t_j = j \tau$, for $j=0,\ldots, n-1$, is given by the $m-$dimensional row vector field
\begin{equation}
\big(
g(t_j,\bx_1;\by), \ldots, g(t_j,\bx_m;\by)\big)
= \bV_o(\by) \bR \be_j ,
\label{eq:In5}
\end{equation}
where $\be_j \in \RR^{nm \times m}$ is the $(j+1)^{\rm th}$ column block of the $nm \times nm$ identity matrix ${\bf I}_{nm}.$
\end{prop}

\vspace{0.03in}\textbf{Proof:} Let us begin with the auxiliary $m-$dimensional column vectors
\begin{equation}
\bsig_l(\by):= \be_l^T \bR^{-1} \bV_o^T(\by), \qquad l = 0,\ldots, n-1,
\label{eq:P1}
\end{equation}
and note that
\begin{align}
\sum_{l=0}^{n-1} \bu_l(\bx)\bsig_l(\by) &= \bU(\bx) \sum_{l=0}^{n-1} \be_l\bsig_l(\by) \nonumber\\
&=\bU(\bx) \left(\sum_{l=0}^{n-1} \be_l \be_l^T \right) \bR^{-1} \bV_o^T(\by) \nonumber\\
&=\bU(\bx) \bR^{-1} \bV_o^T(\by) 
\nonumber\\
&=\bV(\bx) \bV_o^T(\by) = \delta_\by^\RM(\bx).
\label{eq:P2}
\end{align}
Here the first equality is by the definition \eqref{eq:F27} of $\bU(\bx)$,
in the third equality we used that the sum over $l$ equals the identity matrix and the last equality is due to the Gram-Schmidt orthogonalization \eqref{eq:GS}.
Applying the operator $\cos \big(t_j \sqrt{A(c)}\big)$ to both sides of \eqref{eq:P2} and using 
definition \eqref{eq:F22} we get
\begin{align}
\cos\big(t_j \sqrt{A(c)}\big)\delta_\by^\RM(\bx) &= \cos\big(t_j \sqrt{A(c)}\big) \sum_{l=0}^{n-1} \cos\big(t_l \sqrt{A(c)}\big) \bde^f(\bx) \bsig_l(\by) \nonumber \\
&= \frac{1}{2} \sum_{l=0}^{n-1} \left[ \cos\big((j+l)\tau \sqrt{A(c)}\big) +\cos\big(|j-l|\tau \sqrt{A(c)}\big)\right]  \bde^f(\bx) \bsig_l(\by) \nonumber \\ &=
\frac{1}{2} \sum_{l=0}^{n-1} \left[ \bu_{j+l}(\bx) +\bu_{|j-l|}(\bx)\right]\bsig_l(\by) .
\end{align}
Furthermore, using definition \eqref{eq:F20} and the self-adjointness of $A(c)$, we have 
\begin{align}
g(t_l,\bx_r;\by) &= \int_{\Omega} d \bx \, \delta_{\bx_r}(\bx) 
\hat f^{\12}\big(\sqrt{A(c)}\big)\cos\big(t_j \sqrt{A(c)}\big)\delta^\RM(\by) \nonumber\\
&=\int_{\Omega} d \bx \, \delta^f_{\bx_r}(\bx) \frac{1}{2} \sum_{l=0}^{n-1} \left[ \bu_{j+l}(\bx) +\bu_{|j-l|}(\bx)\right] \bsig_l(\by),
\end{align}
for all $r = 1,\ldots,m$. Gathering these results in an $m-$dimensional column vector and recalling definition \eqref{eq:F21} and the expression \eqref{eq:F28} of the data matrices, we obtain
\begin{equation}
\begin{pmatrix}
g(t_j,\bx_1;\by) \\
\vdots \\
g(t_j,\bx_m;\by)
\end{pmatrix}
= \frac{1}{2}  \sum_{l=0}^{n-1}\left( \bD_{j+l} +\bD_{|j-l|}\right)\bsig_l(\by) = 
 \sum_{l=0}^{n-1} \bM_{j,l} \bsig_l(\by), \label{eq:P3}
\end{equation}
where the last equality is by \eqref{eq:F29}. Finally, we substitute \eqref{eq:P1} in \eqref{eq:P3}
and use the Cholesky factorization of the mass matrix to get the result
\begin{align}
\begin{pmatrix}
g(t_j,\bx_1;\by) \\
\vdots \\
g(t_j,\bx_m;\by)
\end{pmatrix}
= \sum_{l=0}^{n-1} \bM_{j,l} \be_l^T \bR^{-1} \bV_o^T(\by) = \be_j^T \bM  \bR^{-1} \bV_o^T(\by) 
= \be_j^T \bR^T \bV_o^T(\by) = \Big( \bV_o(\by) \bR \be_j \Big)^T.
\end{align}
$\Box$

\section{Imaging with the internal wave}
\label{sect:normg}
In this section  we introduce a novel imaging approach
based on the internal wave estimated  in Proposition \ref{prop.1}. The imaging function is strikingly simple: it is the squared norm of this wave evaluated at the sensors
\begin{equation}
\cI(\by) = \sum_{r=1}^m \sum_{j=0}^{n-1} \left| g(j 
\tau,\bx_r;\by)\right|^2, \qquad \by \in \Omega_{\rm im}.
\label{eq:Im1}
\end{equation}
Moreover, $\cI(\by)$  is easy to compute and it does not even require the full ROM. It just uses the
Cholesky factor $\bR$ of the data driven mass matrix $\bM$ with block entries \eqref{eq:F29}, and 
the orthonormal snapshots \eqref{eq:F31pref} calculated  in the reference medium with known 
wave speed $c_o(\bx)$.

Our goal in this section is to analyze  \eqref{eq:Im1} and show that it gives an estimate of the location of reflective structures embedded in the host medium. The analysis is based on the continuum time approximation, where the sum over $j$ is replaced by an integral over time, and assumes a long 
enough duration of the measurements. We also suppose, as is typical in applications, that the wave speed 
is constant near the sensors and therefore the accessible boundary,
\begin{equation}
c(\bx) = c_o(\bx) = \bar c_o, \qquad \bx ~\mbox{near }\partial \Omega_{\rm ac}.
\label{eq:assC}
\end{equation}

We begin in section  \ref{sect:normg1} with the connection between the internal wave and the Green's function of the acoustic wave equation. This is useful for obtaining the main result in section \ref{sect:normg2}, where we relate $\cI(\by)$ to the time reversal process and we 
discuss its resolution. We end in section \ref{sect:normg3} with a comparison of $\cI(\by)$ and the imaging function of the backprojection approach introduced in \cite{druskin2018nonlinear}.
 
\subsection{The Green's function and its connection to the internal wave}
\label{sect:normg1}
The next lemma connects the internal wave \eqref{eq:In3} with the Green's function $G(t,\bx;\bz)$
of the acoustic wave equation, satisfying 
\begin{align}
\left[\frac{1}{c^2(\bx)} \partial_t^2 - \Delta_{\bx}\right] G(t,\bx;\bz) &= \delta(t) \delta_{\bz}(\bx), \qquad
t \in \RR, ~~ \bx \in \Omega, \label{eq:Im2} \\
 G(t,\bx;\bz) &= 0, \qquad t <0, ~~ \bx \in \Omega, \label{eq:Im3} \\
\partial_n G(t,\bx;\bz) &= 0, \qquad t \in \RR, ~~ \bx \in \partial \Omega_{\rm ac}, \label{eq:Im4} \\
G(t,\bx;\bz) &= 0, \qquad t \in \RR, ~~ \bx \in \partial \Omega_{\rm inac},\label{eq:Im5}
\end{align} 
where $\bz$ is an arbitrary point in $\Omega$ and $\Delta_{\bx}$
is the Laplace operator with respect to $\bx$. 

\begin{lem}
\label{lem.1}
Let $\by$ be any point in the imaging domain $\Omega_{\rm im}$, which supports the sought after reflective structures. The internal wave \eqref{eq:In3} evaluated at the sensor locations  satisfies
\begin{equation}
g(t,\bx_r;\by) = \partial_t \check{f}^{\frac{1}{2}}(t) \star_t \int_{\Omega} d \bz \, 
\frac{G(t,\bx_r;\bz)}{c(\bz)\bar c_o} \delta_\by^\RM(\bz), \qquad r = 1, \ldots, m,
\label{eq:Im6}
\end{equation}
for any $t>0$,
where we recall that $\check{f}^{\frac{1}{2}}(t)$ is defined in \eqref{eq:F19} and  $\star_t$ denotes convolution in time $t$. Recall also that $g$ is an even function in $t$.
\end{lem}

\vspace{0.03in}\textbf{Proof:} Let us begin with the even in time wave function
\begin{equation}
\cG_e(t,\bx;\bz) = \cos \big(t \sqrt{A(c)}\big) \delta_{\bz}(\bx),
\label{eq:evenG}
\end{equation}
and use linear superposition to write
\begin{equation}
\cos\big(t \sqrt{A(c)}\big) \delta_{\by}^\RM(\bx) = 
\int_\Omega d \bz \, \cG_e(t,\bx;\bz) \delta_{\by}^\RM(\bz).
\label{eq:Im7}
\end{equation}
The internal wave is
\begin{align}
g(t,\bx;\by) &= \hat f^{\12}\big(\sqrt{A(c)}\big)  \cos\big(t \sqrt{A(c)}\big) \delta_{\by}^\RM(\bx) 
\nonumber \\
&=  \hat f^{\12}\big(\sqrt{A(c)}\big) \int_\Omega d \bz \, \cG_e(t,\bx;\bz) \delta_{\by}^\RM(\bz)
\nonumber \\
&=\sum_{l=1}^{\infty}  \hat f^{\12} \big(\sqrt{\theta_l}\big) \cos\big(t\sqrt{\theta_l}\big) y_l(\bx) \lb y_l, \delta_{\by}^\RM \rb \nonumber \\
&=\sum_{l=1}^{\infty} \int_{-\infty}^{\infty} dt' \, \check{f}^{\12}(t') \cos \big(t'\sqrt{\theta_l}\big) \cos\big(t\sqrt{\theta_l}\big) y_l(\bx) \lb y_l, \delta_{\by}^\RM \rb \nonumber \\
&=\sum_{l=1}^{\infty} \int_{-\infty}^{\infty} dt' \, \check{f}^{\12}(t') \12  \left[
 \cos\big((t-t')\sqrt{\theta_l}\big) +  \cos\big((t+t')\sqrt{\theta_l}\big)\right]  y_l(\bx) \lb y_l, \delta_{\by}^\RM \rb 
 \nonumber \\
&=\sum_{l=1}^{\infty} \int_{-\infty}^{\infty} dt' \, \check{f}^{\12}(t') \cos\big((t-t')\sqrt{\theta_l}\big)y_l(\bx) \lb y_l, \delta_{\by}^\RM \rb,
\label{eq:Im8}
\end{align}
where we used  that operators of $A(c)$ commute,  as well as the spectral decomposition of $A(c)$, definition  \eqref{eq:F19}  and that $ \check{f}^{\12}(t) $ is even. Since we have
\begin{align} 
\cos\big((t-t') \sqrt{A(c)}\big) \delta_{\by}^\RM(\bx) = 
\int_\Omega d \bz \, \cG_e(t-t',\bx;\bz) \delta_{\by}^\RM(\bz) = \sum_{l=1}^{\infty}  \cos\big((t-t')\sqrt{\theta_l}\big)y_l(\bx) \lb y_l, \delta_{\by}^\RM \rb,
 \label{eq:Im8p}
 \end{align}
 pointwise in $t-t'$, and  $\check{f}^{\12}(t')$ has finite support, we can use the dominated convergence theorem to interchange the integral and sum in \eqref{eq:Im8} and get
 \begin{align} 
 g(t,\bx;\by) = \check{f}^{\12}(t) \star_t \int_\Omega d \bz \, \cG_e(t,\bx;\bz) \delta_{\by}^\RM(\bz).
  \label{eq:Im10}
 \end{align}
 
It remains to connect the Green's function $G(t,\bx,\bz)$ to \eqref{eq:evenG}, which is the even extension in time 
\begin{equation}
\cG_e(t,\bx;\bz) = \cG(t,\bx;\bz) + \cG(-t,\bx;\bz),
\label{eq:Im11}
\end{equation}
of the causal Green's function $\cG(t,\bx;\bz)$, satisfying 
\begin{align} 
\left[\partial_t^2 +A(c) \right]\cG(t,\bx;\bz) &= \delta'(t) \delta_\bz(\bx), \qquad
t \in \RR, ~~ \bx \in \Omega, \label{eq:Im12} \\
 \cG(t,\bx;\bz) &= 0, \qquad t <0, ~~ \bx \in \Omega, \label{eq:Im13} \\
\partial_n \cG(t,\bx;\bz) &= 0, \qquad t \in \RR, ~~ \bx \in \partial \Omega_{\rm ac}, \label{eq:Im14} \\
\cG(t,\bx;\bz) &= 0, \qquad t \in \RR, ~~ \bx \in \partial \Omega_{\rm inac}.\label{eq:Im15}
\end{align} 
We are interested in evaluating \eqref{eq:Im10} at the sensor locations, which are far from $\by$, where $\delta_{\by}^\RM(\bx)$ peaks. By causality and the finite speed of propagation we should have 
$g(t,\bx_r;\by)=0$ for time $t = O(t_f)$. For larger time, we conclude from 
\begin{equation}
\cG_e(t,\bx;\bz) = \cG(t,\bx;\bz), \qquad t >0,
\label{eq:Im16}
\end{equation}
and the $t_f$ duration of the pulse, that 
\begin{equation}
g(t,\bx_r;\by) = \check{f}^{\12}(t) \star_t \int_\Omega d \bz \, \cG(t,\bx_r;\bz) \delta_{\by}^\RM(\bz),
\qquad t >O(t_f), ~~ r = 1, \ldots, m.
\label{eq:Im17}
\end{equation}
The solutions of 
\eqref{eq:Im2}--\eqref{eq:Im5} and \eqref{eq:Im12}--\eqref{eq:Im15} are related by
\begin{equation}
\partial_t G(t,\bx;\bz) = c(\bx) c(\bz) \cG(t,\bx;\bz),
\label{eq:Im18}
\end{equation}
where we used the expression \eqref{eq:F3} of the operator $A(c)$ and the assumption \eqref{eq:assC}. The statement of the Lemma follows from the identity \begin{equation}
\check{f}^{\12}(t) \star_t \partial_t G(t,\bx_r;\bz) = \partial_t \check{f}^{\12}(t)  \star_t G(t,\bx_r;\bz).
\end{equation}
$\Box$

\subsection{Analysis of the imaging function}
\label{sect:normg2}
The expression of the imaging function \eqref{eq:Im1} is given in the next proposition,
obtained with the continuum time approximation and for a long duration of the measurements.

\begin{prop}
\label{prop.2}
The imaging function \eqref{eq:Im1} is approximated by
\begin{equation}
\cI(\by) \approx  \int_\Omega d \bz \int_\Omega d \bz' \,
\delta_{\by}^\RM(\bz) \delta_\by^\RM(\bz')\frac{\Gamma(\bz,\bz')}{c(\bz) c(\bz')},
\label{eq:MR}
\end{equation}
where
\begin{equation}
\Gamma(\bz,\bz') := \frac{1}{\tau \bar c_o^2}
\sum_{r=1}^m \int_0^{n\tau} dt  \iint_{\RR^2} ds ds'
\partial_s   \check{f}^{\frac{1}{2}}(s)
\partial_{s'}   \check{f}^{\frac{1}{2}}(s')
G(t-s,\bx_r;\bz) 
G(t-s',\bx_r;\bz') 
.
\label{eq:MR1}
\end{equation} 
\end{prop}

\vspace{0.03in}\textbf{Proof:} Approximating the sum over $j$ in  \eqref{eq:Im1}   by the integral in $t$,
we get 
\begin{align}
\cI(\by) \approx \frac{1}{ \tau} \sum_{r=1}^m \int_{0}^{n\tau} dt \, \left|g(t,\bx_r;\by)\right|^2,
\label{eq:Im20}
\end{align} 
and the result follows by substituting (\ref{eq:Im6}) in this expression.
$\Box$

\vspace{0.03in}
In order to explain why the imaging function  $\cI(\by)$ gives an image of the local changes of the velocity, we interpret its approximate  expression \eqref{eq:MR}  in terms of the result of the following time-reversal experiment:\\
$\bullet$ {\it First step.}
Consider the source function 
\begin{align}
\label{eq:TRSource}
n_1(t,\bx) =\frac{1}{\tau c(\bx) \bar{c}_o}  \partial_t \check{f}^{\frac{1}{2}}(t)  \delta_{\by}^\RM(\bx),
\end{align}
that is localized in space at $\bx$ in the support of $ \delta_{\by}^\RM$, and in time around $0$, at $t$ in the support of  $ \check{f}^{\frac{1}{2}}$.
 The wave field generated by this source satisfies
\begin{align}
\left[\frac{1}{c^2(\bx)} \partial_t^2 - \Delta_{\bx}\right] u_1(t,\bx) &=n_1(t,\bx) , \qquad
t \in \RR, ~~ \bx \in \Omega,  \\
u_1(t,\bx) &= 0, \qquad t \ll 0, ~~ \bx \in \Omega, \\
\partial_n u_1(t,\bx) &= 0, \qquad t \in \RR, ~~ \bx \in \partial \Omega_{\rm ac}, \\
u_1(t,\bx) &= 0, \qquad t \in \RR, ~~ \bx\in \partial \Omega_{\rm inac},
\end{align} 
and suppose we record it at the sensor locations $(\bx_r)_{r=1}^m$ for $t\in [0,n\tau]$.  Using the Green's function $G(t,\bx;\bz)$, the solution of \eqref{eq:Im2}-\eqref{eq:Im5},  we can write these recordings as 
\begin{align}
u_1(t,\bx_r) = \int_\RR ds  \int_\Omega d\bz G(t-s,\bx_r;\bz) n_1(s,\bz).
\end{align}
We then have from \eqref{eq:MR}-\eqref{eq:MR1} that 
\begin{align}
{\cal I}(\by) \approx \sum_{r=1}^m \int_0^{n\tau} dt  \int_\RR ds \int_\Omega \frac{d\bz}{c(\bz) \bar{c}_o}  
G(t-s,\bx_r;\bz) 
 \delta_{\by}^\RM(\bz) 
\partial_s   \check{f}^{\frac{1}{2}}(s)
 u_1(t,\bx_r),
\end{align}
and after  the change of variables $t\mapsto -t$, $s\mapsto -s$ and using that $\partial_s   \check{f}^{\frac{1}{2}}$ is odd, this equation becomes
\begin{align}
{\cal I}(\by) = - \int_\RR ds \int_\Omega \frac{d\bz}{c(\bz) \bar{c}_o}  
 \delta_{\by}^\RM(\bz) 
\partial_s   \check{f}^{\frac{1}{2}}(s)
 \int_\RR dt 
 \sum_{r=1}^m
 G(s-t,\bx_r;\bz) 
 u_1(-t,\bx_r) {\bf 1}_{[-n\tau,0]}(t).
 \label{eq:MR2}
\end{align}
Here ${\bf 1}_{[-n\tau,0]}(t)$ is the indicator function of the interval $[-n \tau,0]$, equal to $1$ when $t$ lies in this interval and $0$ otherwise.
\\
$\bullet$ {\it Second step.}
Consider the source function 
\begin{align}
n_2(t,\bx) =\sum_{r=1}^m 
 u_1(-t,\bx_r) {\bf 1}_{[-n\tau,0]}(t) \delta_{\bx_r}(\bx)  ,
 \end{align}
  that is localized in time in $[-n\tau,0]$ and in space at  the sensor locations $(\bx_r)_{r=1}^m$. This source transmits the time-reversed recorded signals $u_1(-t,\bx_r)$, and the generated wave field 
satisfies 
 \begin{align}
\left[\frac{1}{c^2(\bx)} \partial_t^2 - \Delta_{\bx}\right] u_2(t,\bx) &=n_2(t,\bx) , \qquad
t \in \RR, ~~ \bx \in \Omega,  \\
u_2(t,\bx) &= 0, \qquad t \ll 0, ~~ \bx \in \Omega, \\
\partial_n u_2(t,\bx) &= 0, \qquad t \in \RR, ~~ \bx \in \partial \Omega_{\rm ac}, \\
u_2(t,\bx) &= 0, \qquad t \in \RR, ~~ \bx\in \partial \Omega_{\rm inac} .
\end{align} 
% Assume that we record this wave $u_2(t,\bx)$ around the time $0$, at  $t$ in the support of $ 
% \check{f}^{\frac{1}{2}}$,  and around the position $\by$, at  $\bx$ in the support of $ \delta_{\by^\RM$.
We have, using again the Green's function, that 
\begin{align}
u_2(s,\bz) = \int_\RR dt  \int_\Omega d\bz' G(s-t,\bz;\bz') n_2(t,\bz')  
=
\sum_{r=1}^m  \int_\RR dt \, G(s-t,\bz;\bx_r) u_1(-t,\bx_r) {\bf 1}_{[-n\tau,0]}(t) ,
\end{align}
and  from (\ref{eq:MR2}) we find 
\begin{align}
{\cal I}(\by) = - \int_\RR ds  \int_\Omega  \frac{d\bz}{c(\bz) \bar{c}_o}  
 \delta_{\by}^\RM(\bz) 
\partial_s   \check{f}^{\frac{1}{2}}(s)
 u_2(s,\bz) .
 \label{eq:TRp}
\end{align}
If  $ \delta_{\by}^\RM $ is localized around $\by$, then this expression shows that we observe the time-reversed wave around time $0$,  at $t$ in the support of $\check{f}^{\frac{1}{2}}$,  and around $\by$.

Suppose that both the recording time window $n\tau$ and  the aperture of the sensor array are large enough. Then, the theory of time reversal for waves \cite{ammari} predicts that the refocused wave $u_2$ should be close to the original source \eqref{eq:TRSource}, but time-reversed, and therefore  
\begin{align}\label{eq:EqImagFunc}
{\cal I}(\by) \approx -  \int ds  \int  \frac{d\bz}{c(\bz) \bar{c}_o}  
 \delta_{\by}^\RM(\bz) 
\partial_s   \check{f}^{\frac{1}{2}}(s)
n_1(-s,\bz)  = 
 \int ds \big[\partial_s\check{f}^{\frac{1}{2}}(s)\big]^2
 \Big[  \int_\Omega \frac{d\bz}{c(\bz)^2 \bar{c}_o^2 \tau}  
 \delta_{\by}^\RM(\bz)^2  \Big].
\end{align}
This expression shows that  the imaging function $\cI(\by)$ is related to the local velocity at $\by$, provided  $\delta_{\by}^\RM$ is peaked at $\by$. If there are sharp and significant changes $c(\bz)-c_o(\bz)$ around $\by$, which correspond to reflective structures in the non-scattering host-medium,
they appear in $\cI(\by)$ with a resolution that depends on  $\delta_\by^\RM$. The more focussed 
this is at $\by$, the better the resolution. The other resolution controlling factors are the pulse width 
(support of $\check f^\12$), the recording time $n \tau$ and the array aperture, which determine how well the time reversal wave $u_2$ refocusses.

{
We display $\delta_\by^\RM(\bz)$ in the numerical results section \ref{sect:num} to show that it is indeed peaked at $\by$
if the time sample interval $\tau$ and the sensor separation are chosen properly. We also note 
that the ROM point spread function is insensitive to the variations $c(\bz) - c_o(\bz)$, so there is no cancellation of the 
wave speed  in \eqref{eq:TRp} and \eqref{eq:EqImagFunc}. To show this, we compute explicitly the $L^2(\Omega)$ norm  of $\delta^\RM_\by(\bz)$ using its definition \eqref{eq:In4} and the orthonormality of the components 
of $\bV(z)$, 
\[
\bv_j (\bz) = \big(  v_j^{(1)}(\bz), \ldots,  \bv_j^{(m)}(\bz) \big), \qquad j = 0, \ldots, n-1,
\]
which gives 
\[
\int_{\Omega} d \bz \, v_j^{(s)}(\bz) v_{j'}^{(s')}(\bz) = \delta_{j,j'} \delta_{s,s'}, \qquad \forall ~j,j' = 0, \ldots, n-1, ~~ s,s' = 1, \ldots, m,
\]
where $\delta_{j,j'}$ is the Kronecker delta. 
We obtain that 
\begin{align}
\big\| \delta^\RM_\by \|^2_{L^2(\Omega)} = \int_{\Omega} d \bz \, \left[\sum_{j=0}^{n-1} \sum_{s=1}^m v_j^{(s)}(\bz) 
v_{o,j}^{(s)}(\by) \right]^2 = \sum_{j=0}^{n-1} \sum_{s=1}^m \big[v_{o,j}^{(s)}(\by)\big]^2,
\label{eq:normROM}
\end{align}
where the right hand side depends on the orthonormal snapshots in the reference medium. 
}
%We display $\delta_\by^\RM(\bz)$ in the numerical results section \ref{sect:num} to show that it is indeed peaked at $\by$
%if the time sample interval $\tau$ and the sensor separation are chosen properly. We also plot there the norm $\|\delta_\by^\RM\|_{L^2(\Omega_{\rm im})}$ as a function of $\by \in \Omega_{\rm im}$, to  illustrate that it is insensitive to the sharp variations of $c$. Therefore, the dependence on the wave speed  in \eqref{eq:TRp} and \eqref{eq:EqImagFunc} does not cancel out and $\cI(\by)$ images  as explained above.}

\subsection{Comparison with backprojection imaging}
\label{sect:normg3}
The backprojection imaging function introduced in \cite{druskin2018nonlinear}  is given by 
\begin{equation}
\cI^\BP(\by) = \bV_o(\by) \big( \cPR - \boldsymbol{\cP}_o^\RM\big) \bV_o^T(\by),
\label{eq:IM30}
\end{equation}
where $\boldsymbol{\cP}_o^\RM$ is the ROM propagator calculated in the reference medium 
with known wave speed $c_o(\bx)$. To compare it with our imaging function $\cI(\by)$, let us rewrite 
\eqref{eq:IM30}  using equation \eqref{eq:F32} for the ROM propagator and the analogue of \eqref{eq:In2p} in the reference medium,
\begin{equation}
 \delta_{o,\by}(\bx) := \bV_o \bV_o^T \delta_{\by}(\bx) = \sum_{j=0}^{n-1} \bv_{o,j}(\bx) \bv_{o,j}^T(\by) = \bV_o(\bx) \bV_o^T(\by).
 \label{eq:IM32b}
\end{equation}
We obtain that
\begin{align}
\cI^\BP(\by) &= \bV_o \bV^T \cP  \bV \bV_o^T(\by) -  \bV_o \bV_o^T \cP_o  \bV_o \bV_0^T(\by)
\nonumber \\
 &= \sum_{j=0}^{n-1} \bv_{o,j}(\by) \langle \langle \bv_j, \cP \delta_{\by}^\RM \rangle \rangle 
 -  \sum_{j=0}^{n-1} \bv_{o,j}(\by) \langle \langle \bv_{o,j}, \cP_o \delta_{o,\by}\rangle \rangle,
\label{eq:IM31}
\end{align}
where $\cP_o$ is the wave propagator operator in the reference medium.  

Let us explain the meaning of the two terms in the right hand side of \eqref{eq:IM31}: 
The last term models the wave $\cP_o \delta_{o,\by}(\bx)$ with 
initial state $\delta_{o,\by}(\bx) := \bV_o \bV_o^T \delta_\by(\bx)$  peaked around $\by$, and propagated in the reference medium for the duration $\tau$. This wave is then projected on the reference space $\mathscr{S}_o$ using $\bV_o \bV_o^T$, and the result is evaluated at $\by$. 
The first term in \eqref{eq:IM31} involves the internal wave 
\begin{equation}
\cP \delta_{\by}^\RM(\by) = \cos \big( \tau \sqrt{A(c)} \big) \delta_{\by}^\RM(\bx),
\label{eq:IM32}
\end{equation}
that is similar to our wave $g(\tau,\bx;\by)$ given in \eqref{eq:Im3}. Ideally, this wave would be 
projected in the space $\mathscr{S}$, but since $\bV(\bx)$ is unknown, the projection $\bV \bV^T$ is replaced by $\bV_o \bV^T$, based on the expectation that the  approximation 
\eqref{eq:approxV} holds. 

By %causality 
{hyperbolicity}
and the short duration $\tau$ of propagation of the waves involved in \eqref{eq:IM31}, 
both terms described above should be affected mostly by the medium in the vicinity of $\by$. Therefore, by taking the difference of the terms, the backprojection imaging function is designed to sense changes $c(\bx)-c_o(\bx)$ in  the vicinity of the imaging point, like $\cI(\by)$. 

The numerical results in section \ref{sect:num} show that $\cI(\by)$ and $\cI^\BP(\by)$ perform similarly,  although $\cI(\by)$ has better cross-range resolution. They both outperform the reverse-time migration imaging approach, in the sense that they do not suffer from multiple scattering  artifacts. However, $\cI^\BP(\by)$ has the disadvantage that it is expensive to compute, because it 
involves the ROM propagator $\cPR$. The calculation of $\cPR$ is given in \eqref{eq:F32} and involves the unstable step of taking the inverse of $\bR$.  This requires careful additional 
processing to mitigate noise in the data, as explained in  \cite{borcea2019robust}. The imaging function $\cI(\by)$ 
is easy to compute and is also more robust to noise. 
It only needs the Cholesky factor $\bR$ of the mass matrix $\bM$, and not its inverse.  While the data driven $\bM$ may not be symmetric and positive definite due to noise, it can be easily transformed to such a matrix using for example the singular value decomposition,  and then $\bR$ can be obtained using the block Cholesky factorization algorithm described for example in \cite[Appendix B]{DtB}.

\section{Pixel scanning type imaging}
\label{sect:pixel}
In this section we use the internal wave estimated in Proposition \ref{prop.1} as a steering control 
at the array, for focusing the wave at the imaging points $\by \in \Omega_{\rm im}$. Such steering can be implemented experimentally and can be used for imaging in a pixel scanning manner. 

The basic idea is the principle of time reversal: Since we know the internal wave $\big(g(t,\bx_r;\by)\big)_{r=1,\ldots, m}$ originating from the vicinity of $\by 
\in \Omega_{\rm im}$,  we can just time reverse it and re-emit it into the medium, where it refocuses near $\by$.  If $\by$ lies near a reflector, the refocussed wave will be reflected back towards the array, where it can be measured. The reflector location can then be estimated from the peaks of the 
``pixel scanning" imaging function $\cI^{\PS}(\by)$ defined below, which  matches the reflected wave measured at $\bx_r$ with $g(t,\bx_r;\by)$, for $r = 1, \ldots, m$.

\subsection{Imaging algorithm}
The calculation of the imaging function $\cI^\PS(\by)$ is carried out with the following steps:
\begin{enumerate}
\item[(1)] Compute the mass matrix $\bM$ from the data collected at the array, as given in 
equation \eqref{eq:F29}.
\item[(2)] Compute the Cholesky factorization \eqref{eq:F31} and store $\bR$. 
\item[(3)] Compute $\bV_o(\bx)$ by carrying out the Gram-Schmidt orthogonalization of $\bU_o(\bx)$
computed by solving the wave equation in the reference medium with known wave speed $c_o(\bx)$. This is especially easy to do if $c_o(\bx) = \bar c_o$ for all $\bx \in \Omega$. 
\item[(4)] For each $\by \in \Omega_{\rm im}$ compute the internal wave $g(t,\bx_r; \by)$ using equation \eqref{eq:In5}, for  $r = 1, \ldots, m$.
\item[(5)] Define the control at the array for focusing at $\by$
\begin{equation}
\mathfrak{F}(t,\bx_s;\by):= {\bf 1}_{[0,n \tau]}(t) g(n\tau-t,\bx_s;\by), \qquad s = 1, \ldots, m.
\label{eq:PS1}
\end{equation}
Measure the wave $\gamma(t,\bx;\by)$ at the the sensors $\bx=\bx_r$,  after using the illumination 
\eqref{eq:PS1}.
\item[(6)] Calculate the imaging function
\begin{equation}
\cI^{\PS}(\by):= \sum_{r=1}^m \int_0^{n \tau} dt \, \gamma(n \tau + t,\bx_r;\by) g(t,\bx_r;\by).
\label{eq:PS2}
\end{equation} 
\end{enumerate}

Note that equation \eqref{eq:In5}  gives  the internal wave at the 
discrete time instants $t = j \tau$, for $j = 0, \ldots, n-1$. If $\tau$ is small enough, we can use interpolation to get $\mathfrak{F}(t,\bx_s;\by)$ at $ t \in [0, n \tau]$. 
Note also that at step (4),  the measurements should be for the acoustic pressure $c(\bx) \gamma(t,\bx;\by)$. Since the wave speed at the sensors equals the known constant $\bar c_o$, those measurements determine $\gamma(t,\bx_r;\by)$, for $r = 1, \ldots, m$.

\subsection{Expression of the refocusing and imaging functions}
\label{subsec:computegamma}
 The mathematical model of the wave $\gamma(t,\bx;\by)$ measured at step (4) of the algorithm is the solution of the wave equation 
\begin{align}
\partial_t^2 \gamma(t,\bx;\by) + A(c) \gamma(t,\bx;\by) &= \partial_t \sum_{s=1}^m \mathfrak{F}(t,\bx_s;\by), \qquad t > 0, ~~\bx \in \Omega, \label{eq:PS3} \\
\gamma(t,\bx;\by) &= 0, \qquad t < 0, ~~ \bx \in \Omega, \label{eq:PS4} \\
\partial_n \gamma(t,\bx;\by) &= 0, \qquad t > 0, ~~ \bx \in \partial \Omega_{\rm ac}, 
\label{eq:PS5} \\
\gamma(t,\bx;\by) &= 0, \qquad t > 0, ~~ \bx \in \partial \Omega_{\rm inac}.
\label{eq:PS6}
\end{align}
We now show that this wave focuses near $\by$ at time $t = n \tau$. 

Using the Green's function $\cG(t,\bx;\bz)$ defined in equations 
\eqref{eq:Im12}--\eqref{eq:Im15}, we can write using linear superposition that 
\begin{equation}
\gamma(t,\bx;\by) = \sum_{s=1}^m \cF(t,\bx_s;\by) \star_t \cG(t,\bx;\bx_s),
\label{eq:PS8}
\end{equation}
where $\cF(t,\bx_s;\by)$ is given by \eqref{eq:PS1} in terms of the internal wave \eqref{eq:In3}.
A calculation similar to that in the proof of Lemma \ref{lem.1} gives that 
\begin{equation}
g(t,\bx_s;\by) = \int_\Omega d \bz \,  \cG(t,\bx_s;\bz) \delta_\by^{f,\RM}(\bz),
\label{eq:PS7}
\end{equation}
and substituting the result in \eqref{eq:PS8} we get 
\begin{align}
\nonumber
\gamma(t,\bx;\by) 
&= \int_{\Omega} d \bz \, \delta_\by^{f,\RM}(\bz) \int_0^{n \tau} dt' \sum_{s=1}^m 
\cG(n \tau - t',\bx_s;\bz) \cG(t - t',\bx;\bx_s) \\
&= \int_{\Omega} d \bz \, \delta_\by^{f,\RM}(\bz) \int_0^{n \tau} dt' \sum_{s=1}^m 
\cG(n \tau - t',\bz;\bx_s) \cG(t - t',\bx;\bx_s) ,
\label{eq:PS9}
\end{align}
where
we have used the reciprocity relation
$\cG(t',\bx_s;\bz) = \cG(t',\bz;\bx_s)$. 
This clearly peaks at the instant $t = n \tau$, when the two Green's functions are in sync, 
and at points $\bx \approx \bz$ in the support of $\delta_\by^{f,\RM}(\bz)$ defined in 
\eqref{eq:In3}. Similar reasoning to that used in section \ref{sect:setup.1} for the sensor functions
\eqref{eq:F18} gives that $\delta_\by^{f,\RM}(\bz)$  has a slightly larger support than $\delta_\by^\RM(\bx)$, by an $O(c(\by) t_f)$ radius.

\vspace{0.03in}
The expression of the imaging function follows once we use  \eqref{eq:PS7}-\eqref{eq:PS9} in \eqref{eq:PS2}
\begin{align}
\cI^\PS(\by) &= \sum_{r=1}^m \sum_{s=1}^m \int_{\Omega} d \bz \, \delta^{f,\RM}_\by(\bz) \int_{\Omega} d \bz' \, \delta^{f,\RM} _\by(\bz') \int_0^{n \tau} dt \int_0^{n \tau} dt' \, \cG (t,\bx_r;\bz)\cG(t',\bz';\bx_s) \cG(t +t',\bx_r;\bx_s) \nonumber \\
&\approx \sum_{r=1}^m \sum_{s=1}^m \int_{\Omega} d \bz \, \delta^{f,\RM}_\by(\bz) \int_{\Omega} d \bz' \, \delta^{f,\RM} _\by(\bz')  \cG (-t,\bx_r;\bz) \star_t \cG(-t,\bz';\bx_s) \star_t \cG(t,\bx_r;\bx_s)  \big|_{t=0},
\label{eq:PS10}
\end{align}
where the approximation is for large enough $n\tau$.
% and we used the reciprocity relation \[\cG(t',\bz';\bx_s) = \cG(t',\bx_s;\bz') .\] 
As was the case in the previous section,  the ROM point spread function $\delta_\by^\RM(\bz)$ plays an important role in the imaging function. If $\delta_\by^\RM(\bz)$ is sharply peaked at $\by$,  so is $\delta_\by^{f,\RM}(\bz)$ and we have a contribution 
to \eqref{eq:PS10} from points $\bz \approx \bz' \approx \by$.  Then, we can interpret the terms in \eqref{eq:PS10} as follows: The first time convolution
\[
\cG (t,\bx_r;\bz) \star_t \cG(t,\bz';\bx_s) \approx \cG (t,\bx_r;\by) \star_t \cG(t,\by;\bx_s)
\]
models the wave propagating from the source at $\bx_s$  to $\by$, where we suppose there is a reflector, it presumably scatters there and then propagates back to the  receiver at $\bx_r$ in the array. The second time convolution matches this wave with $\cG(t,\bx_r;\bx_s)$, which models 
the echoes received at $\bx_r$, due to the illumination from $\bx_s$. If indeed there is a scatterer at 
$\by$, then there should be an arrival in $\cG(t,\bx_r;\bx_s)$ that is synchronous to that in  $ \cG (t,\bx_r;\by) \star_t \cG(t,\by;\bx_s)$, and we will get a large contribution to $\cI^\PS(\by)$. 

\begin{rem}
\label{rem.1}
The imaging function $\cI^\PS(\by)$ resembles that of the reverse-time migration approach,
where the array data, modeled by $ f(t) \star_t \cG (t,\bx_r;\bx_s), $ are migrated 
to the imaging point $\by$ in the reference medium
\begin{equation}
\cI^\RTM(\by) \approx \sum_{r=1}^m \sum_{s=1}^m \cG_o (-t,\bx_r;\by) \star_t \cG_o(-t,\by;\bx_s) \star_t f(t) \star_t \cG(t,\bx_r;\bx_s) \big|_{t=0}.
\label{eq:PS11}
\end{equation} 
In \eqref{eq:PS10} we use the Green's function $\cG$ in the true medium and not the reference one,
which should give a better result. However, we cannot obtain the ideal ``time-reversal" function 
\begin{equation}
\cI^\TR(\by) \approx \sum_{r=1}^m \sum_{s=1}^m \cG(-t,\bx_r;\by) \star_t \cG(-t,\by;\bx_s) \star_t f(t) \star_t \cG(t,\bx_r;\bx_s) \big|_{t=0}.
\label{eq:PS13}
\end{equation} 
Instead, we have  the blurrier version \eqref{eq:PS10},  where we integrate over points in the support of $\delta_\by^{f,\RM}(\bx)$.
\end{rem}

\begin{rem}
\label{rem.2}
The imaging functions $\cI(\by)$ and $\cI^\BP(\by)$ discussed in section \ref{sect:normg} are 
quite different than $\cI^{\PS}(\by)$ and $\cI^\RTM(\by)$. They are designed to be sensitive only to changes of the wave speed in the vicinity of the imaging point $\by$, and are not affected by 
the arrivals of the multiply scattered echoes in the medium. Such echoes are the cause of ghost reflectors present in the images formed with all three functions \eqref{eq:PS10}-\eqref{eq:PS13}, as we show with numerical simulations 
in section \ref{sect:num}.

\end{rem}

%------------------------
\section{Numerical results}
\label{sect:num}
In this section we present numerical results in two-dimensions. The setup mimics that in 
Fig. \ref{fig:setup}, with a rectangular domain $\Omega$ and the accessible boundary near the array, modeled as sound hard. The inaccessible boundary is sound soft and consists of two side boundaries aligned with the range direction, and a remote boundary, parallel to the array, which does not affect the waves over the duration $(2n-1)\tau$ of the data gather. The side  boundaries are close enough to each other to  play a role in  the simulations shown in sections \ref{sect:num1}-\ref{sect:num3} and thus cause a waveguide effect. We also present in section \ref{sect:num2} simulations for well separated side boundaries, that have no effect  for $t \in (0, (2n-1)\tau)$.

The reference (host) medium is homogeneous, with constant wave speed $\bar c_o$. The unknown wave speed $c(\bx)$ varies with the simulation and is displayed in the figures below.  All length scales 
are in units of the central wavelength $\la_c = 2 \pi \bar c_o/\om_c$. 
The probing pulse is 
\begin{equation}
f(t) = \frac{\sqrt{2\pi}}{2} \exp \Big(-\frac{t^2 B^2}{2} \Big) \cos(\om_c t), \qquad B = 0.25 \om_c.
\label{eq:N1}
\end{equation}
The data are generated by solving the wave equation for the acoustic pressure $\big(p^\ss(t,\bx)\big)_{s=1, \ldots, m}$, using a time domain,  second order centered
finite differences scheme, on a square mesh with size $\la_c/16$. The time steps are chosen to satisfy the Courant Friedrichs Lewy (CFL) condition.

The ROM construction is as described in \cite{borcea2019robust}. The mass matrix $\bM$ may be ill conditioned, especially if the time sample interval $\tau$ and the sensor separation is too small. Since the ROM computation involves the inverse of the Cholesky factor $\bR$ of $\bM$, even for noiseless data it requires regularization. For the calculation of the internal wave we only need $\bR$, so less regularization is needed.  

For noisy data,  formula \eqref{eq:F29} gives a symmetric data driven mass matrix that may not be positive definite. 
Let us call this matrix  $\tilde \bM$ and consider its eigenvalue decomposition $\tilde \bM = \bW \tilde \bLa \bW^T,$ where $\tilde \bLa$ is the diagonal matrix of the eigenvalues $\big(\tilde \Lambda_j \big)_{j=1}^{nm}$ and $\bW$ is the orthogonal matrix of the eigenvectors. 
We transform $\tilde \bM$ to a positive definite matrix $\bM$  used in the computation of the internal wave and the ROM as follows:  Set a threshold of the lowest acceptable positive eigenvalue $\Lambda_{\rm min}$ and define 
\begin{equation}
\bM =  \bW \bLa \bW^T,
\label{eq:regulariz}
\end{equation}
where $\bLa = \mbox{diag}(\Lambda_1, \ldots,  \Lambda_{nm})$ and $\Lambda_j = \max \{\tilde \Lambda_j, \Lambda_{\rm min}\}$, 
for $j = 1, \ldots, nm.$

\subsection{Imaging in a waveguide setting}
\label{sect:num1}
\begin{figure}[t]
\centering
\raisebox{0.3in}{\includegraphics[width=0.45\textwidth]{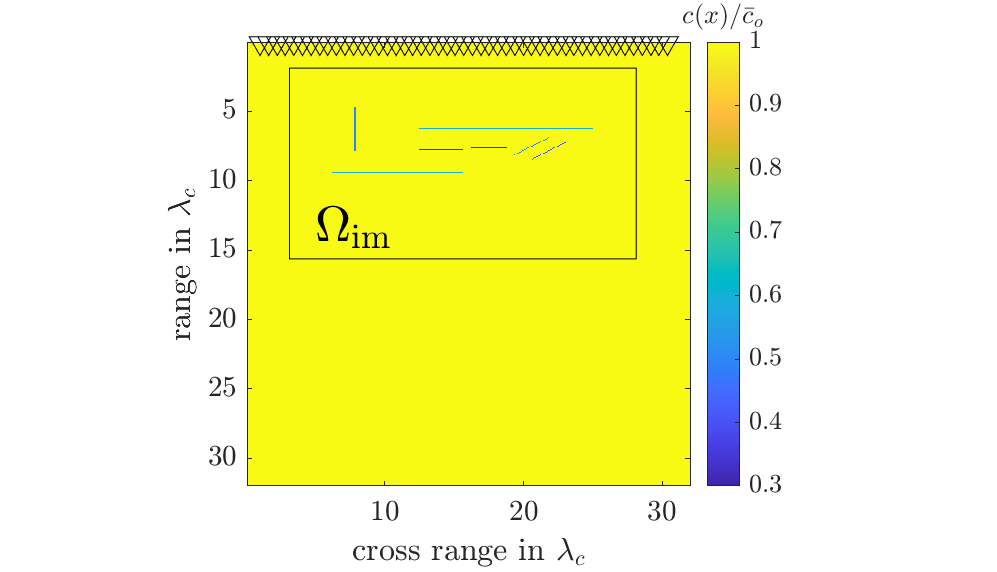}}
\includegraphics[width=0.54\textwidth]{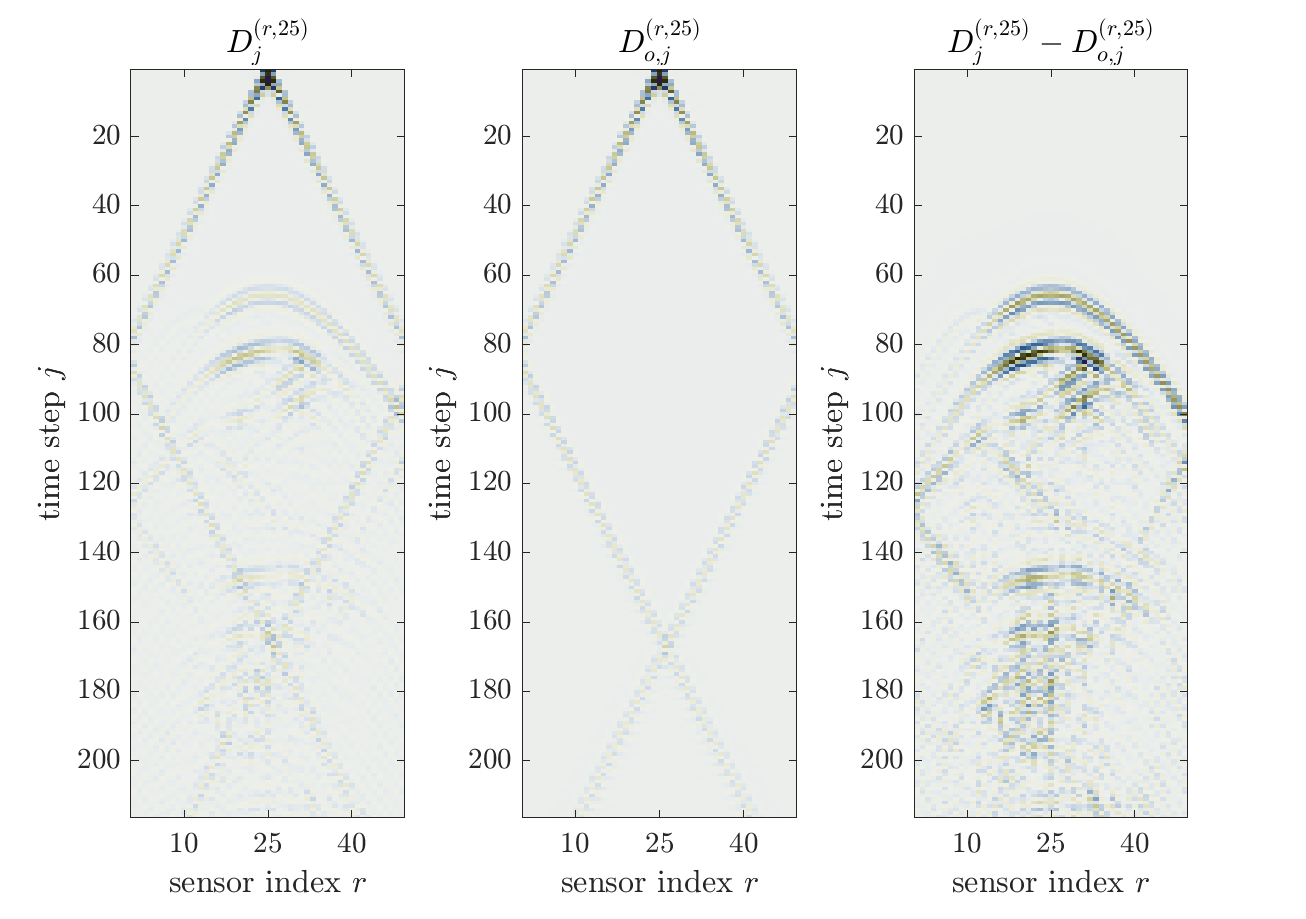}
\vspace{-0.25in}\caption{Left: Illustration of the setup: The array of $m = 49$ sensors (indicated with triangles) lying near the accessible boundary probes a medium with wave speed $\bar c_o$, containing a few thin 
reflecting structures, modeled by the low velocity shown in the color bar. Right: The data corresponding to the illumination from the center element in the array. We show it for the medium with the reflectors, the reference medium and the difference between the two.
}
\label{fig:Sim1}
\end{figure}

The numerical results in this section are for the setup illustrated in Fig. \ref{fig:Sim1}, where the side boundaries are close enough to play a role over the duration of the experiment, hence the name 
waveguide setting.  We consider first a large aperture size $a = 30 \la_c$, with beginning and end  at distance $\la_c$ from the side walls, and  containing $m = 49$ equidistantly placed sensors.  The time 
sample interval is $\tau = 0.4 \pi/\om_c,$ corresponding to $5$ points per carrier period.  
However, we also test how the aperture size, the separation between the sensors and $\tau$ affect the results, so we give the values of $a$, $m$ and $\tau$ in the captions of the figures. 

We display in the right plot of Fig. \ref{fig:Sim1} the data $D_j^{(r,s)}$ for $j = 0, \ldots, n-1$, $ r = 1, \ldots, m$ and $s = 25$, which indexes the center sensor in the array. We also show for comparison 
the data  in the reference medium. Note the echoes from the side walls that are present in the true and the reference medium, and the echoes from the sought after reflectors that are emphasized in the data differences.

\begin{figure}[t]
\centering
\includegraphics[width=0.4\textwidth]{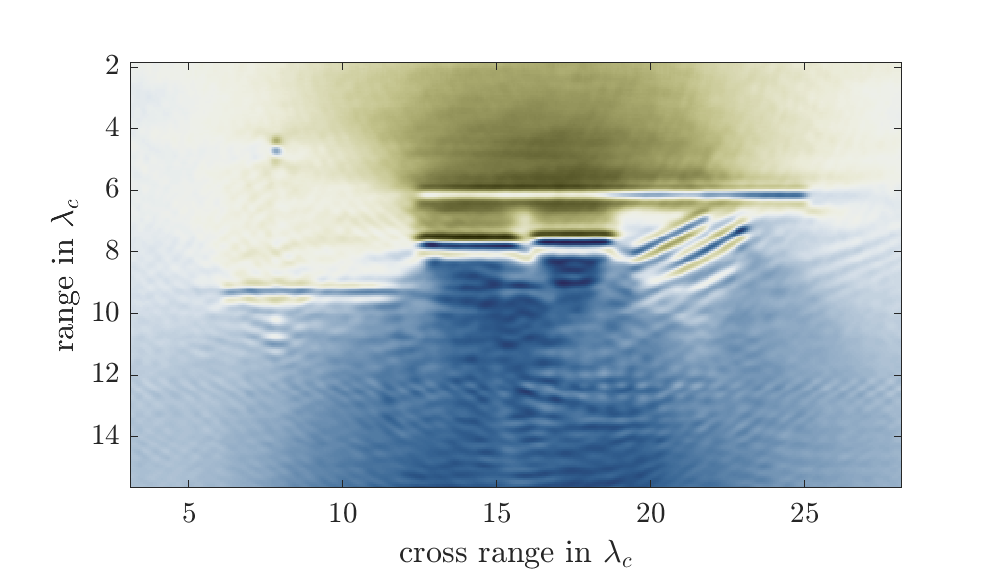}
\includegraphics[width=0.4\textwidth]{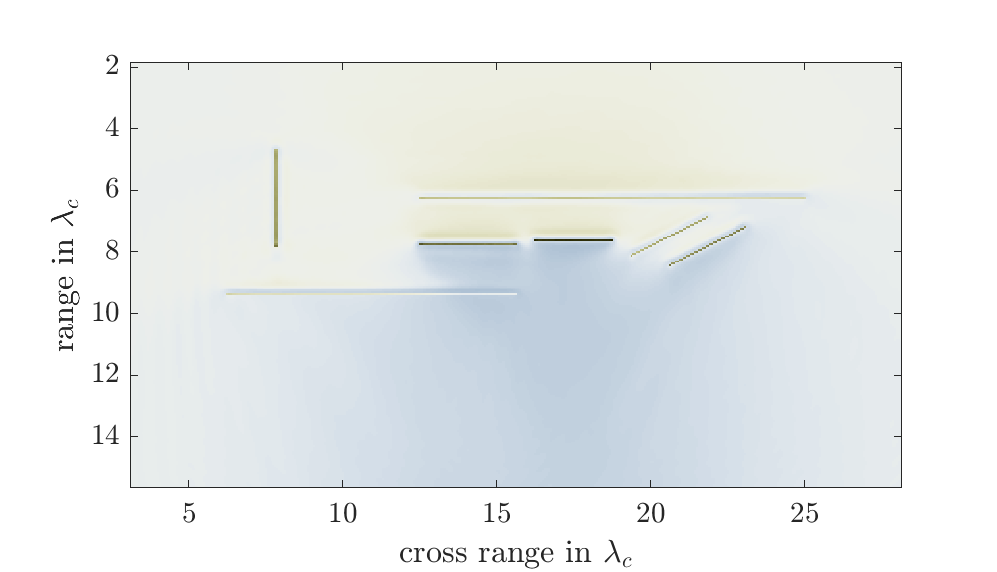}
\vspace{-0.1in}
\caption{Imaging function $\cI(\by)$ (left) and $\cI^{\rm ideal} (\by)$ (right) for the setup shown in Fig. \ref{fig:Sim1}. The aperture length is $a = 30 \la_c$ and the array has $m = 49$ sensors. The data are sampled in time at interval $\tau = 0.4 \pi/\om_c$. }
\label{fig:Sim2}
\end{figure}

In Fig. \ref{fig:Sim2} we display the imaging function $\cI(\by)$ defined  in \eqref{eq:Im1} and  the analogue function 
\begin{equation}
\cI^{\rm ideal}(\by) = \sum_{r=1}^m \sum_{j=0}^{n-1} \left| g^{\rm ideal}(j
\tau,\bx_r;\by)\right|^2, \qquad \by \in \Omega_{\rm im},
\label{eq:Num2}
\end{equation}
defined in terms of the ``ideal" internal wave \eqref{eq:In1} that cannot be computed. We can infer from Proposition \ref{prop.2} that there is  only one difference between  these functions:
The ROM point spread function $\delta_\by^\RM$ in the expression \eqref{eq:MR} of $\cI(\by)$ is replaced by 
the projection \eqref{eq:In2p} of $\delta_\by$ in the expression of $\cI^{\rm ideal}(\by)$. 
Due to the excellent focussing of \eqref{eq:In2p}, we see  that $\cI^{\rm ideal}(\by)$ gives a very sharp  (photo-like) estimate of the reflectors, whereas the image $\cI(\by)$ is a blurrier estimate. Moreover, $\cI(\by)$ captures only the top of the vertical reflector, and the unobstructed part of the bottom reflector.

Note that both plots in Fig. \ref{fig:Sim2}  display shadows of the reflectors and  have larger values near the array, because of the energy trapped there. This is less visible in $\cI^{\rm ideal}(\by)$, because its peak values are higher. To remove this effect, we display henceforth the derivative of the images in the range direction. This derivative is computed after smoothing the image in range by a convolution with a Gaussian function, with standard deviation $0.05\lambda_c$.

\begin{figure}[h]
\centering
\includegraphics[width=0.4\textwidth]{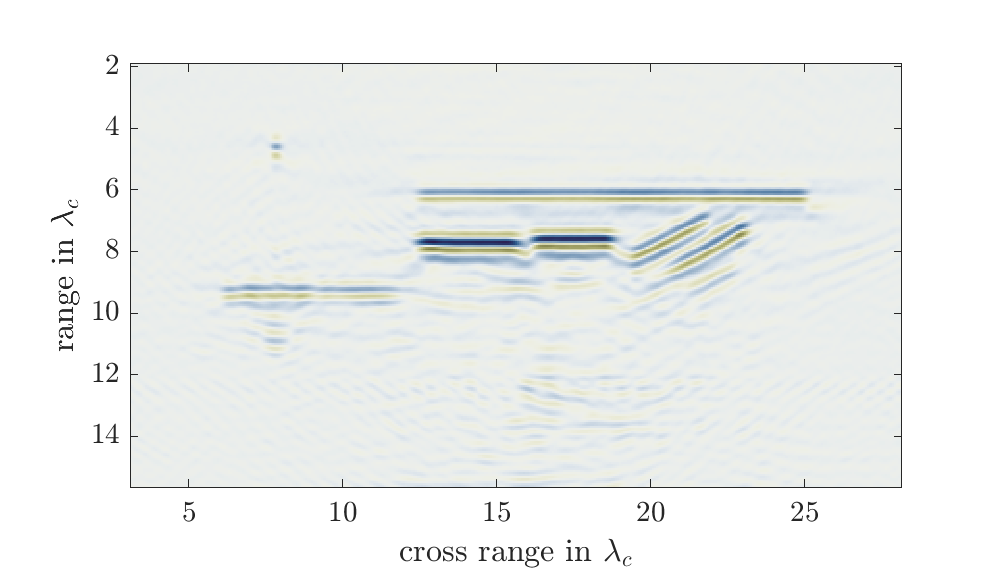}
\includegraphics[width=0.4\textwidth]{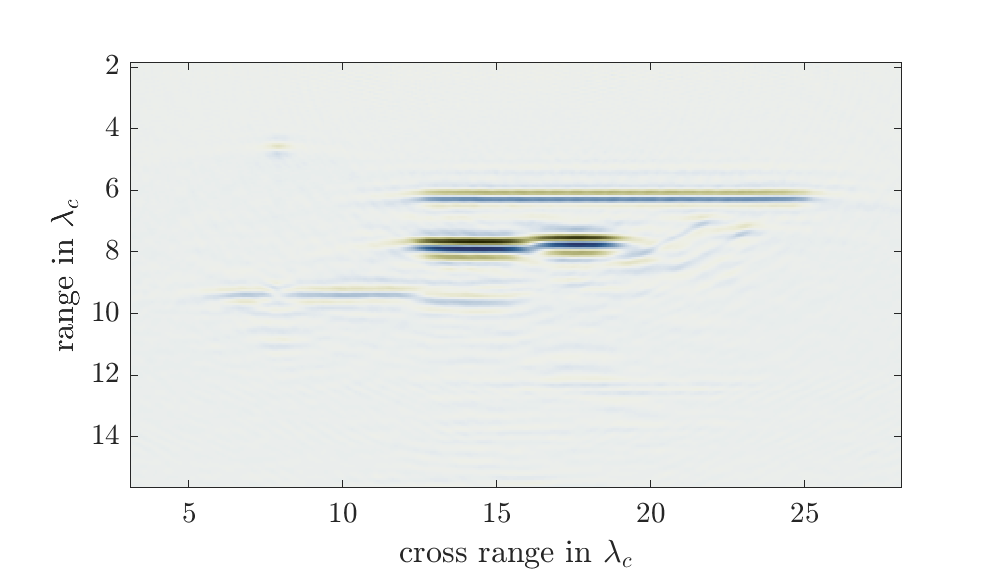}\\

\includegraphics[width=0.4\textwidth]{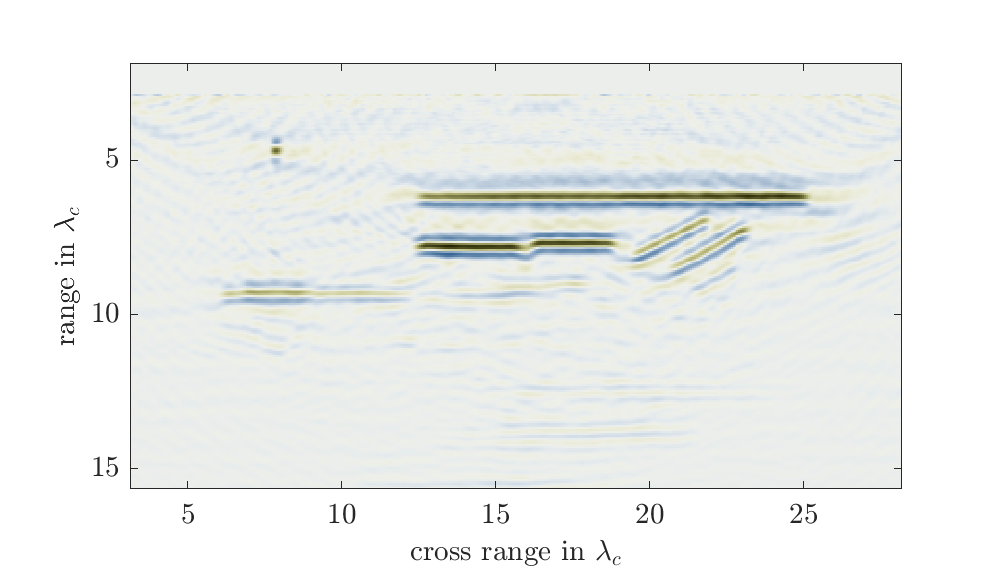}
\includegraphics[width=0.4\textwidth]{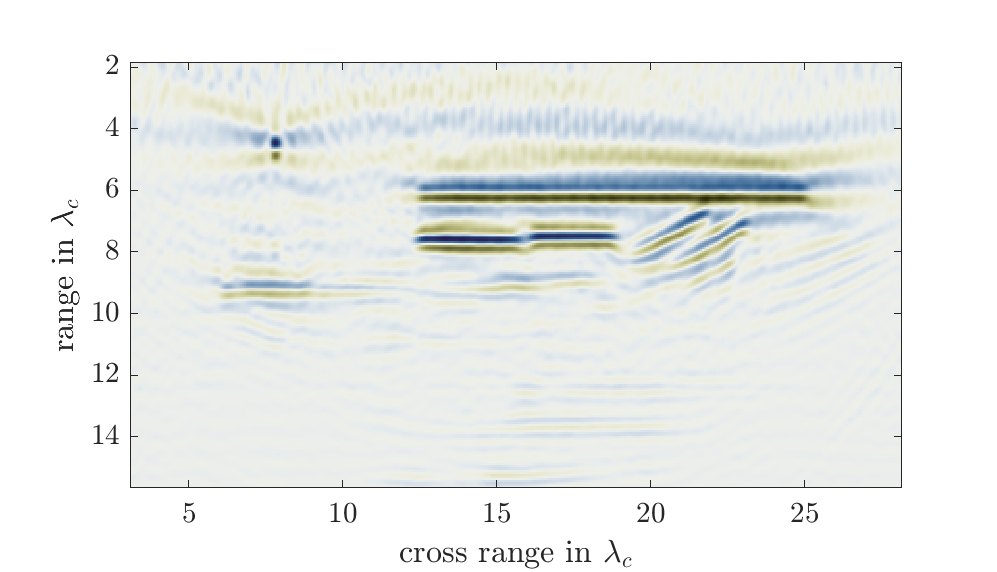}
\vspace{-0.1in}
\caption{Top plots: Range derivative of the imaging function $\cI(\by)$ (left), the backprojection image $\cI^{\rm BP}(\by)$ (right). Bottom plots: The range derivative of the pixel scanning image  $\cI^{\rm PS}(\by)$ (left) and 
the reversed time migration image  $\cI^{\small {\rm RTM}}(\by)$ (right) for the setup shown in Fig. \ref{fig:Sim1}. The aperture length is $a = 30 \la_c$ and the array has $m = 49$ sensors. The time sample interval is $\tau =0.4 \pi/\om_c.$}
\label{fig:Sim3}
\end{figure}

The plots in Fig. \ref{fig:Sim3} compare the four imaging functions discussed in the paper: 
$\cI(\by)$, $\cI^\BP(\by)$, $\cI^\PS(\by)$ and $\cI^\RTM(\by)$. They all localize the reflectors, 
with the exception of the vertical one, whose top is the only visible part, and the obstructed part of the bottom one. However, the result given by 
the computationally inexpensive imaging function $\cI(\by)$ is the better one, because: (1) It gives a better separation of the two nearby horizontal reflectors; (2) It displays clearly the oblique reflectors; 
(3) It does not have the ghost reflector seen in $\cI^\PS(\by)$ and especially $\cI^\RTM(\by)$, due to the reverberation between the top reflector and the accessible boundary. The backprojection image is 
also free of the ghost, but its cross-range resolution is worse and it barely sees the oblique reflectors.

\begin{figure}[tb]
\centering
\includegraphics[width=0.4\textwidth]{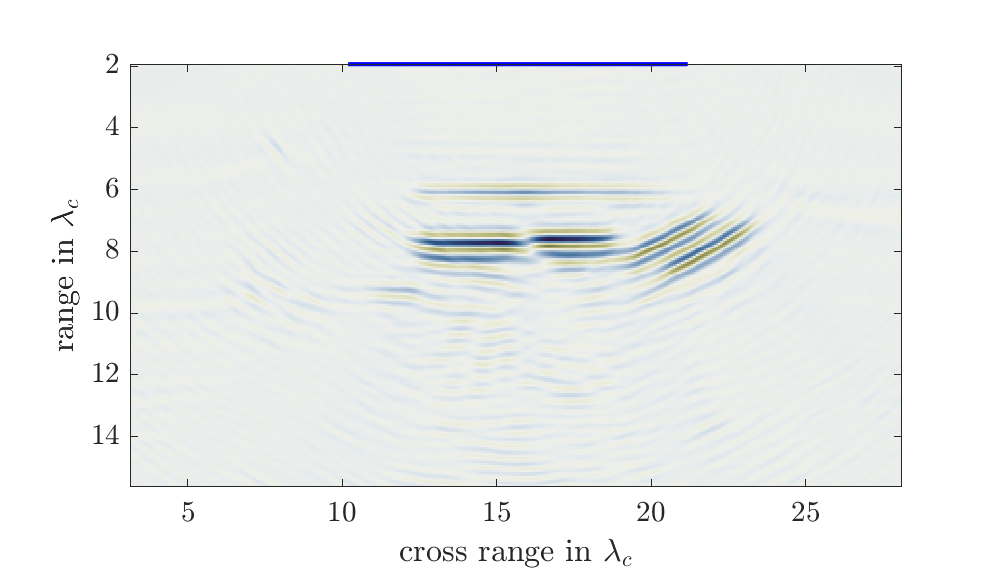}
\includegraphics[width=0.4\textwidth]{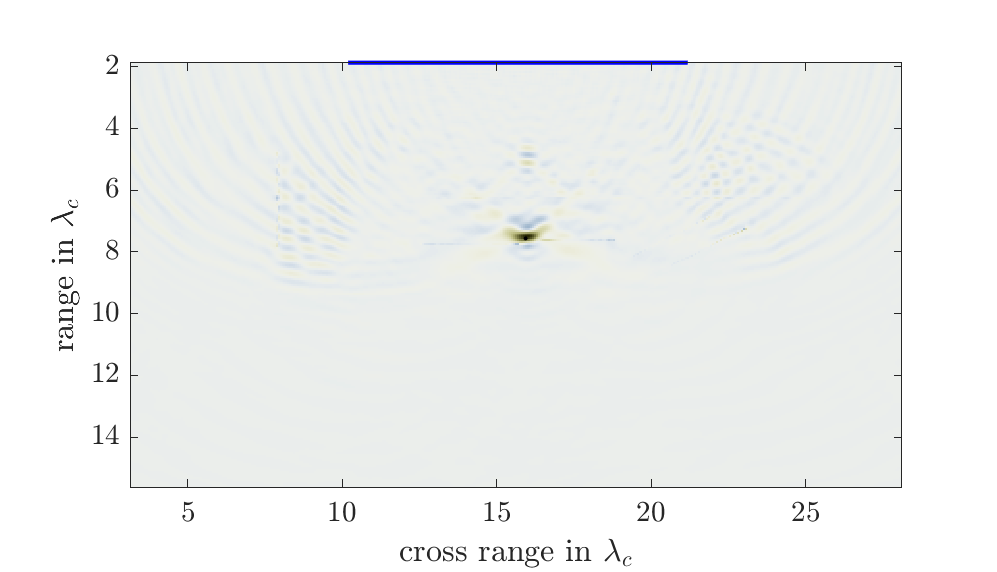}
\\
\includegraphics[width=0.4\textwidth]{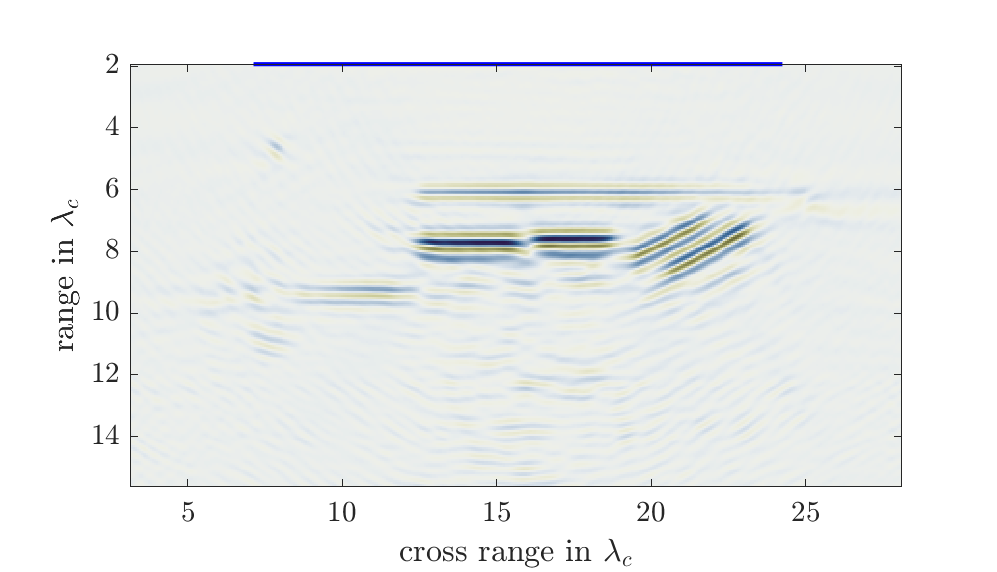}
\includegraphics[width=0.4\textwidth]{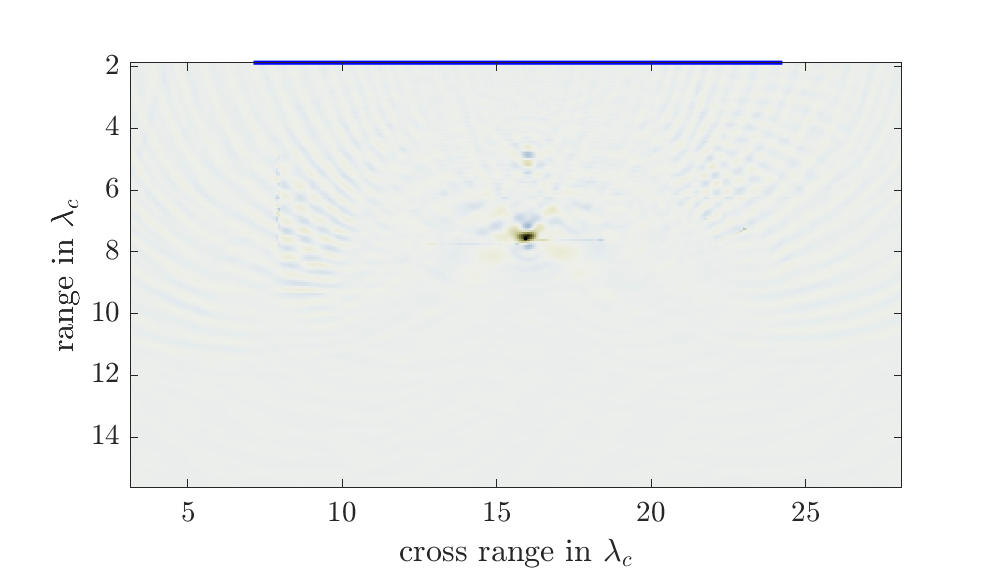}
\\
\includegraphics[width=0.4\textwidth]{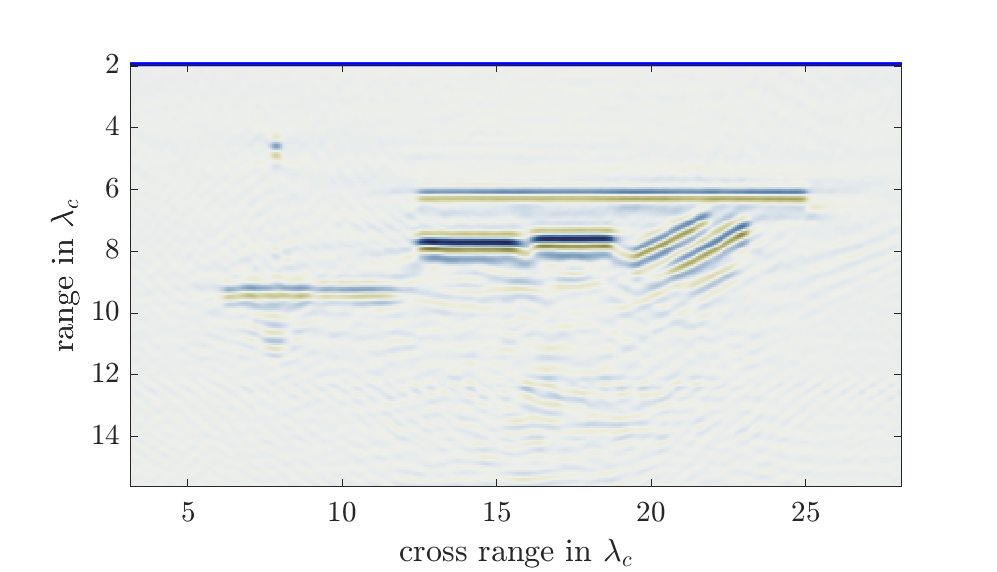} 
\includegraphics[width=0.4\textwidth]{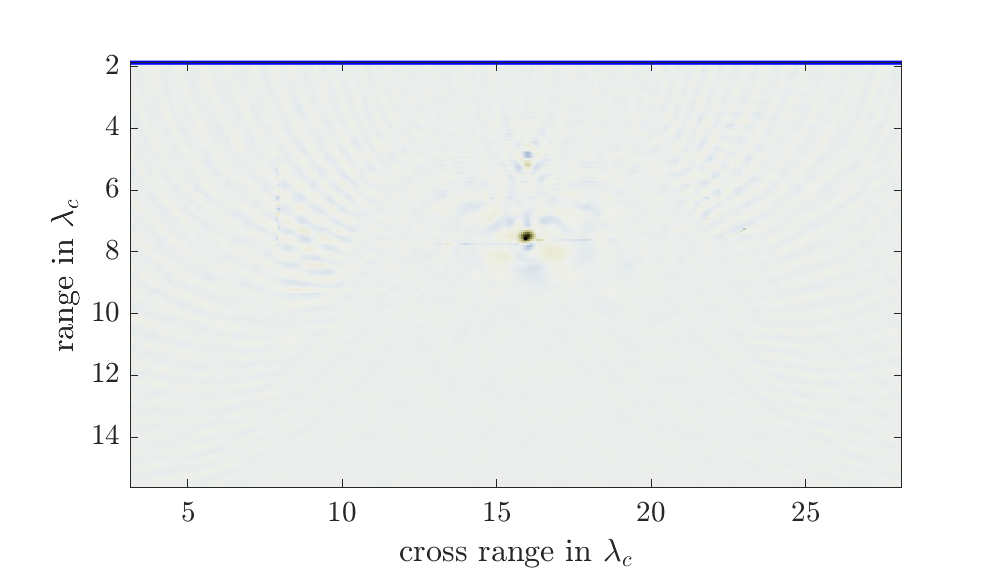}
\caption{Illustration of the effect of the aperture size. Left column: Range derivative of the imaging function $\cI(\by)$. Right column: The ROM point spread function $\delta_\by^{\rm ROM}(\bx)$ for the point $\by$ between the two nearby horizontal reflectors. The aperture of the array is shown in blue at the top of the plots. Top row for $40\%$ aperture, middle row for $60\%$ aperture and bottom row for 
the full aperture $a = 30 \la_c$ and $m = 49$ sensors. The separation between the sensors is kept the same, so the smaller the aperture, the fewer sensors. The time sample interval is $\tau =0.4 \pi/\om_c.$ }
\label{fig:Sim4}
\end{figure}

In Fig. \ref{fig:Sim4} we illustrate the effect of the aperture size on the ROM point spread function 
$\delta_\by^\RM$ and the image $\cI(\by)$. The larger the aperture, the better the focussing of $\delta_\by^\RM$ in cross-range and  the better the image.

\begin{figure}[htb]
\centering
\includegraphics[width=0.4\textwidth]{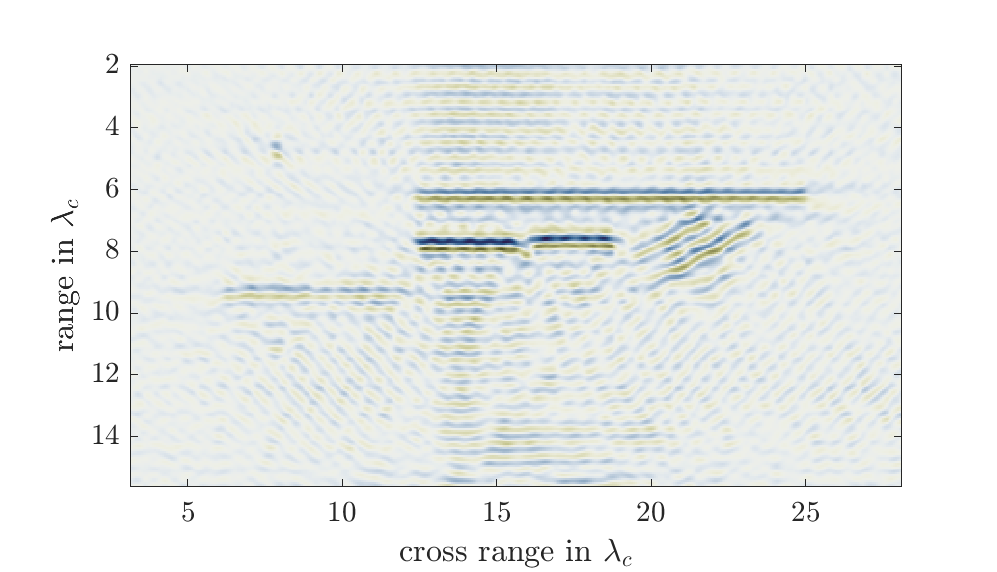}
\includegraphics[width=0.4\textwidth]{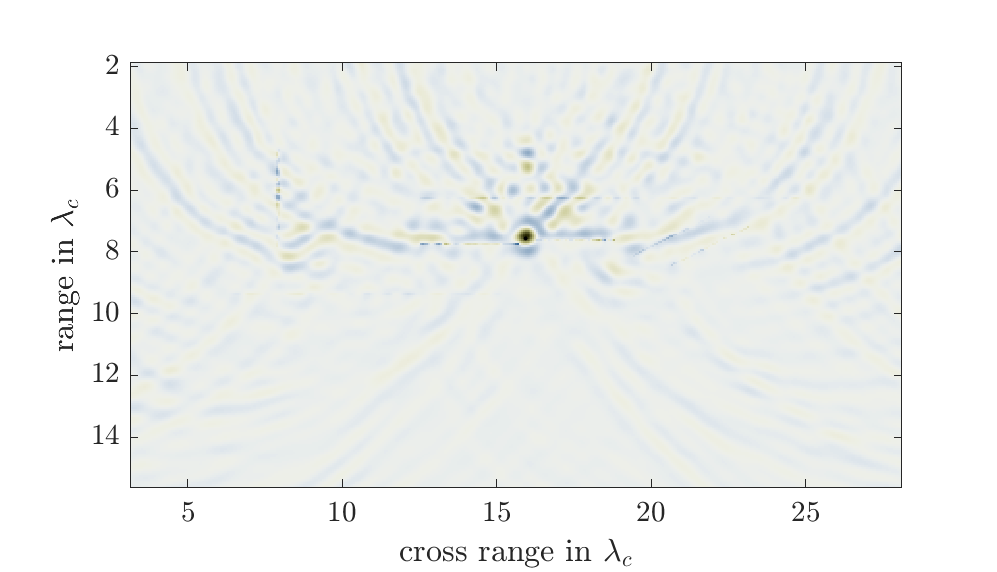}
\\
\includegraphics[width=0.4\textwidth]{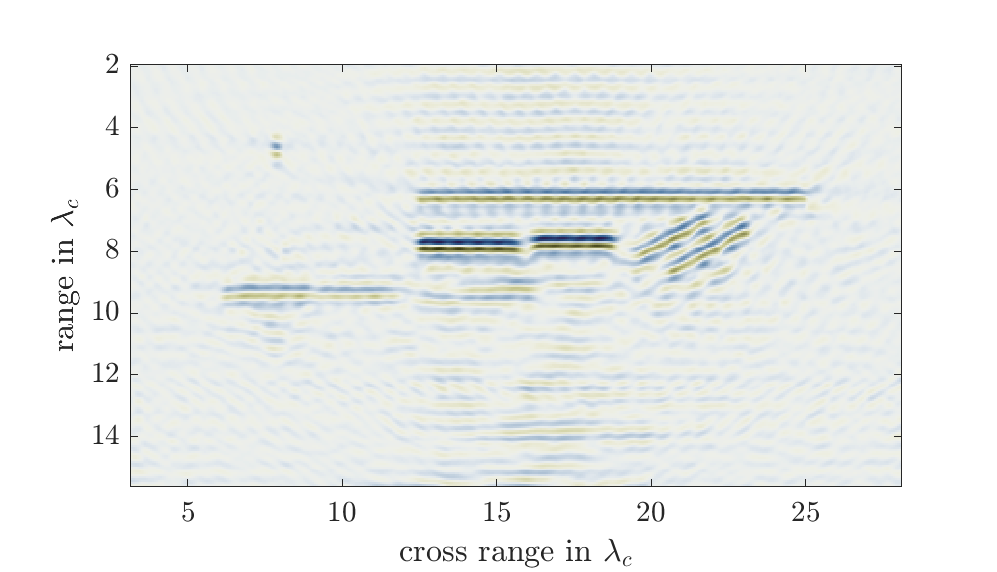}
\includegraphics[width=0.4\textwidth]{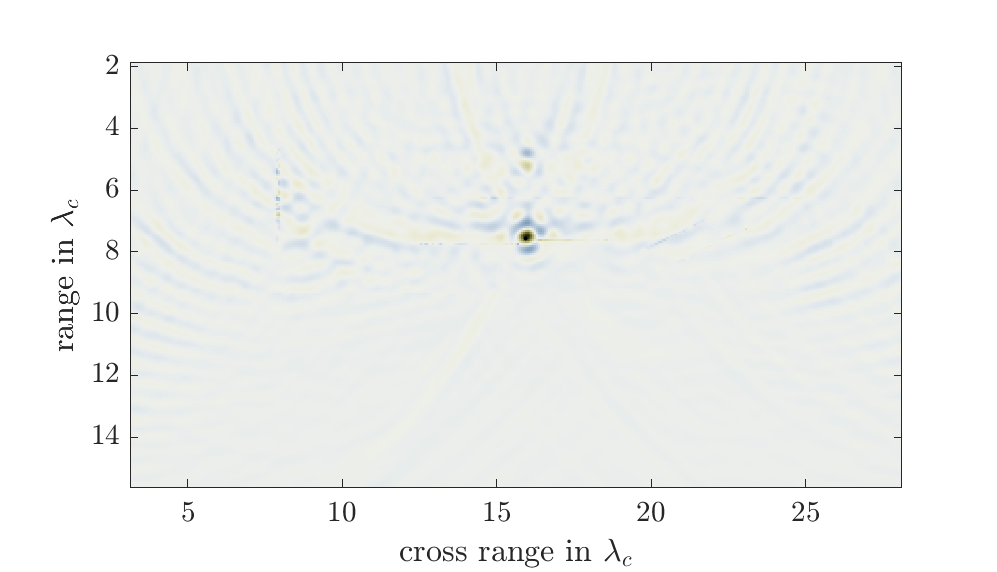}
\\
\includegraphics[width=0.4\textwidth]{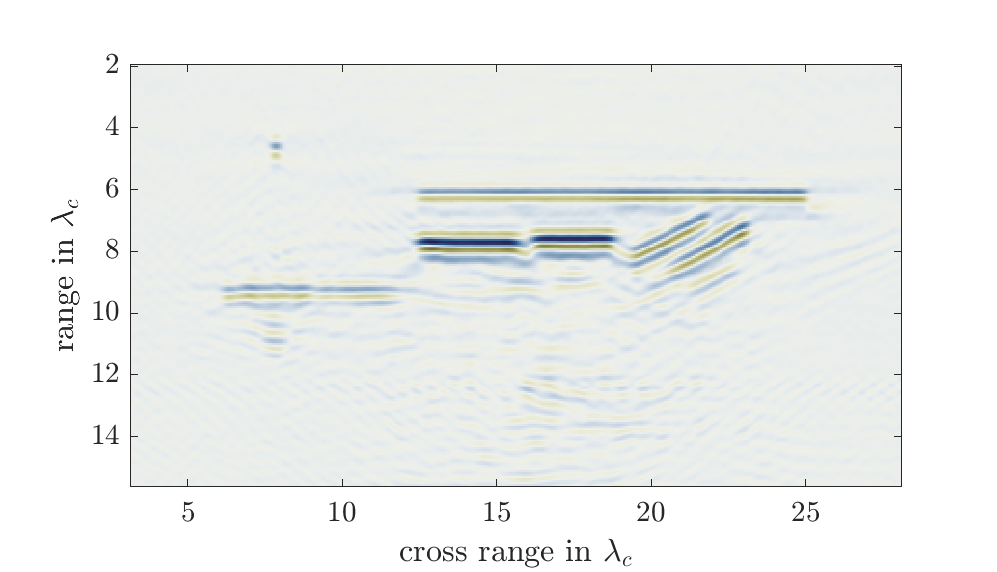}
\includegraphics[width=0.4\textwidth]{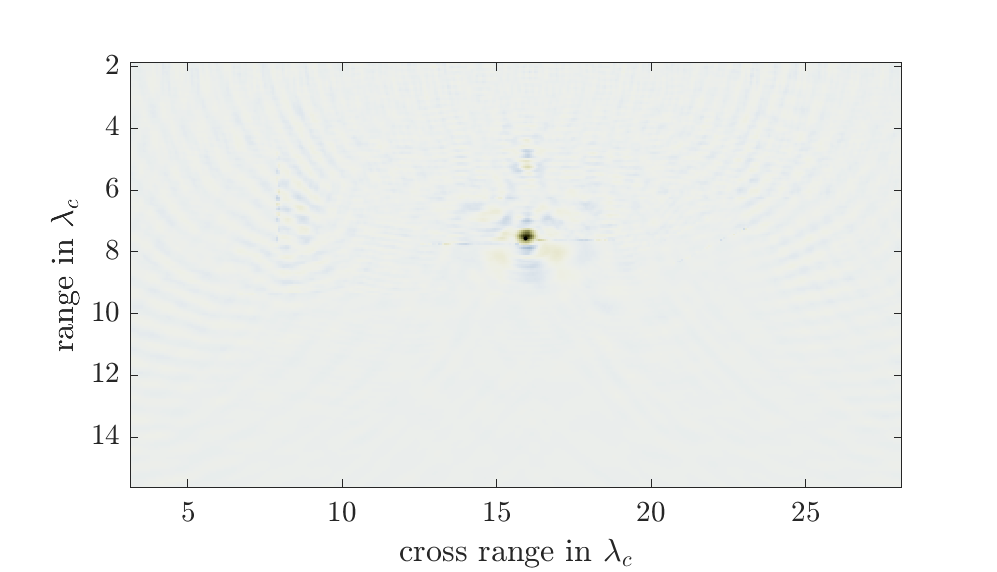}
\caption{Illustration of the effect of $\tau$.  Left column: Range derivative of the imaging function $\cI(\by)$. Right column: The ROM point spread function $\delta_\by^{\rm ROM}$ for the point $\by$ between the two nearby horizontal reflectors. The reference $\tau$ is $\tau_{\rm ref} = 0.4 \pi/\om_c$. Top row for $\tau = 3\tau_{\rm ref}$, middle row for $\tau = 1.8 \tau_{\rm ref}$ and bottom row for $\tau = 0.8 \tau_{\rm ref}$ (the case with $\tau=\tau_{\rm ref}$ is shown in the bottom row of Figure \ref{fig:Sim4}). The aperture is $a = 30 \la_c$, with $m = 49$ sensors. The duration of the experiment is kept the same, so the larger $\tau$, the fewer time steps.}
\label{fig:Sim5}
\end{figure}

Fig. \ref{fig:Sim5} shows the effect of the time sampling interval $\tau$. The reference value is 
as in the previous experiments $\tau = 0.4 \pi/\om_c$. For larger $\tau$ the focus of the ROM point spread function  deteriorates and the image becomes noisy.
For the smaller $\tau $ the results are basically the same as  in the bottom plots of Fig. \ref{fig:Sim4}. In practice $\tau$ should not be reduced too much, because the 
snapshots become too close to each other and consequently, the Cholesky factorization and the Gram-Schmidt orthogonalization become ill conditioned. 

\begin{figure}[htb]
\centering
\includegraphics[width=0.4\textwidth]{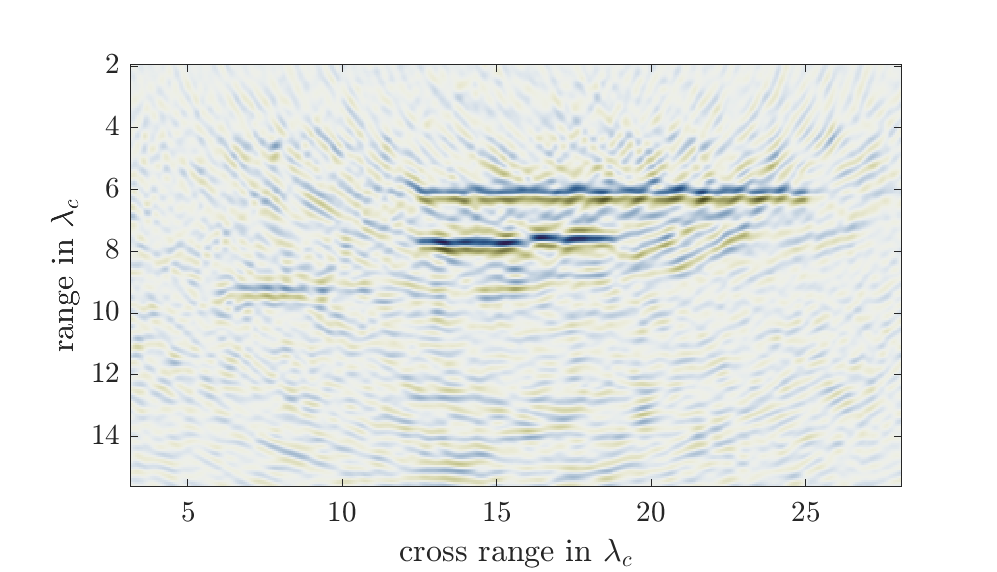}
\includegraphics[width=0.4\textwidth]{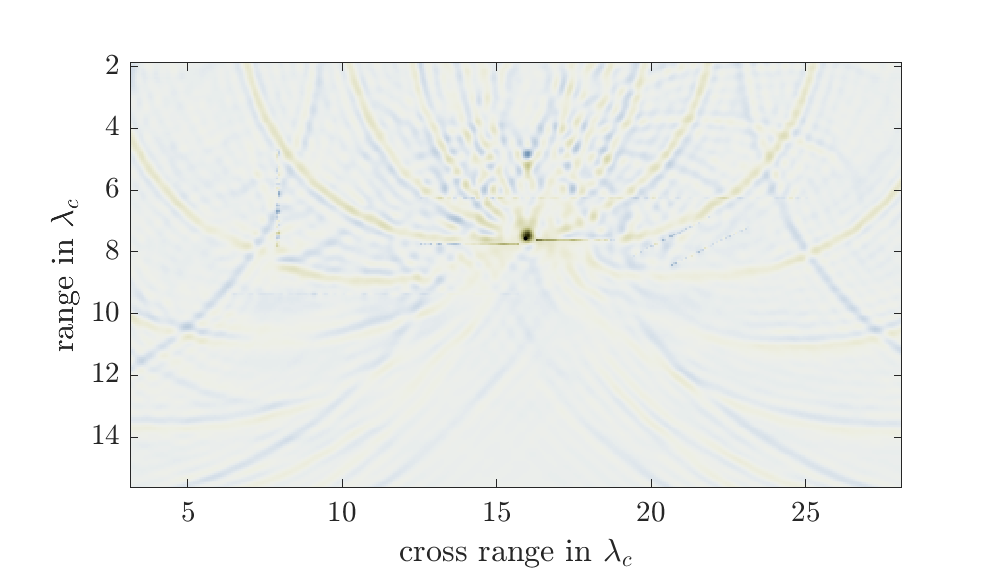}
\\
\includegraphics[width=0.4\textwidth]{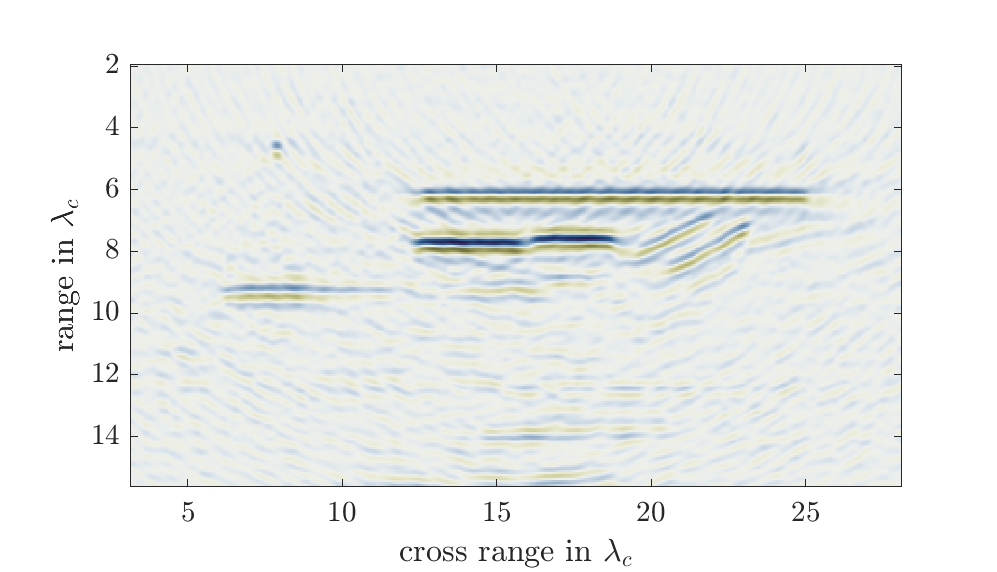}
\includegraphics[width=0.4\textwidth]{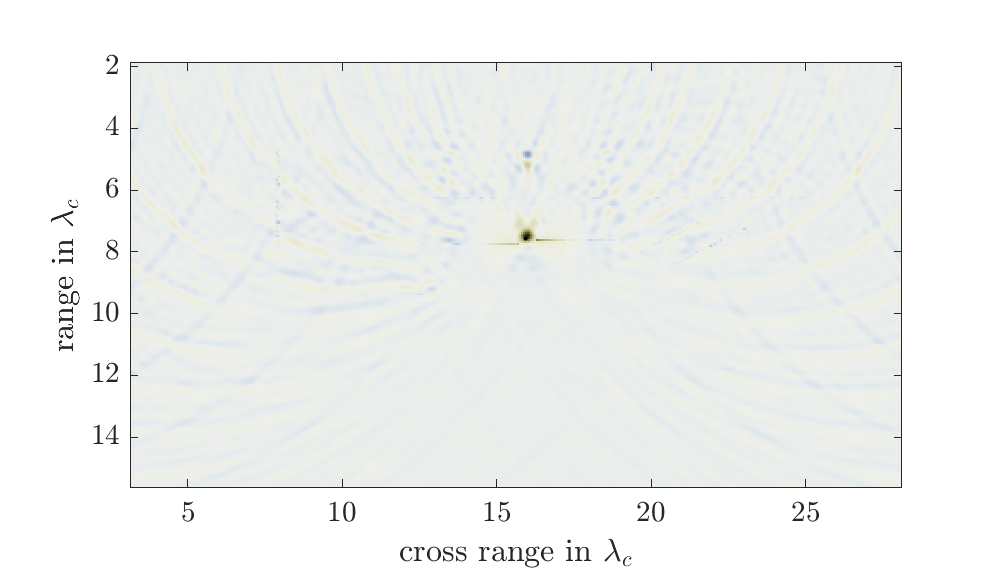}
\\
\includegraphics[width=0.4\textwidth]{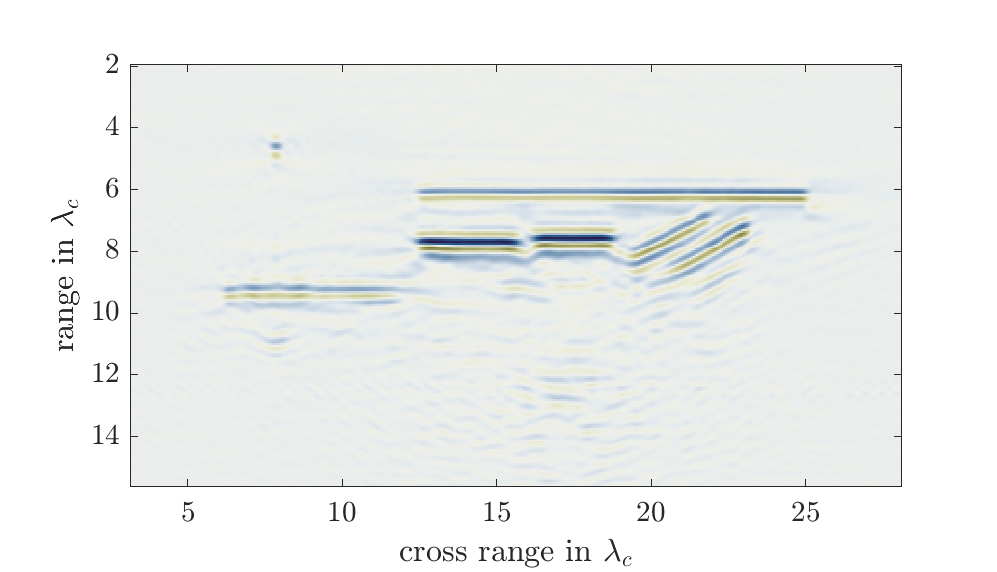}
\includegraphics[width=0.4\textwidth]{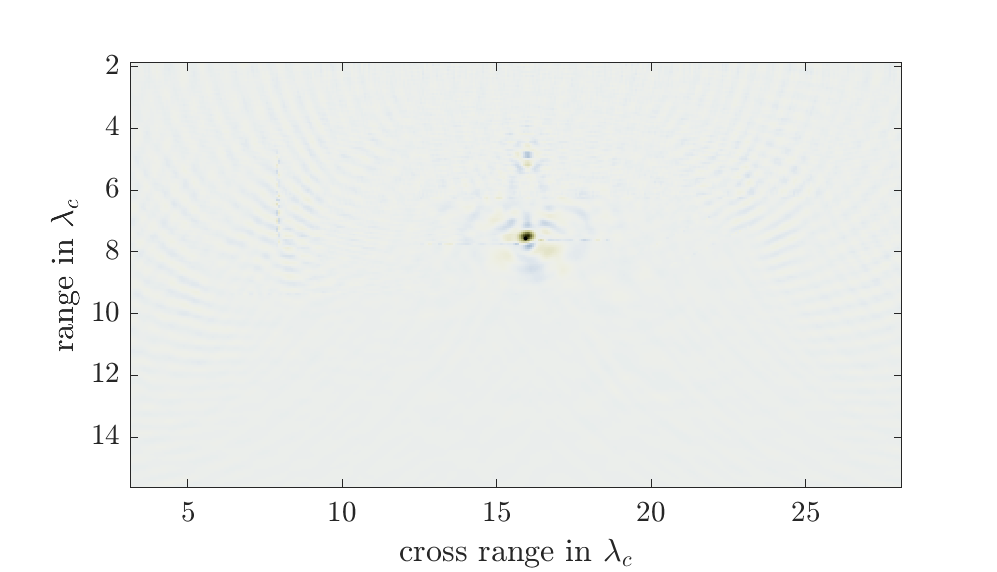}
\caption{Illustration of the effect of $m$.  Left column: Range derivative of the imaging function $\cI(\by)$. Right column: The ROM point spread function $\delta_\by^{\rm ROM}(\bx)$ for the point $\by$ between the two nearby horizontal reflectors. The reference $\tau$ is $\tau_{\rm ref} = 0.4 \pi/\om_c$. Top row for $m = 10$, middle row for $m = 20$ and bottom row for $m = 60$ (the case with $m=49$ is shown in the bottom row of Figure \ref{fig:Sim4}). The aperture is $a = 30 \la_c$, so the smaller $m$ is, the larger the separation between the sensors. The time sample interval is $\tau = 0.4 \pi/\om_c$.}
\label{fig:Sim6}
\end{figure}

Finally, we illustrate in Fig. \ref{fig:Sim6} the effect of the separation between the sensors. 
The aperture is fixed at $a = 30 \la_c$ and we display results for $m = 10, 20$ and $60$ equidistant sensors. We see that if the sensors are too far apart, the focus of the ROM point spread function  deteriorates and the image becomes noisy. The bottom plots obtained with $m = 60$ are basically the same as those for $m = 49$ (shown in the bottom row of Figure \ref{fig:Sim4}). In practice one should not take $m$ too large (i.e., sensors that are too close), because the Cholesky factorization of the mass matrix and the Gram-Schmidt orthogonalization become ill conditioned.

%{\color{blue}
%We display in Fig. ~\ref{fig:Delta1} the norm of $\|\delta_{\by}^\RM\|_{L^2(\Omega_{\rm im})}$ as a function of 
%$\by \in \Omega_{\rm im}$. Note that it is a smooth function that decreases slowly with range, due to the finite duration of the measurements, and is thus not affected by the sharp variations of the wave speed which model the reflectors in the host medium. }
% 
%\begin{figure}[t]
%\centering
%\includegraphics[width=0.4\textwidth]{./Figures/horizontal_inclusions/PixelScan_deltaNormobtnotitle.png}
%\vspace{-0.2in}\caption{The $L^2$ norm of the ROM point spread function $\delta_{\by}^{\mbox{ROM}}$ as a function of 
%$\by \in \Omega_{\rm im}$. The time sampling interval is $\tau = 0.4 \pi/\om_c$, the aperture length is $a = 
%30 \la_c$ and the array has $m = 49$ sensors. {\color{blue} Should we display a colorbar here?}}
%\label{fig:Delta1}
%\end{figure}

\subsection{Focussing with the internal wave}
\label{sect:num3}
In this section  we illustrate the focusing of the wave $\gamma(t,\bx;\by)$, the solution of 
\eqref{eq:PS3}-\eqref{eq:PS6} with the illumination \eqref{eq:PS1} defined in terms of the internal wave computed  as in Proposition \ref{prop.1}. The setup is as in Fig. \ref{fig:Sim1} and we use 
the large aperture $a = 30 \la$, with $m = 49$ sensors and the time sample $\tau = 0.4 \pi/\om_c$. 

\begin{figure}[htb]
\centering
\includegraphics[width=0.31\textwidth]{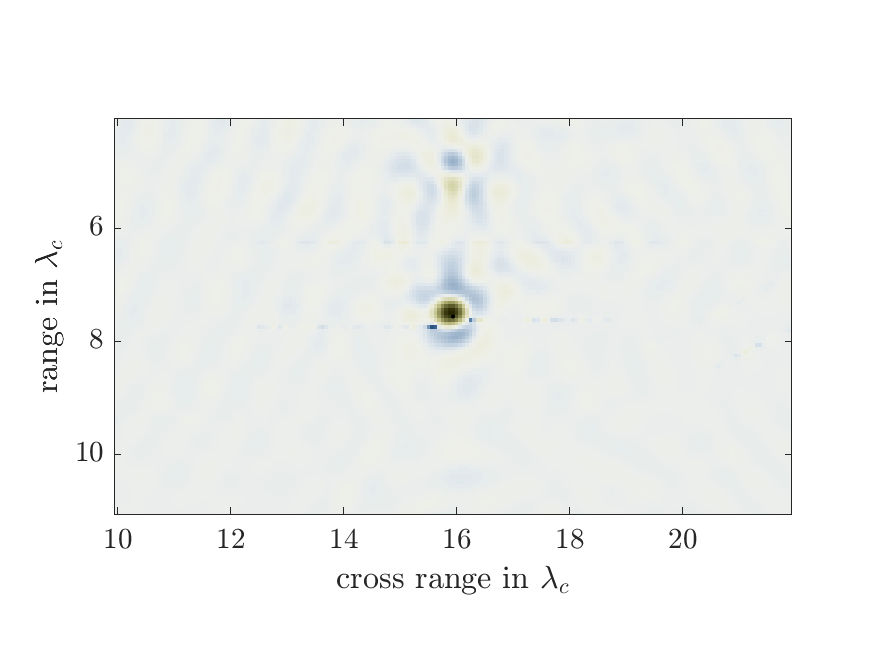}
\includegraphics[width=0.31\textwidth]{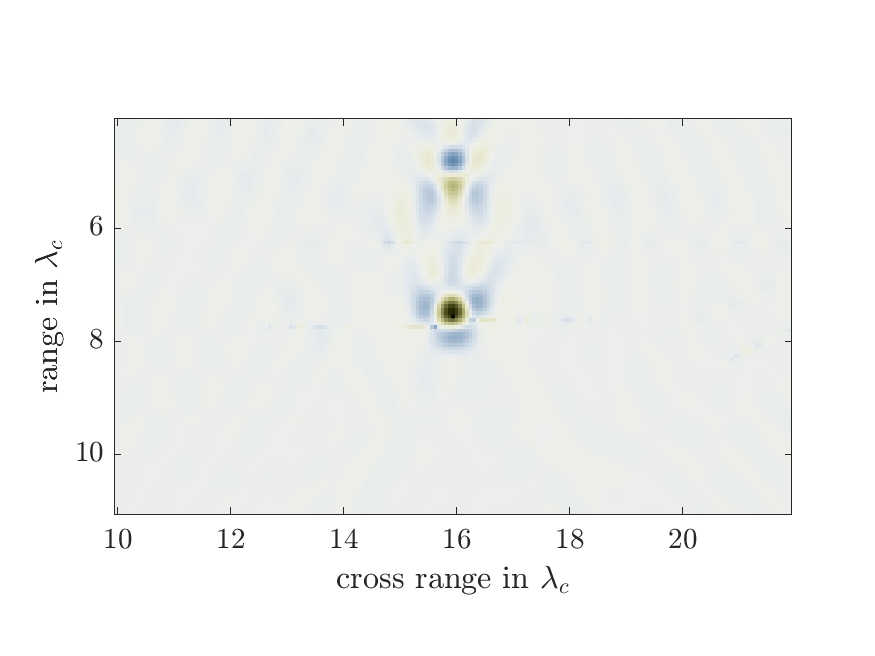}
\includegraphics[width=0.31\textwidth]{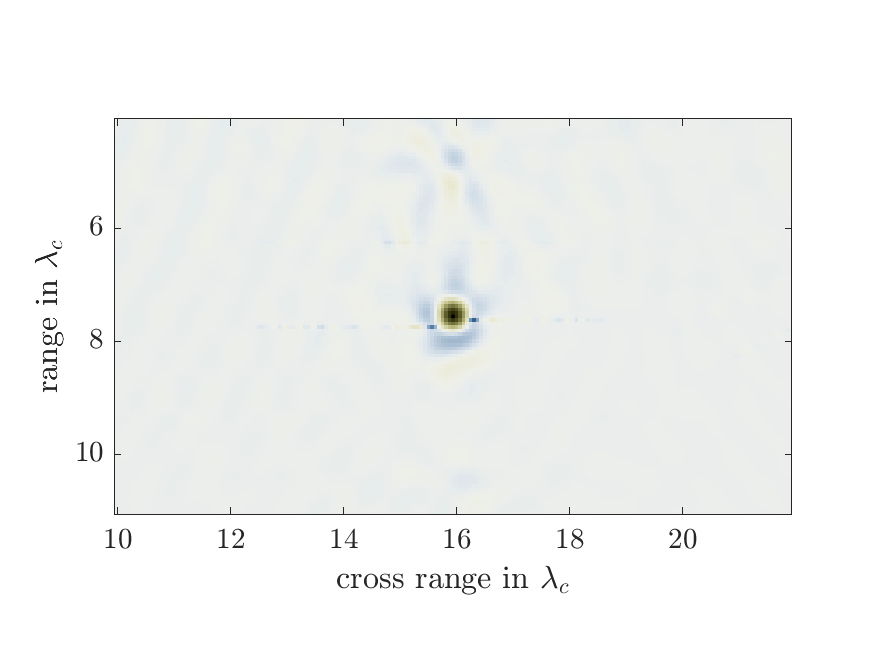}
\\
\includegraphics[width=0.31\textwidth]{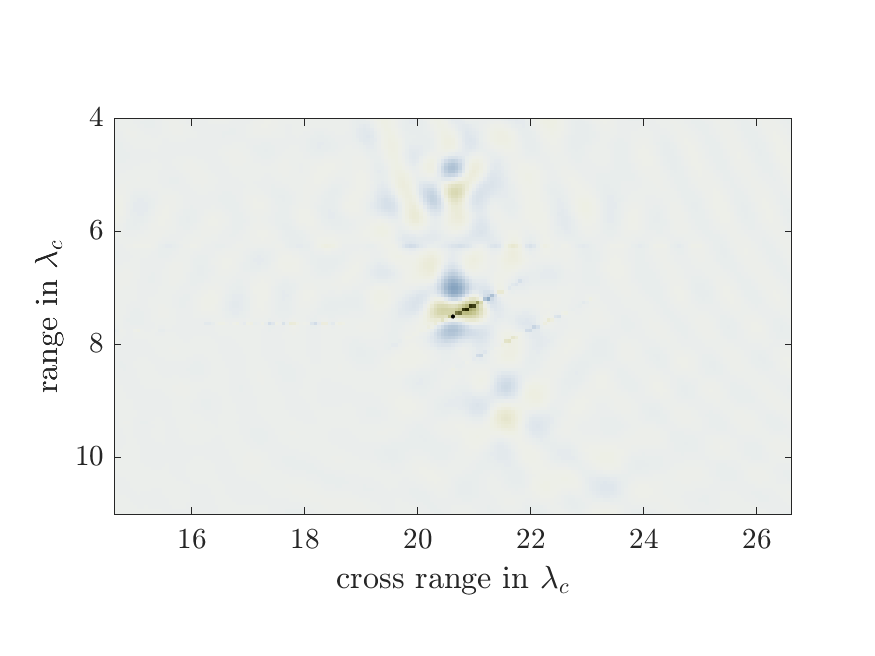}
\includegraphics[width=0.31\textwidth]{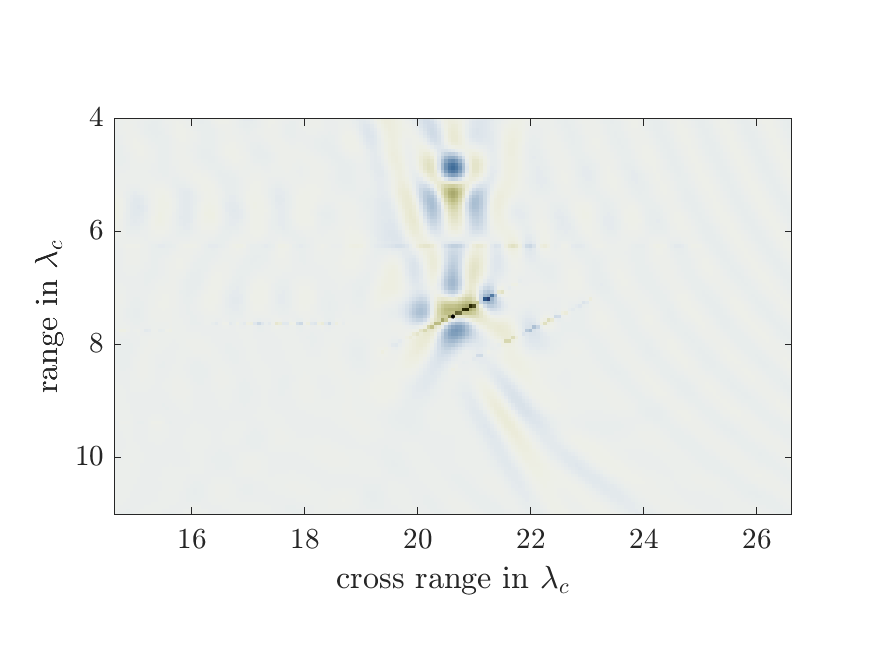}
\includegraphics[width=0.31\textwidth]{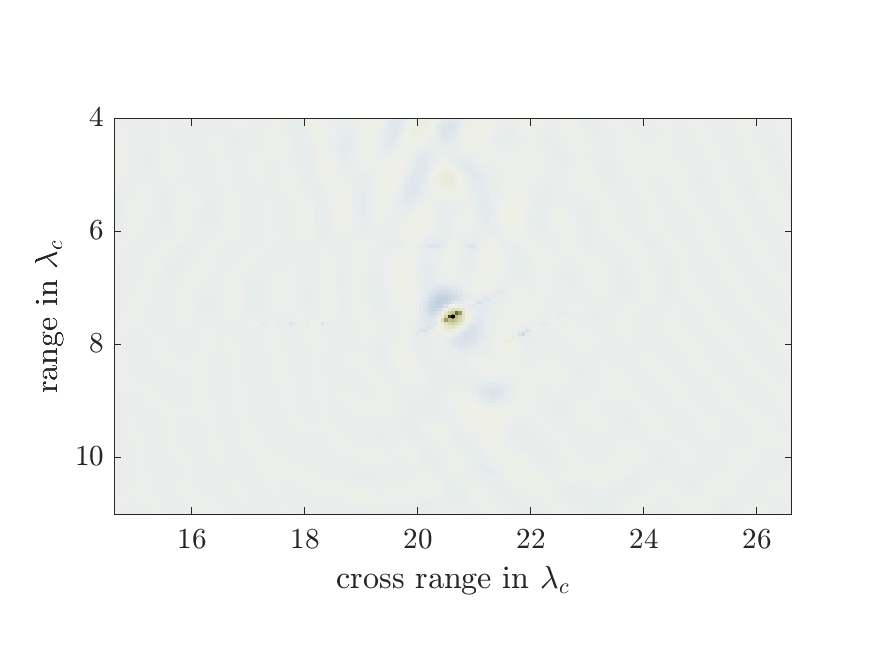}
\\
\includegraphics[width=0.31\textwidth]{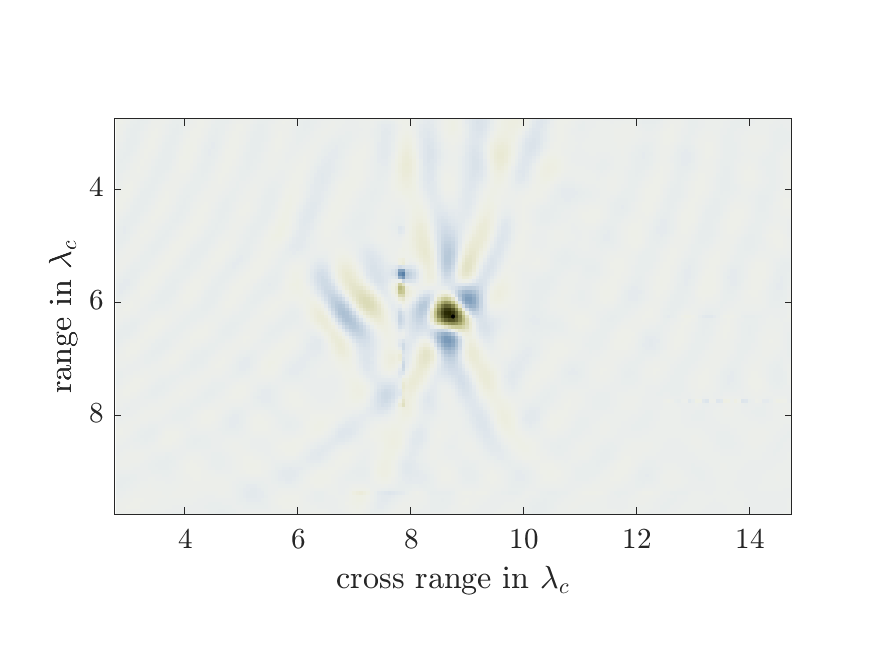}
\includegraphics[width=0.31\textwidth]{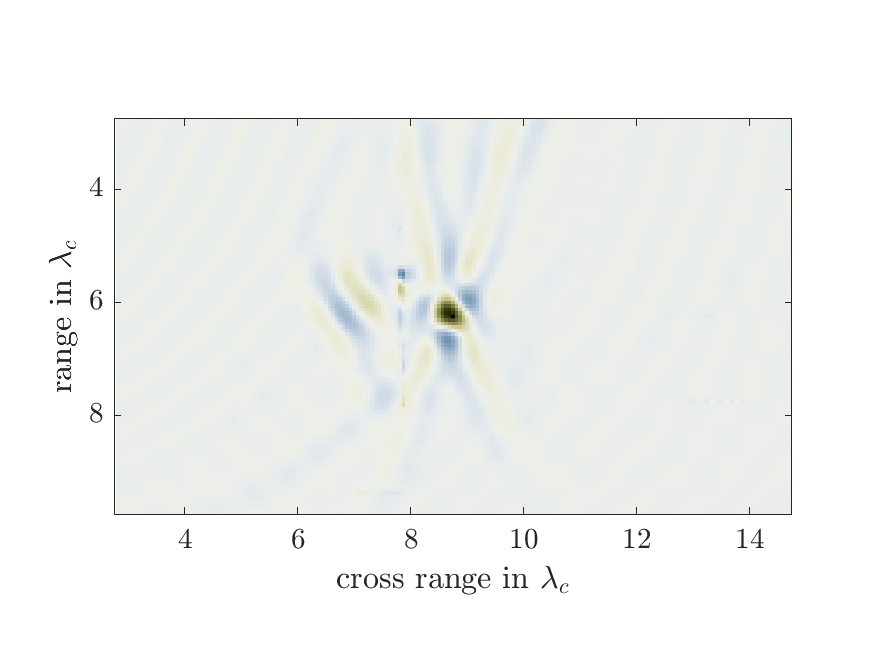}
\includegraphics[width=0.31\textwidth]{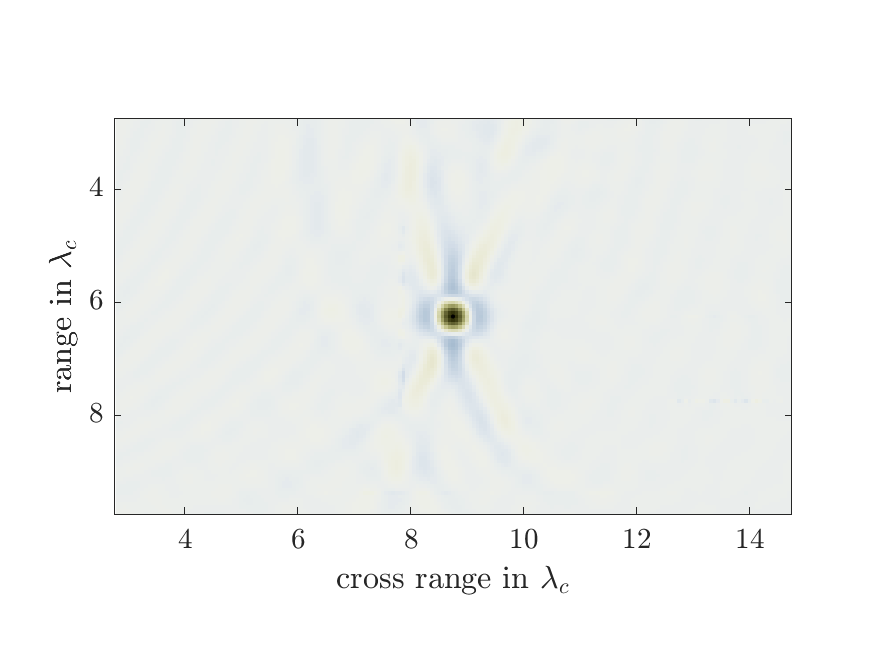}
\vspace{-0.15in}
\caption{Illustration of focusing at three different points in the imaging domain. Left column: 
The refocused wave after the illumination computed with the internal wave $g(t,\bx;\by)$. 
Middle column: The refocused wave after the illumination computed with $g_o(t,\bx;\by)$ calculated in the reference medium. Right column: The refocused wave after the illumination computed with  $g^{\rm ideal}(t,\bx;\by)$.
The aperture is $a = 30 \la_c$,  the number of sensors is $m = 49$ and  $\tau = 0.4 \pi/\om_c$.}
\label{fig:Sim7}
\end{figure}

We display in the left column of Fig. \ref{fig:Sim7} the wave $\gamma(t,\bx;\by)$  at the time of focus, 
for three different points $\by \in \Omega_{\rm im}$: Between the two nearby horizontal reflectors, 
on one of the oblique reflectors and near the hard to see vertical reflector. For comparison, we also display the waves given by the illumination calculated as in \eqref{eq:PS1}, with $g(t,\bx;\by)$ 
replaced by:  $g_o(t,\bx;\by)$ calculated in the reference medium (middle column plots) and by  $g^{\rm ideal}(t,\bx;\by)$ that cannot be computed in practice (right column plots). While the refocusing is not as 
good as the unattainable one obtained with $g^{\rm ideal}(t,\bx;\by)$, we see that using  $g(t,\bx;\by)$ 
 is better than $g_o(t,\bx;\by)$ for the two first points $\by$. 

\subsection{Simulations with noisy data}
\label{sect:num4}
\begin{figure}[htb]
\centering
\includegraphics[width=0.4\textwidth]{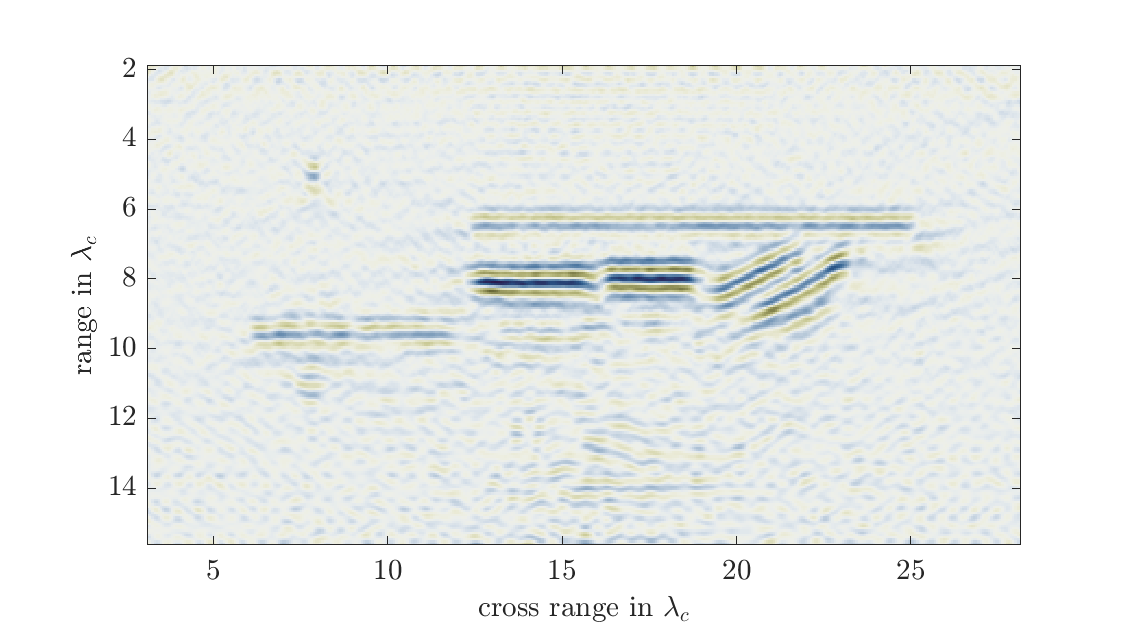}
\includegraphics[width=0.4\textwidth]{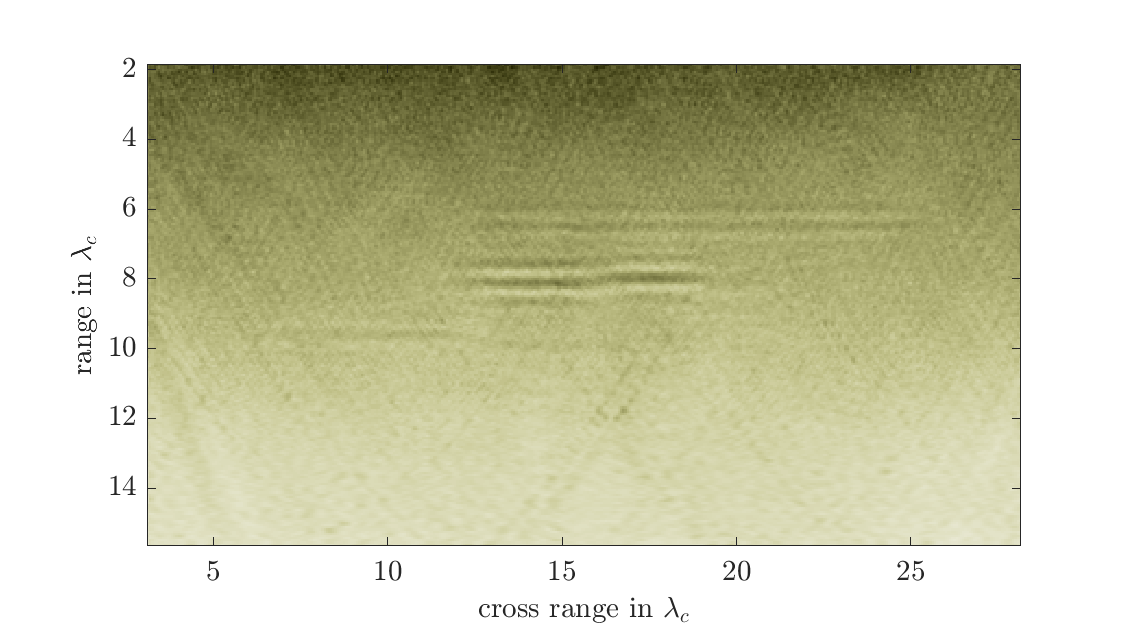}
\vspace{-0.1in}
\caption{The range derivative imaging function $\cI(\by)$ (left) and the  backprojection image $\cI^{\rm BP}(\by)$ (right) when imaging with data that has $20\%$ additive white Gaussian noise. The aperture is $a = 30 \la_c$,  the number of sensors is $m = 49$ and $\tau =0.67 \pi/\om_c.$ }
\label{fig:Sim10}
\end{figure}

We now show results for the same waveguide setting as above, and for data contaminated with white Gaussian, additive noise. The noisy data set is $ \left\{w^\ss(j\tau,\bx_r)+\varepsilon_{s,r,j}, ~s,r=1,\ldots, m, ~~ j =0,\ldots, 2n-1\right\}$,
where the $\varepsilon_{s,r,j}$ are independent and identically distributed Gaussian random variables with zero-mean and variance $\sigma^2_{\rm noise} = 0.2^2 \max_{s,r,j}
\{ w^\ss(j\tau,\bx_r)^2 \}$.
As we explained above, the conditioning of the noiseless mass matrix depends on the time sampling interval $\tau$. The smaller $\tau$ is, 
the worse the conditioning. Thus,  for noisy data it is beneficial to increase $\tau$ a little, as we do in the results in Fig. \ref{fig:Sim10}. 
However, even for such $\tau$, the mass matrix needs to be regularized, as explained at the beginning of the section.

In Fig.~\ref{fig:Sim10} we show images obtained at $20\%$  noise level. 
Because the computation of the  backprojection image $\cI^{\rm BP}(\by)$ involves the unstable step of inverting the Cholesky factor $\bR$ of the mass matrix, the noise effect  is much worse than in $\cI(\by)$ when using the same regularization \eqref{eq:regulariz}. A more involved stabilization of the backprojection image is needed, as explained in \cite{druskin2018nonlinear,borcea2019robust}. We do not repeat that regularization strategy here, but note that it typically leads to ghost multiples at high noise levels.

\subsection{Imaging in the half space}
\label{sect:num2}
\begin{figure}[htb]
\centering
\raisebox{0.3in}{\includegraphics[width=0.45\textwidth]{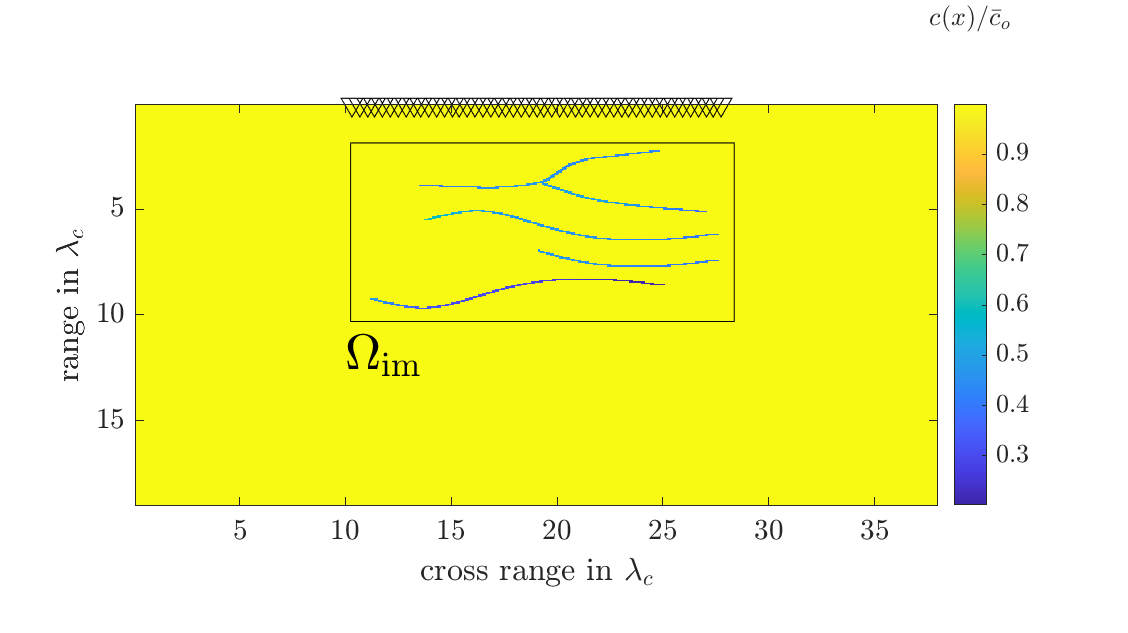}}
\includegraphics[width=0.54\textwidth]{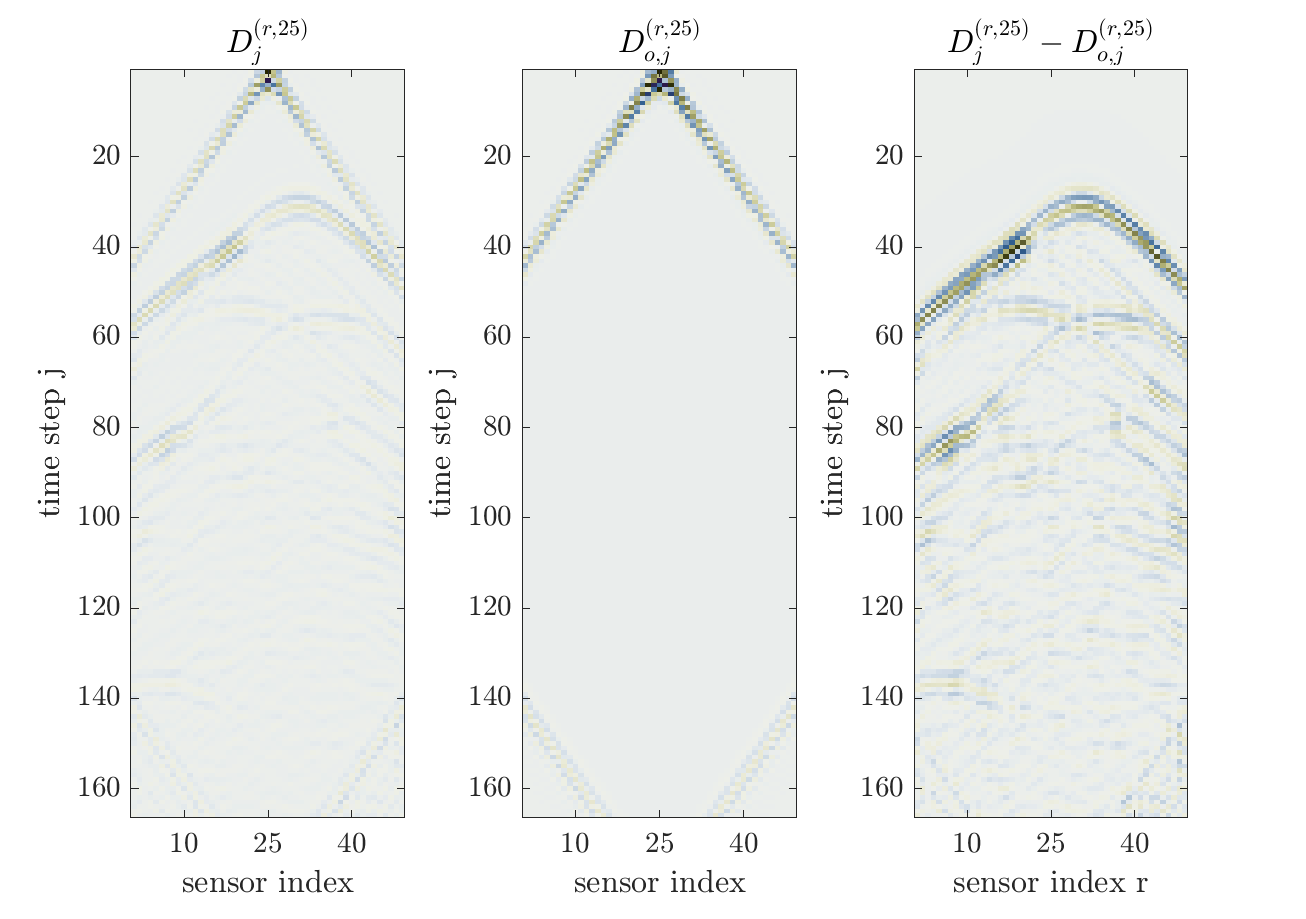}
\vspace{-0.25in}\caption{Left: Illustration of the setup: The array of $m = 49$ sensors (indicated with triangles) lying near the accessible boundary probes a medium with wave speed $\bar c_o$, containing a few thin 
reflecting structures (cracks), modeled by the low velocity shown in the color bar. Right: The data corresponding to the illumination from the center element in the array. We show it for the medium with the reflectors, the reference medium and the difference between the two.
}
\label{fig:Sim8}
\end{figure}

Here we present numerical results for the setup shown in Fig. \ref{fig:Sim8}, where the side boundaries are sufficiently far to have no effect on the data displayed in the right plots. 

\begin{figure}[htb]
\centering
	\includegraphics[width=0.4\textwidth]{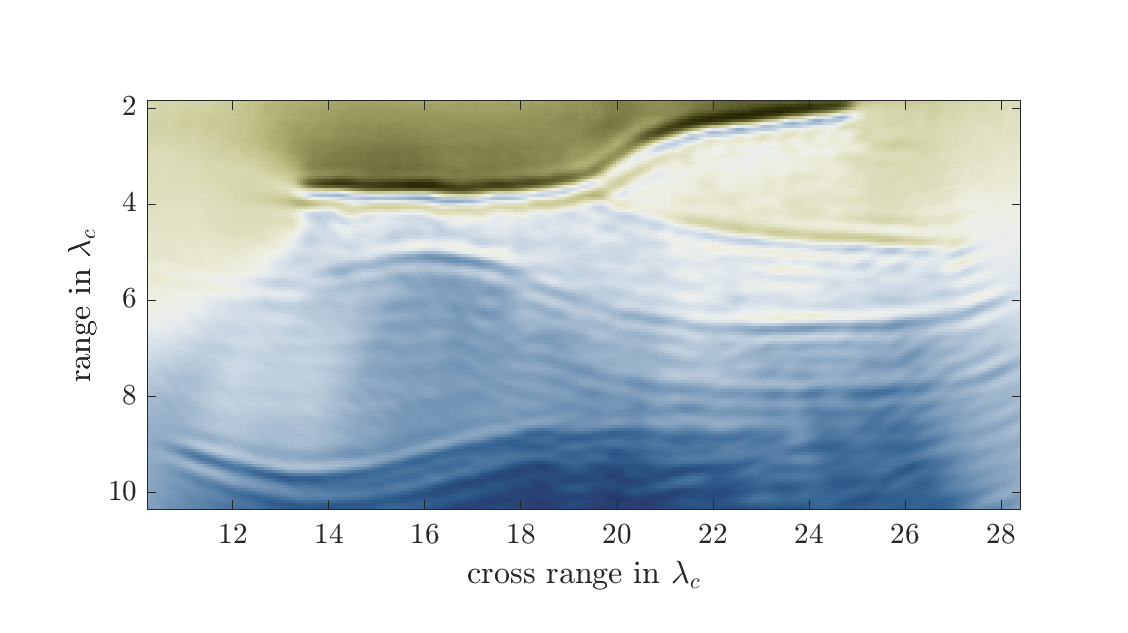}
\includegraphics[width=0.4\textwidth]{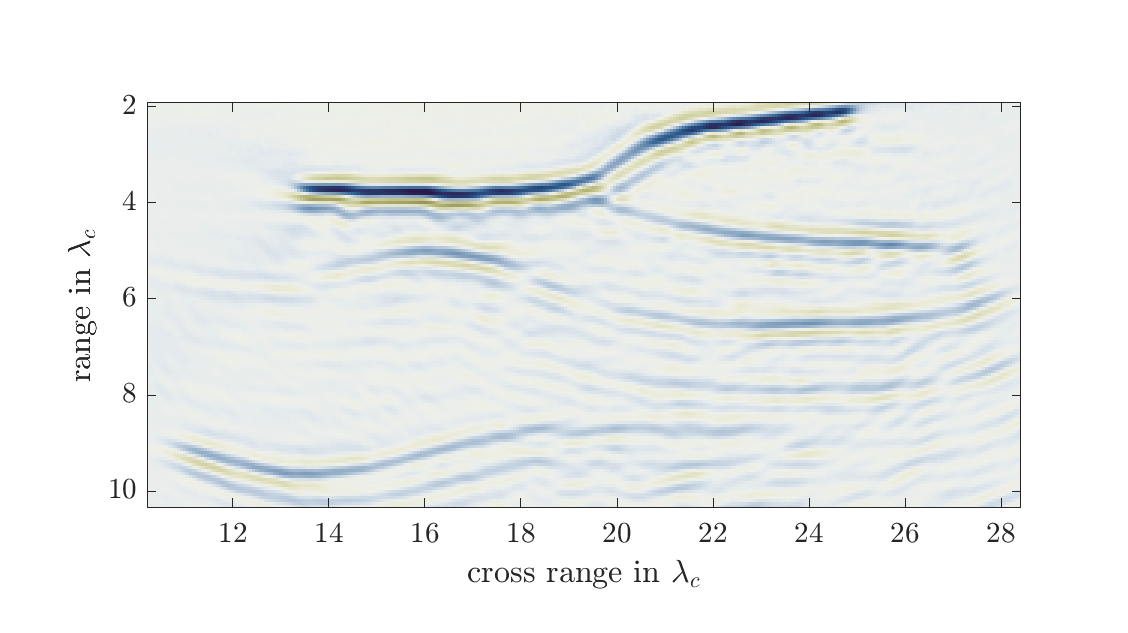}\\
\includegraphics[width=0.4\textwidth]{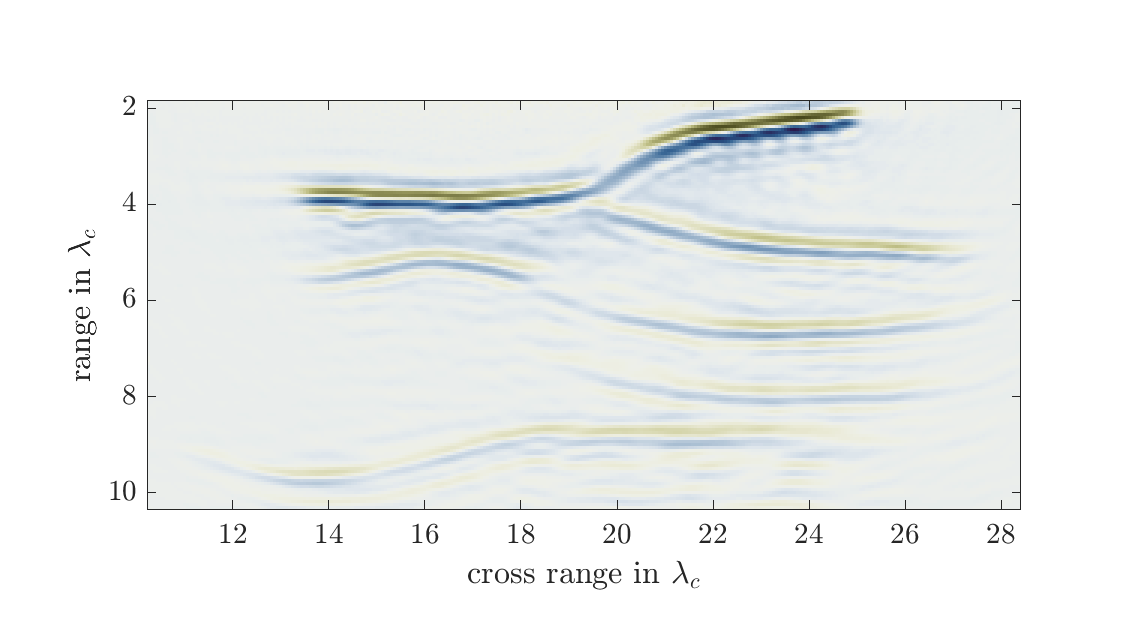}
\includegraphics[width=0.4\textwidth]{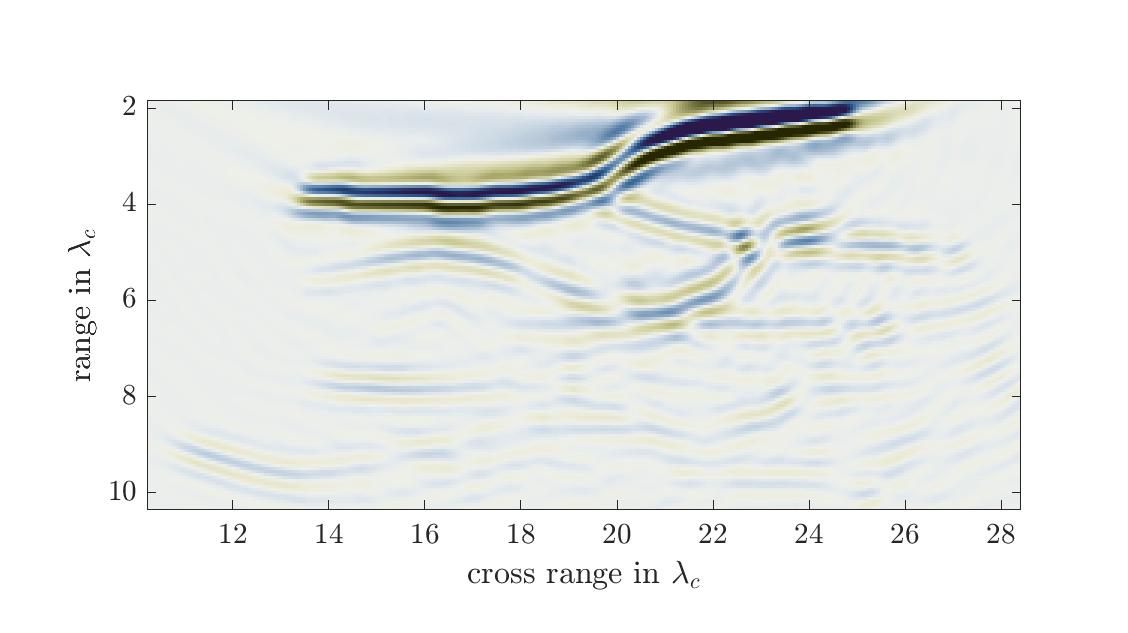}\\
\includegraphics[width=0.4\textwidth]{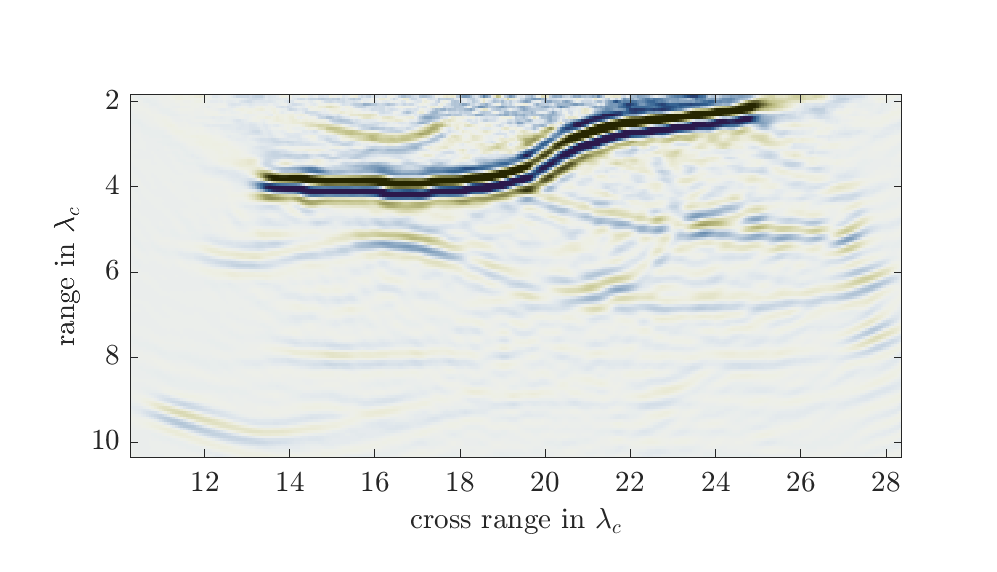}

\vspace{-0.1in}
\caption{Top plots: The imaging function $\cI(\by)$ (left) and its range derivative (right). Middle plots: The  backprojection image $\cI^{\rm BP}(\by)$ (left) and the reversed time migration image  $\cI^{\small {\rm RTM}}(\by)$ (right) for the setup shown in Fig. \ref{fig:Sim8}. Bottom plot: The range derivative of the pixel scanning image  $\cI^{\rm PS}(\by)$. The aperture length is $a = 18 \la_c$ and the array has $m = 49$ sensors. The time sample interval is $\tau =0.42 \pi/\om_c.$ }
\label{fig:Sim9}
\end{figure}

We show in Fig. \ref{fig:Sim9} images obtained with the  aperture size $a = 18 \la_c$,  
containing $m = 49$ equidistantly spaced sensors. The time sampling interval is $\tau = 0.42 \pi/\om_c$. Note that the  multiple scattering artifacts  in the reversed time migration and the pixel scanning images are more pronounced than in Fig. \ref{fig:Sim3}. The backprojection image $\cI^\BP(\by)$ and $\cI(\by)$ do not have such artifacts and they both localize well the crack-like reflectors. Arguably, $\cI(\by)$ does a slightly better job at localizing the sloped part of the middle crack.

%----------------------------------
\newpage

\section{Summary}
\label{sect:sum}

We introduced and studied with analysis and numerical simulations a novel, computationally inexpensive approach for imaging reflectors in a host, non-scattering medium, with 
an active array of $m$ sensors, which probe the medium with a pulse $f(t)$  and measure the generated waves. 
The measurements are for a finite duration $2n \tau$, at instants spaced by $\tau$, chosen to satisfy the Nyquist sampling 
requirement for $f(t)$.  The imaging is based on a data driven reduced order model (ROM) of the wave propagator, the operator that maps the wave from one instant to the next. Specifically, it uses the ROM 
to estimate an ``internal wave" $g(t,\by,\bx_r)$ originating from the vicinity of the imaging point $\by$ and propagating through the unknown medium to the sensors at $\bx_r$, for $r = 1,\ldots, m$. 

We introduced two kinds of imaging functions: The first, denoted by $\cI(\by)$, has a very simple expression,
given by the squared norm of the internal wave at the sensors. The second, denoted by $\cI^\PS(\by)$,  can be implemented experimentally. It is a pixel scanning imaging approach which uses the internal wave to define a control of the illumination of the medium from the  array, for improved focusing of a probing wave at the  pixel (imaging point) $\by$. It then uses a matched field approach to obtain  $\cI^\PS(\by)$ from the resulting measured backscattered wave.
 
The functions $\cI(\by)$ and $\cI^\PS(\by)$ use a different imaging principle: The first one is designed to be sensitive to
variations of the wave speed locally, near the imaging point, so it is not affected by the arrivals of the multiply scattered 
echoes in the medium. The second one matches time arrivals at the array and is therefore affected by multiple scattering.

Both imaging functions use the time reversal refocusing principle. In particular, $\cI(\by)$ is a blurry version 
of the function that models the refocusing of the wave in the time reversal experiment with a source at the imaging point $\by$. 
The sharpness of the refocusing depends on the bandwidth of the probing pulse $f(t)$, 
the duration $2 n \tau$ of the measurements and the aperture size of the array.  These affect the resolution of $\cI(\by)$, 
but the blur, quantified by the ROM  point spread function $\delta_\by^\RM$, is the main factor. The better this peaks at 
the imaging point $\by$, the better the image.  We showed with analysis and numerical simulations that $\delta_\by^\RM$ is highly peaked at $\by$ if: (1) the kinematics (the smooth wave speed in the host medium) is known accurately;  (2) the time sampling interval $\tau$ and the sensor separation are small enough; and (3) the array has large enough aperture.  Of all these requirements, knowing the kinematics 
may be harder to achieve in some applications.

Since  the imaging function $\cI(\by)$ is easily computed and it is unfocused when the assumed kinematics is wrong, it could be possible to carry out an estimation of the smooth part of the wave speed based on an optimization of the sharpness of $\cI(\by)$ quantified properly by some norm. This could be the subject of future work. 

\section*{Acknowledgements}
This research is supported in part by the  ONR award N00014-21-1-2370 and by the AFOSR award
FA9550-21-1-0166.

\appendix
\section{The approximation of the snapshots in the reference space}
\label{ap:A}
In this appendix we discuss two setups where we can analyze explicitly the 
approximation \eqref{eq:approxV} of the orthonormal snapshots. Specifically, 
we quantify how well  the snapshots $u_j^\ss(\bx)$,
 for $j = 0, \ldots, n-1$ and $s = 1, \ldots, m$, which span the projection space $\mathscr{S}$, can be 
 approximated in the reference space $\mathscr{S}_o$. If the approximation error  is small, 
since the Gram Schmidt orthogonalization 
 is a stable process, we have $\bV(\bx) \approx \bV_o(\bx)$.

The first setup is for a layered medium and it is discussed in \ref{ap:A1}. The second setup, discussed in \ref{ap:A2}, is for a waveguide. 

\subsection{Snapshots in a layered medium}
\label{ap:A1}
The analysis is simplest in the one-dimensional case, where $m =1$, so we consider it first. 
The higher dimensional case is discussed after that.

\subsubsection{One-dimensional case.}
It is well known \cite[Chapter 3]{fouque2007wave} that in one-dimension scattering occurs due to changes of the acoustic impedance, and that the wave speed $c(z)$ can be eliminated from the wave equation by transforming  to the travel time coordinate
\begin{equation}
T(z):= \int_0^z \frac{dz'}{c(z')}.
\label{eq:A1}
\end{equation}
Thus, we consider, only in this section, the more general acoustic wave equation corresponding to variable mass density $\rho(z)$ and bulk modulus $K(z)$, which define the wave speed
$c(z) = \sqrt{K(z)/\rho(z)}$ and acoustic impedance $\zeta(z) = \sqrt{K(z) \rho(z)}$. The wave is modeled by the acoustic pressure $p(t,z)$, the solution of 
\begin{align}
\partial_t^2 p(t,z) - \zeta(z) c(z) \partial_z \left[ \frac{c(z)}{\zeta(z)} \partial_z p(t,z) \right] &= f'(t) c(0) 
\delta(z), \qquad t \in \RR, ~~ z \in (0-,L), 
\label{eq:A2} \\
\partial_z p(t,0-) &= p(t,L) = 0, \qquad t \in \RR, \label{eq:A3} \\
p(t,z) &= 0, \qquad t \ll 0, ~~ z \in (0-,L),  \label{eq:A4}
\end{align}
where $z = 0-$ is the range coordinate of the accessible boundary, just above the sensor at $z =0$, and $L>0$ is the  range of the inaccessible boundary, assumed large enough so it does not affect the wave over the duration of the measurements. 

After the  travel time coordinate transformation \eqref{eq:A1}, we get that the even in time 
wave
\begin{equation}
P(t,T):= p(t,z(T)) + p(-t,z(T)),
\label{eq:A5}
\end{equation}
satisfies
\begin{align}
\partial_t^2 P(t,T) - \zeta(z(T)) \partial_T \left[ \frac{1}{\zeta(z(T))} \partial_T P(t,T) \right] &= 0, \qquad 
t >0 , ~~ T \in (0-,T(L)), 
\label{eq:A6} \\
\partial_T P(t,0-) &= P(t,T(L)) = 0, \qquad t >0, \label{eq:A7} 
\end{align}
with initial conditions
\begin{align}
P(0,T) = \varphi(T) \approx 2 f(T),  \quad \partial_t P(0,T) = 0, \qquad  T \in (0-,T(L)). \label{eq:A8}
\end{align}
The layered medium is modeled by the piecewise constant impedance 
\begin{equation}
\zeta(z(T)) = \zeta_j, \qquad T \in (T_{j-1}, T_j], ~~ T_j := T(z_j), ~~ j = 0, \ldots, \ell +1,
\label{eq:A9}
\end{equation}
whose jumps at range coordinates $z_j$, ordered as   $0-=z_{-1} < 0 < z_0 \ldots < z_{\ell+1} = L$, give the reflection and transmission coefficients \cite[Chapter 3]{fouque2007wave} 
\begin{equation}
\mathfrak{R}_j:= \frac{\zeta_j - \zeta_{j+1}}{\zeta_j + \zeta_{j+1}}, \quad 
\mathfrak{T}_j:= \frac{2 \sqrt{\zeta_j \zeta_{j+1}}}{\zeta_j + \zeta_{j+1}}, \qquad j = 0, \ldots, \ell.
\label{eq:ART}
\end{equation}

In the reference medium with constant impedance $\zeta_o = \zeta(0)$, the wave is given by d'Alembert's solution
\begin{equation}
P_o(t,T) = \frac{1}{2} \big[\varphi(T-t)+\varphi(T+t)\big].
\label{eq:A10}
\end{equation}
It is the sum of a forward and a backward wave, due to the accessible boundary.

If there is a single scattering layer ($\ell = 0$) in the medium, we obtain after a standard 
calculation as described in  \cite[Chapter 3]{fouque2007wave} that
\begin{align}
P(t,T) &= \sum_{q=0}^\infty (-\mathfrak{R}_0)^q \12\left[ \varphi(T-t+2 q T_0) + \varphi(T+t-2 q T_0)\right]= \sum_{q=0}^\infty (-\mathfrak{R}_0)^q P_o(t-2 q T_0,T), \label{eq:A11}
\end{align}
if $T \in (0-,T_0),$ whereas for $T >T_0$ we have
\begin{align}
P(t,T) &= \frac{\mathfrak{T}_0 \sqrt{\zeta_1}}{\sqrt{\zeta_0}} \sum_{q=0}^\infty (-\mathfrak{R}_0)^q \12 \varphi(T-t+2 q T_0) = \frac{\mathfrak{T}_0 \sqrt{\zeta_1}}{\sqrt{\zeta_0}} \sum_{q=0}^\infty (-\mathfrak{R}_0)^qP_o(t-2 q T_0,T).  \label{eq:A12}
\end{align}
In the last equation we used that 
$\varphi(T+t-2 q T_0)=0$ for $T >T_0 \gg t_f$ and time $t >2q T_0$ at  which the $q^{\rm th}$ transmitted 
wave can be  observed. 

The  series over $q$  in equations  \eqref{eq:A11}--\eqref{eq:A12}  account for the multiple reflections at the interface $T=T_0$. We have a train of waves that look just like the wave in the reference medium, with delays $2 q T_0$ corresponding to the number of roundtrips between the accessible boundary and the interface.
Using causality, we conclude  that
\begin{equation}
P(j\tau,T) \in \mathscr{S}_{o,j}:=\mbox{span}\left\{ P_o(j' \tau,T) , ~ j'=0,\ldots, j \right\}, \qquad \mbox{\rm if} ~~ \frac{2 T_0}{\tau} \in \mathbb{N}.
\label{eq:A13}
\end{equation}
Otherwise, $P(j\tau,T) $ is approximated in $\mathscr{S}_{o,j}$ with some error, which is small 
if $\tau$ is small with respect to the scale of variation of $\varphi(t)$ and therefore $f(t)$.

If the medium has multiple layers ($\ell \ge 1$), the expression of $P(t,T)$ is given by a more complicated series, with each term corresponding to a sequence of scattering events \cite[Chapter 3]{fouque2007wave}. Nevertheless, the conclusion is similar to the above: If the travel time 
 between the interfaces is an integer multiple of $\tau$, which  corresponds to a ``Goupillaud medium" \cite[Section 3.5.4]{fouque2007wave}, then the snapshots $P(j \tau,T)$ are represented exactly in the span of the snapshots in the reference medium. Otherwise, we have an error that is small
if $\tau$ is small with respect to the scale of variation of $f(t)$.

In conclusion, in the one-dimensional case, as long as $\tau$ is small enough, the orthonormal snapshots are  approximately the same as in the reference medium, 
in the travel time coordinate. Furthermore, if we have an accurate estimate of the smooth part of the wave speed, 
called $c_o(z)$, we can transform to the range coordinate and obtain \eqref{eq:approxV}.

\subsubsection{Higher dimensions.}
Here we suppose that the waves generated by a source at range $z =0$ propagate in the half space $z > 0-$ filled with a layered medium with wave speed $c(z)$ and impedance $\zeta(z)$. Consider the system of coordinates $\bx = (\bx^\perp,z) $, 
with cross-range $\bx^\perp \in \RR^d$, for $d=1$ or $2$,   and let the source be $f'(t) S(\bx^\perp) \delta(z)$, with cross-range profile $S(\bx^\perp)$. Then, if we Fourier transform the acoustic wave equation for the pressure $p(t,\bx)$ with respect to $t$ and $\bx^\perp$,  we obtain a family of one-dimensional Helmholtz equations 
\begin{equation}
\om^2 \hat p(\om,\bka,z) + \zeta^{\bka} (z) c^{\bka }(z) \partial_z \left[ \frac{c^{\bka}(z)}{\zeta^{\bka}(z)}
\partial_z \hat p(\om,\bka,z)\right] = i \om \hat f^{\bka}(\om) c^{\bka}(0) \delta(z),
\label{eq:A14}
\end{equation}
for the time harmonic plane waves 
\begin{equation}
\hat p(\om,\bka,z):= \int_{\RR} dt \int_{\RR^d} d \bx^\perp \, p(t,\bx^\perp,z) e^{i \om(t-\bka \cdot \bx^\perp)}.
\end{equation}
Here $\bka$ is the slowness vector, with units of time over length, which defines the plane wave speed and impedance
\begin{equation}
c^{\bka}(z):= \frac{c(z)}{\sqrt{1 - c^2(z) |\bka|^2}}, \quad \zeta^{\bka}(z):= \frac{\zeta(z)}{\sqrt{1 - c^2(z) |\bka|^2}},
\label{eq:A16}
\end{equation}
and
\begin{equation}
\hat f^{\bka}(\om):= \frac{\hat f(\om) \hat S(\om \bka)}{\sqrt{1 - c^2(0) |\bka|^2}}.
\label{eq:A17}
\end{equation}
We assume that the source excites propagating waves only i.e., $\bka$ in the support of $S(\om \bka)$ satisfies $|\bka| < \min_z c^{-1}(z)$, so that equations \eqref{eq:A16}-\eqref{eq:A17}   return real values.

Now we can use as in the previous section the travel time transformation 
\begin{equation}
T^{\bka}(z):= \int_0^z \frac{dz'}{c^{\bka}(z')},
\label{eq:A18}
\end{equation}
and obtain that
\begin{equation}
P^{\bka}(t,T^{\bka}) := \int_{\RR} \frac{d \om}{2 \pi} \, \hat p(\om,\bka,z(T^{\bka})) e^{-i \om t} =
\int_{\RR^d} d \bx^\perp \, p(t + \bka \cdot \bx^\perp,\bx^\perp, z(T^{\bka})),
\end{equation}
satisfies an equation like \eqref{eq:A6}, with $\zeta(z)$ replaced by $\zeta^{\bka}(z)$ and
$f(t)$ replaced by $f^{\bka}(t)$. Thus, we can use the results in the previous section to conclude 
that if $\tau$ is small enough, the snapshots of $P^{\bka}(t,T^{\bka})$ can be approximated by those in the reference medium. 
Note however that the wave speed cannot be removed completely via the travel time 
transformation, as in the one-dimensional case, because $c(z)$ appears in the 
expression of the impedance $\zeta^{\bka}(z)$. Thus, knowing the kinematics (the smooth part of $c(z)$) is very important for getting the alignment of the wavefronts in the true  layered medium and the reference medium. 

\subsection{Snapshots in a waveguide}
\label{ap:A2}
Here we return to the wave equation in a medium with constant density, and assume for simplicity a two-dimensional waveguide $\bx = (x^\perp,z) \in (0,D) \times (0-,\infty)$, with sound hard wall 
at $z=0$, representing the accessible boundary $\partial \Omega_{\rm ac}$, and sound soft side walls at $x^\perp \in \{0,D\}$ that are part of the inaccessible boundary $\partial \Omega_{\rm inac}$. The 
remaining part of $\partial \Omega_{\rm inac}$ is an artificial sound soft boundary at $z = L$, for large enough $L$ so that the waves do not reach it over the duration of the measurements. Note that this is the setup for the numerical simulations in section \ref{sect:num1}.

The waveguide is filled with a homogeneous medium with wave speed $\bar c_o$, and contains a thin reflector localized for simplicity at the range $z = z_0$, modeled by the reflectivity $\mathfrak{r}(x^\perp) \delta_{z_0}(z)$ as follows
\begin{equation}
\frac{1}{c^2(\bx)} = \frac{1}{\bar c_o^2} \left[ 1 + \mathfrak{r}(x^\perp) \delta_{z_0}(z)\right].
\label{eq:W1}
\end{equation}
We analyze  the  acoustic pressure $p^\ss(t,\bx)$  in the waveguide, related 
to the wave $w^\ss(t,\bx)$ as explained in section \ref{sect:setup}. The excitation is as in \eqref{eq:F1},
and the pulse $f(t)$ is given by an even envelope function $F$ supported in the interval $(-1,1)$ and modulated at the central frequency $\om_c$,
\begin{equation}
f(t):=F\Big(\frac{t}{t_f}\Big) \cos(\om_c t).
\label{eq:W2}
\end{equation}
The snapshots at $z \ne z_0$ are defined by the even extension in time of the pressure, divided by the constant speed $\bar c_o$,
\begin{equation}
u^\ss(t,\bx) := \big[p^\ss(t,\bx) + p^\ss(-t,\bx)\big]/\bar c_o.
\end{equation}

The analysis uses the mode decomposition of $u^\ss(t,\bx)$, based on its expansion in the $L^2(0,D)$ orthonormal basis $\{\psi_j(x^\perp), ~j \ge 1\}$, where 
\begin{equation}
\psi_j(x^\perp) =\sqrt{ \frac{2}{D}} \sin(\alpha_j x^\perp), \quad \alpha_j:= \frac{\pi j}{D}, \qquad j \ge 1,
\label{eq:W3}
\end{equation}
are the eigenfunctions of the operator $\partial_{x^\perp}^2$ acting on functions of $x^\perp \in (0,D)$,  with homogeneous Dirichlet boundary conditions. We are interested in the propagating modes,
indexed by $j = 1, \ldots, N = \lfloor {k_c D}/{\pi}\rfloor$, because the evanescent modes 
generated by the reflectivity at $z = z_0$ are negligible by the time they reach the array. Here 
$k_c = \om_c/\bar c_o$ is the wave number at the central frequency, and we  assume that the bandwidth $B = O(1/t_f)$ of the probing pulse is small enough, so that 
\begin{equation}
\Big\lfloor \frac{\om D}{\pi \bar c_o}\Big\rfloor \approx \Big\lfloor \frac{k_c D}{\pi }\Big\rfloor, \qquad 
\forall \,  \om \in (\om_c-B, \om_c+B).
\end{equation}

The expression of the snapshots in the empty (reference) waveguide is obtained after a standard calculation, as explained for example in \cite[Chapter 20]{fouque2007wave},
\begin{align}
u^\ss_o(t,\bx) \approx \sum_{j=1}^N \psi_j(x^\perp) u_{o,j}^\ss(t,z)  + \mbox{ evanescent}.  \label{eq:W4p}
\end{align}
It is a superposition of one-dimensional propagating waves (modes)
\begin{align}
u_{o,j}^\ss(t,z):=
\frac{2 k_c  \psi_j(x^\perp_s)}{\bar c_o \beta_j(\om_c)}
\left\{ F \Big( \frac{t-z/c_{o,j}}{t_f}\Big) \cos \big[ \beta_j(\om_c) z -\om_c t \Big] 
+F \Big( \frac{t+z/c_{o,j}}{t_f}\Big) \cos \big[ \beta_j(\om_c) z +\om_c t \Big] \right\},  \label{eq:W4}
\end{align}
with wave numbers  
\begin{equation}
\beta_j(\om):= \mbox{sign}(\om) \sqrt{\frac{\om^2}{\bar c_o^2} - \alpha_j^2}, \qquad 
j = 1, \ldots, N,
\label{eq:W5}
\end{equation}
and the approximation in \eqref{eq:W4p} is due to the small bandwidth assumption that allows us to write
\begin{equation}
\beta_j(\om)z \approx \beta_j(\om_c)z + (\om-\om_c)\beta'_j(\om_c) z, \qquad \forall \,  \om \in (\om_c-B, \om_c+B).
\label{eq:W6}
\end{equation}
Again, we see that due to the accessible boundary, we have both forward and backward 
going waves in  \eqref{eq:W4}. The backward waves are observed only at small $z$ and time $t = O(t_f)$. The waveguide is dispersive, so the propagation is at mode dependent group speed 
\begin{equation}
c_{o,j}:= \frac{1}{\beta'_j(\om_c)} = \frac{\bar c_o \beta_j(\om_c)}{k_c},
\label{eq:W7}
\end{equation}
which is different than the phase speed $\om_c/\beta_j(\om_c)$, for $j = 1, \ldots, N$.

The expression of the snapshots in the waveguide with the reflectivity given in \eqref{eq:W1} involves 
a series of multiple scattering events at the reflector. For our purposes it suffices to look at the first two terms in this series, corresponding to the single scattering, Born approximation. The analysis of the higher order terms is  similar and does not bring anything new. A standard calculation that uses approximations like \eqref{eq:W6} gives that the snapshots are 
\begin{equation}
u^\ss(t,\bx) \approx \sum_{j=1}^N \psi_j(x^\perp) \big[u^\ss_{o,j}(t,z) + u^\ss_{{\rm Born},j}(t,z)\big] +O(\mathfrak{r}^2) + \mbox{evanescent},
\label{eq:W8}
\end{equation}
where for $z \in (0-,z_0)$ we have 
\begin{align}
u^\ss_{{\rm Born},j}(t,z) \approx  \frac{k_c^2 }{
2 \bar c_o^2 \beta_j }  \sum_{l=1}^N \mathfrak{r}_{j,l} \frac{\psi_l(x_s)}{
\beta_l} & \partial_t 
\Big\{F \Big[ \frac{t- \beta_j' z-z_0(\beta_j'+\beta_l')}{t_f} \Big]
\cos \big[ \beta_j z - \om_c t + z_0 ( \beta_j + \beta_l) \big] \nonumber \\
&+ F \Big[ \frac{t+ \beta_j' (z-z_0) - z_0\beta_l')}{t_f} \Big]
\cos \big[ \beta_j (z-z_0) + \om_c t - z_0  \beta_l \big]\Big\}
\label{eq:W9}
\end{align}
and for $z > z_0$ we have
\begin{align}
u^\ss_{{\rm Born},j}(t,z) \approx  \frac{k_c^2}{
\bar c_o^2 \beta_j } \sum_{l=1}^N \mathfrak{r}_{j,l} \frac{\psi_l(x_s)}{
 \beta_l}&  \partial_t 
\Big\{F \Big[ \frac{t-\beta_j' z-z_0(\beta_j'+\beta_l') }{t_f} \Big]
\cos \big[ \beta_j z - \om_c t + z_0 ( \beta_j + \beta_l) \big] \nonumber \\
&+ F \Big[ \frac{t- \beta_j'( z-z_0)-z_0\beta_l'}{t_f} \Big]
\cos \big[ \beta_j (z -z_0)- \om_c t + z_0 \beta_l ) \big]\Big\}. \label{eq:W10}
\end{align}
In these equations we simplified the notation by dropping the $\om_c$ arguments of $\beta_j, \beta_l$ and their derivatives, and we introduced the reflectivity matrix 
\begin{equation}
\mathfrak{r}_{j,l} := \int_{0}^D d x^\perp \, \mathfrak{r}(x^\perp)\psi_j(x^\perp) \psi_l(x^\perp).
\label{eq:W11}
\end{equation}

Note that the terms in \eqref{eq:W9} model two kinds of waves: The first kind strikes the reflector as mode $l$, it is converted to mode $j$, travels to the accessible boundary,
it is reflected there and then travels forward. The second kind strikes the reflector as mode $l$, it is converted to mode $j$ and then travels backward. 
Similarly, the first term in \eqref{eq:W10} models the wave that starts as mode $l$, it is converted 
to mode $j$, travels to the accessible boundary where it reflects and then propagates forward.
The second term models the wave that strikes the reflector as mode $l$, it is converted to mode $j$ and then propagates forward.
 We now  show that these waves can be approximated in the span 
of the time delayed reference waveguide modes \eqref{eq:W4}.

Let us introduce the travel times 
\begin{equation}
t_{j,l}:= z_0 ( \beta'_j + \beta'_l)= \frac{z_0}{c_{o,j}} + \frac{z_0}{c_{o,l}},
\label{eq:W12}
\end{equation}
corresponding to the propagation of the envelope of the wave at group speeds 
\eqref{eq:W7}, and 
\begin{equation}
T_{j,l}:= \frac{z_0 (\beta_j+\beta_l)}{\bar c_o k_c},
\label{eq:W13}
\end{equation}
corresponding to the propagation of the phase. Then, expanding the cosine in \eqref{eq:W9} 
we get 
\begin{align}
u^\ss_{{\rm Born},j}(t,z) \approx  &\frac{k_c^2 }{
2 \bar c_o^2 \beta_j }  \sum_{l=1}^N \frac{\psi_l(x_s)}{
\beta_l} \mathfrak{r}_{j,l} \left\{ \cos \big[\om_c(T_{j,l}-t_{j,l})\big] \partial_t 
\Big\{F \Big[ \frac{t- t_{j,l} - z/c_{o,j}}{t_f} \Big]
\cos \big[ \beta_j z - \om_c (t-t_{j,l}) \big]\right.  \nonumber \\
&\hspace{1.2in}+ F \Big[ \frac{t-t_{j,l} + z/c_{o,j})}{t_f} \Big]
\cos \big[ \beta_j z + \om_c (t - t_{j,l}) \big]\Big\}\nonumber \\
&\hspace{1.2in} -\sin \big[\om_c(T_{j,l}-t_{j,l})\big] \partial_t 
\Big\{F \Big[ \frac{t- t_{j,l} - z/c_{o,j}}{t_f} \Big]
\sin \big[ \beta_j z - \om_c (t-t_{j,l}) \big] \nonumber \\
&\hspace{1.2in}\left. - F \Big[ \frac{t-t_{j,l} + z/c_{o,j})}{t_f} \Big]
\sin \big[ \beta_j z + \om_c (t - t_{j,l}) \big]\Big\}\right\}
\label{eq:W14}
\end{align}
for $z \in (0-,z_0)$. Recalling equation \eqref{eq:W4}, we note that the first curly bracket is proportional to 
$\partial_t u_{o,l}^\ss(t-t_{j,l},z)$. The second curly bracket  is approximately proportional to $\partial_t^2 u_{o,l}^\ss(t-t_{j,l},z)$, because 
\begin{align*}
\partial_t \left\{ F \Big[ \frac{t- t_{j,l} - z/c_{o,j}}{t_f} \Big]
\cos \big[ \beta_j z - \om_c (t-t_{j,l}) \big]\right\} = \om_c F \Big[ \frac{t- t_{j,l} - z/c_{o,j}}{t_f} \Big]
\sin \big[ \beta_j z - \om_c (t-t_{j,l}) \big]\left[ 1 + O \left(\frac{1}{\om_c t_f}\right)\right] \\
\partial_t \left\{ F \Big[ \frac{t- t_{j,l} + z/c_{o,j}}{t_f} \Big]
\cos \big[ \beta_j z + \om_c (t-t_{j,l}) \big]\right\} = -\om_c F \Big[ \frac{t- t_{j,l} + z/c_{o,j}}{t_f} \Big]
\sin \big[ \beta_j z + \om_c (t-t_{j,l}) \big]\left[ 1 + O \left(\frac{1}{\om_c t_f}\right)\right]
\end{align*} 
and we have assumed 
\[
\frac{1}{\om_c t_f} = O \left(\frac{B}{\om_c}\right) \ll 1.
\]
For a small enough time sample interval $\tau$, the time derivatives of $u_{o,l}^\ss(t-t_{j,l},z)$ can be approximated with finite differences, so we conclude that the snapshots \eqref{eq:W9} evaluated at 
range $z \in (0-,z_0)$ can be approximated by linear combinations of the time delayed reference snapshots \eqref{eq:W4}. Similarly, it follows that the result also holds  for the snapshots \eqref{eq:W10} evaluated at $z > z_0$. 

In the ROM construction we do not use  the mode decomposition. However, if the sensors are closely 
spaced in the array, so that we can approximate the sum over them by an integral over the 
array aperture $\cA \subseteq (0,D)$, we can write
\begin{equation}
\int_{\cA} d x^{\perp} u_o^\ss(t,\bx) \psi_j(x^\perp) \approx \sum_{l=1}^N Q_{j,l} u_{o,l}^\ss(t,z) ,
\end{equation}
where $\bQ = \big(Q_{j,l}\big)_{j,l = 1, \ldots, N}$ is the mode coupling matrix 
\begin{equation}
Q_{j,l} := \int_{\cA} d x^\perp \, \psi_j(x^\perp) \psi_l(x^\perp), \qquad j,l = 1, \ldots, N.
\end{equation}
If $\cA$ is large enough, $\bQ$ is invertible, so the snapshots at the array carry the same information 
as the modes \eqref{eq:W4}. This is what we need for our approximation, in addition 
to the small $\tau$ required to deal with the discrete time samples of the wave, as in the previous section.

%========================================================================
\section*{References}
\bibliography{biblio} \bibliographystyle{siam}
\end{document}